\documentclass[12pt,a4paper,final]{article}

\RequirePackage[backend=biber,style=authoryear,natbib,
  uniquename=false]{biblatex}
\RequirePackage{csquotes}
\usepackage{amsmath}    % AMS math style
\usepackage{amsfonts}   % Include AMS fonts
\usepackage{amssymb}    % Include AMS symbols
\usepackage{amsthm}     % For theorem-like environments
\usepackage{mathrsfs}   % Formal Script Math Symbol font
\usepackage{parskip}    % For whitespaced paragraphs
\usepackage{calc}       % For use of \widthof in def of \convas macro
\usepackage{enumitem}   % For nicer spacing in itemized environments
\usepackage{graphicx}   % For \includegraphics
\usepackage{lpic}       % For latex-annotated praphics
\usepackage{tikz-cd}    % For commutative diagrams
  \usetikzlibrary{matrix,arrows,decorations.pathmorphing}
\usepackage[english]{babel}
\RequirePackage{hyperref}
\hypersetup{
  colorlinks=true, %set true if you want colored links
  linktoc=all,     %set to all for both sections and subsections linked
  linkcolor=black,  %choose some color if you want links to stand out
  citecolor=black,
  urlcolor=black,
  }
\usepackage{cleveref}   % For more versatile referencing
\usepackage{multicol}   % To display TOC in two columns
\usepackage{mathptmx}   % For Times Roman font
\usepackage{titlesec}   % For nicer typesetting of chapter, section, etc headers
\titleformat{\chapter}[display]
  {\normalfont\rmfamily\huge\bfseries}
  {\chaptertitlename\ \thechapter}{20pt}{\Huge}
\titleformat{\section}
  {\normalfont\rmfamily\Large\bfseries}
  {\thesection}{1em}{}
\titleformat{\subsection}
  {\normalfont\rmfamily\large\bfseries}
  {\thesubsection}{1em}{}

\renewcommand{\baselinestretch}{1.00}

\newlength{\vslength}
\newcommand{\ie}{{\it i.e.}}
\newcommand{\cf}{{\it cf.}}
\newcommand{\eg}{{\it e.g.}}

\newcommand{\etc}{{\it etcetera}}
\newcommand{\vv}{{\it vice versa}}

\newcommand{\QQ}{{\mathbb Q}}
\newcommand{\RR}{{\mathbb R}}

\newcommand{\NN}{{\mathbb N}}

\newcommand{\scrA}{{\mathscr A}}
\newcommand{\scrB}{{\mathscr B}}
\newcommand{\scrC}{{\mathscr C}}

\newcommand{\scrE}{{\mathscr E}}
\newcommand{\scrF}{{\mathscr F}}

\newcommand{\scrP}{{\mathscr P}}

\newcommand{\scrS}{{\mathscr S}}
\newcommand{\scrT}{{\mathscr T}}
\newcommand{\scrU}{{\mathscr U}}

\newcommand{\scrX}{{\mathscr X}}
\newcommand{\scrY}{{\mathscr Y}}

\newcommand{\al}{{\alpha}}
\newcommand{\be}{{\beta}}

\newcommand{\ep}{{\epsilon}}
\newcommand{\vep}{{\varepsilon}}
\renewcommand{\emptyset}{{\varnothing}}

\newcommand{\ft}[2]{{\textstyle{\frac{#1}{#2}}}}
\newcommand{\supp}{{\mathrm{supp}}}
\newcommand{\dd}{{\,{\mathrm d}}}

\newcommand{\GTI}[4]% Bourbaki, General Topology I, {ch}, {sec}, {no}, {...}
    {\citep{BourGTI}, Ch.\,{#1}, \S\,{#2}, No.\,{#3}, {#4}}
\newcommand{\GTIthree}[3]% Bourbaki, General Topology I, {ch}, {sec}, {no}
    {\citep{BourGTI}, Ch.\,{#1}, \S\,{#2}, No.\,{#3}}
\newcommand{\GTII}[4]% Bourbaki, General Topology II, {ch}, {sec}, {no}, {...}
    {\citep{BourGTII}, Ch.\,{#1}, \S\,{#2}, No.\,{#3}, {#4}}
\newcommand{\GTIIthree}[3]% Bourbaki, General Topology II, {ch}, {sec}, {no}
    {\citep{BourGTII}, Ch.\,{#1}, \S\,{#2}, No.\,{#3}}
\newcommand{\TVS}[4]% Bourbaki, Topological vector spaces, {ch}, {sec}, {no}, {...}
    {\citep{BourTVS}, Ch.\,{#1}, \S\,{#2}, No.\,{#3}, {#4}}
\newcommand{\TVSthree}[3]% Bourbaki, Topo vector spaces, {ch}, {sec}, {no}
    {\citep{BourTVS}, Ch.\,{#1}, \S\,{#2}, No.\,{#3}}
\newcommand{\II}[4]% Bourbaki, Integration I, {ch}, {sec}, {no}, {...}
    {\citep{BourII}, Ch.\,{#1}, \S\,{#2}, No.\,{#3}, {#4}}
\newcommand{\IIthree}[3]% Bourbaki, Integration I, {ch}, {sec}, {no}
    {\citep{BourII}, Ch.\,{#1}, \S\,{#2}, No.\,{#3}}
\newcommand{\III}[4]% Bourbaki, Integration II, {ch}, {sec}, {no}, {...}
    {\citep{BourIII}, Ch.\,{#1}, \S\,{#2}, No.\,{#3}, {#4}}
\newcommand{\IIIthree}[3]% Bourbaki, Integration II, {ch}, {sec}, {no}
    {\citep{BourIII}, Ch.\,{#1}, \S\,{#2}, No.\,{#3}}

\newcommand{\comment}[1]{{}}
\renewcommand{\qedsymbol}{$\Box$}

\newtheoremstyle{customtheorem}% name of the style to be used
  {0.5em}% measure of space to leave above the theorem. E.g.: 3pt
  {0.2em}% measure of space to leave below the theorem. E.g.: 3pt
  {\itshape}% name of font to use in the body of the theorem
  {}% measure of space to indent
  {\scshape}% name of head font
  {}% punctuation between head and body
  {1ex}% space after theorem head; " " = normal interword space
  {}% Manually specify head

\theoremstyle{customtheorem}
\newtheorem{theorem}{Theorem}[section]
\newtheorem{lemma}[theorem]{Lemma}
\newtheorem{proposition}[theorem]{Proposition}
\newtheorem{corollary}[theorem]{Corollary}
\newtheorem{definition}[theorem]{Definition}

\newtheoremstyle{customremark}% name of the style to be used
  {0.5em}% measure of space to leave above the theorem. E.g.: 3pt
  {0.2em}% measure of space to leave below the theorem. E.g.: 3pt
  {}% name of font to use in the body of the theorem
  {}% measure of space to indent
  {\scshape}% name of head font
  {}% punctuation between head and body
  {1ex}% space after theorem head; " " = normal interword space
  {}% Manually specify head

\theoremstyle{customremark}
\renewenvironment{proof}{\par\noindent{\scshape Proof}\;}{\hfill\qedsymbol\par}
\newtheorem{remark}[theorem]{Remark}
  \AtBeginEnvironment{remark}{%
    \pushQED{\qed}%\renewcommand{\qedsymbol}{$\triangle$}%
  }
  \AtEndEnvironment{remark}{\popQED\endremark}
\newtheorem{example}[theorem]{Example}
  \AtBeginEnvironment{example}{%
    \pushQED{\qed}%\renewcommand{\qedsymbol}{$\triangle$}%
  }
  \AtEndEnvironment{example}{\popQED\endexample}

\AtEndDocument{
  \printbibliography
}
\begin{document}
\thispagestyle{empty}
\renewcommand{\baselinestretch}{0.80}
\title
{%\vspace*{-26mm}
  % Existence and phase structure of random histogram limits
  % Existence and phase structure of\\[-1mm] random projective limit measures
  Existence and phase structure of\\[2mm] random inverse limit measures
}
%
%\vspace*{-10mm}}
%
\author
% { B.~J.~K. {Kleijn}${}^{1}$ and H. {de~With}${}^{2}$\\[1mm]
%   {\small\it ${}^{1}$ Korteweg-de~Vries Institute for Mathematics,
%     Univ. of Amsterdam\thanks{corresponding author, email: b.j.k.kleijn@uva.nl}}\\
%   {\small\it ${}^{2}$ Affiliation ABN Bank?}  
% }
{
  B.~J.~K. {Kleijn}\\[1mm]
  {\small\it Korteweg-de~Vries Institute for Mathematics,
    Univ. of Amsterdam}\\
  {\footnotesize %corresponding author, B.~J.~K. {Kleijn},
  email: {\tt b.j.k.kleijn@uva.nl}, orcid: {0000-0003-2495-2922}}
}
\date{May 2025}
\maketitle
\renewcommand{\baselinestretch}{1.00}
\begin{abstract}\noindent
Analogous to Kolmogorov's theorem for the existence of stochastic
processes describing random functions, we consider theorems for
the existence of stochastic processes describing random measures,
as limits of inverse measure systems. Specifically, given a
coherent inverse system of random (bounded/signed/positive/probability)
histograms on refining partitions, we study conditions for the
existence and uniqueness of a corresponding random inverse limit,
a Radon probability measure on the space of
(bounded/signed/positive/probability) measures.
Depending on the topology (vague/tight/weak/total-variational)
and Kingman's notion of complete randomness, the limiting random
measure is in one of four phases, distinguished by their
degrees of concentration (support/domination/discreteness). Results
are applied in the well-known Dirichlet and Polya tree families of random
probability measures and in a new Gaussian family of signed inverse 
limit measures. In these three families, examples of all four phases
occur and we describe the corresponding conditions on defining
parameters.
\end{abstract}
%
%%%%%%%%%%%%%%%%%%%%%%%%%%%%%%%%%%%%%%%%%%%%%%%%%%%%%%%%%%%%%%%%%%%%%%%%%%%%%%%
%
% Table of contents
%
%%%%%%%%%%%%%%%%%%%%%%%%%%%%%%%%%%%%%%%%%%%%%%%%%%%%%%%%%%%%%%%%%%%%%%%%%%%%%%%
%
\vspace*{3mm}
\begin{multicols}{2}
\setcounter{tocdepth}{2}
{\scriptsize \tableofcontents}
\end{multicols}
\newpage
%
%%%%%%%%%%%%%%%%%%%%%%%%%%%%%%%%%%%%%%%%%%%%%%%%%%%%%%%%%%%%%%%%%%%%%%%%%%%%%%%
%
% Main text
%
%%%%%%%%%%%%%%%%%%%%%%%%%%%%%%%%%%%%%%%%%%%%%%%%%%%%%%%%%%%%%%%%%%%%%%%%%%%%%%%
\section{Introduction}
\label{sec:intro}

In many parts of science probability distributions
on function spaces play a central role: random functions of a
time-variable, for example, describe \emph{evolution with random
influences} of systems in all fields, supported by
the mathematical theory of \emph{stochastic processes}. The
central example, Brownian motion, serves as a primitive for
the modelling of a wide range of phenomena. Conceptual developments
(like the introduction of martingales and the Markov property)
deepen the theory and extend modelling applicability greatly.

Underpinning theories of random functions, is \emph{Kolmogorov's
theorem for the existence of stochastic processes}: given a domain
$\scrX$ for real valued functions $f$, we depart from
the collection of projections $f\mapsto f_S=(f(s):s\in S)$ (where
$S$ is any finite subset of $\scrX$), and for every $S$ we
provide a probability distribution $\Pi_S$ for the projection:
\[
  f_S:=\bigl(f(s):s\in S\bigr)\sim\Pi_S.
\]
Consistency among projections dictates the necessary condition
that $S'\subset S$ implies that $\Pi_{S'}$ is marginal to $\Pi_S$.
Kolmogorov's theorem says that this consistency condition is also
sufficient: if the chosen $\Pi_S$ are consistent, there exists
a probability distribution $\Pi$ for a random function $f$ such
that the projected random $f_S$ are distributed according to
$\Pi_S$, for all finite $S\subset\scrX$.

Kolmogorov's theorem has advantages and disadvantages: on one
hand, mere consistency is enough without further conditions,
making existence remarkably easy to prove. This enhances
applicability greatly, as does the simplicity of the approach:
properties of $\Pi$ derive from those of the (user-defined)
$\Pi_S$, and calculations for random $f_S$ are feasible
because they take place in \emph{finite-dimensional}
probability spaces. On the other hand, the limiting $\Pi$ is
a Borel measure for the product topology on $\RR^\scrX$, which
dissociates $\Pi$ from some interesting finer (\eg\ metrizable)
topologies, while the indistinct nature of $\RR^\scrX$ forces
the imposition of extra conditions on the choices of the
$\Pi_S$ to induce properties like ($\Pi$-almost-sure)
continuity, differentiability or integrability.

Probability measures on spaces of measures are (less common,
but) also of interest in various parts of science. For
example, in non-parametric statistics (\cite{LeCam86}),
particularly of the Bayesian type (\cite{Ghosal17,Kleijn21}),
probability distributions on spaces of probability measures play
a central role; machine learning shares much with statistics,
and random (probability) measures also feature prominently
there. 
%\eg\ with \emph{Latent Dirichlet Allocation} (\cite{Blei03}),
% (applied in natural language processing,
% evolutionary genetics, and many other settings)
%to mention just one of many examples.

In those disciplines the most straightforward approach
is often to define random measures by mapping of random
functions: it is commonplace, for example, to normalize an
integrable positive random function, in order to define
a random probability density function. But it remains
attractive to think about constructions of random measures
that take place directly on the space of measures in its
full generality.
In this sense, based on the Poisson-type family of
\emph{completely random measures} (\cite{Kingman67,Kingman75}),
a well-developed theory of point-processes (\eg\ the family
of Dirichlet random probability measures) exists
(for an overview, see \cite{daley03,daley07}) and is
applied widely. Further examples exist, but a 
comprehensive mathematical theory for the construction and
study of probability measures on measure spaces is lacking to
date.
% (See, for example, \cite{daley07}, p.42: \emph{``it
% appears to be an open problem to find simple sufficient
% conditions [...] for the realizations of a random measure
% to be absolutely continuous''}.)

Ideally such a theory of random measures would
be based on an existence theorem like Kolmogorov's theorem
for random functions. In this paper, we formulate several
new existence theorems of this kind. For a collection
$\scrA$ of finite, measurable partitions $\al$ of $\RR$,
we choose distributions $\Pi_\al$ and define, for all
restrictions of $\mu$ to $\al$,
\[
  \mu_\al:=\bigl(\mu(A):A\in\al\bigr)\sim\Pi_\al.
\]
These projections $\mu_\al$ are called \emph{random
histograms} in what follows, and their consistency
follows from the additivity of $\mu$: if $\be$ refines
$\al$, then for any $A\in\al$,
\[
  \Bigl(\sum_{B\subset A}\mu_\be(B):A\in\al\Bigr)\sim\Pi_\al.
\]
when $\mu_\be\sim\Pi_\be$. Under which conditions does
such a \emph{a random histogram system} have a (unique)
\emph{histogram limit} $\Pi$? (by which we mean a
probability measure $\Pi$ on the space of measures
with $\al$-restrictions that match the $\Pi_\al$).

This question is, of course, not new: both Bayesian non-parametric
statistics and stochastic analysis have formulated a wide
variety of conditions for existence, more or less
independently. First explorations of the subject in
stochastic analysis date back to (\cite{Bochner55}) and
(\cite{Choksi58}), who formulate the classical
Bochner-Kolmogorov conditions for existence of a random
distribution function on $\RR$. Other approaches based on
inner regularity are considered in
(\cite{Metivier63,Schwartz73,BourIII}) and discussed
comprehensively in (\cite{Rao81,Bogachev07}).
Definitions of measure-theoretic inverse-limits
are presented in (\cite{Mallory71,Rao71,Pinter10,Beznea14}).
Limits of random histogram systems in the Bayesian
non-parametric literature are first discussed in
(\cite{Kraft64}), which introduces the P\'olya-tree
family of histogram systems. Many further developments
are based on Kingman's completely random measures,
most prominently in the form of the Dirichlet process
(\cite{Ferguson73,Ferguson74}). For overviews of these and
further developments on non-parametric Bayesian priors,
see (\cite{DeBlasi15,Ghosal17}). Regarding existence, most
noteworthy is (\cite{Orbanz11}), which formulates the
so-called \emph{mean-measure condition} for the existence of a
limit for a system of random probability histograms
$P_\al\sim\Pi_\al$: Orbanz requires that there exists a
Borel probability measure $G$ on $\scrX$ with histogram
projections $G_\al$ that match histogram expectations:
$G_\al(A) = \int P_\al(A)\,d\Pi_\al(P_\al)$ for all
partitions $\al$ and all $A\in\al$.

In this paper we prove and apply several existence
theorems for limits of random histogram systems, and we
analyse the variety of ways in which the corresponding
random measures manifest, in theory and examples. After
introductory remarks and a discussion of the
Bourbaki-Prokhorov-Schwartz theorem in
section~\ref{sec:histograms}, we consider
spaces of probability measures $M^1(\scrX)$
with the tight, weak or total-variational topology,
and we derive conditions that guarantee the existence
and uniqueness of a limiting Radon probability measure $\Pi$
in those cases, in sections~\ref{sec:weakhistogramlimits}
and~\ref{sec:prokhorov}. As it turns out, the manifestations
of their respective random probability measures
are quite different: a limit $\Pi$ that is Radon for
the weak or total-variational topology is
supported on the subset of $M^1(\scrX)$ of \emph{measures
dominated by $G$}, while a $\Pi$ that is Radon
for the tight topology is supported on the
subset of $M^1(\scrX)$ of \emph{measures with support
contained in that of $G$}.

Combined with the Poisson-like manifestation of
completely random measures, in section~\ref{sec:phasestructure}
we distinguish \emph{four phases} for random probability
histogram limits: \emph{absolutely-continuous, fixed-atomic,
continuous-singular and random-atomic}. The results are
applied to known examples like the Dirichlet (in
section~\ref{sec:dirichletlimits}) and P\'olya-tree
families (in section~\ref{sec:polya}) and we re-derive
and sharpen some of the existing results.

In section~\ref{sec:gauss}, we consider spaces of signed
measures with the vague and weak topologies and
derive conditions that guarantee the existence and uniqueness of
Radon histogram limits $\Pi$. The generalization to signed
measures accommodates a new family of Gaussian
probability distributions, defined as limits of random
histogram systems of the form:
\[
  \bigl(\Phi_\al(A):A\in\al\bigr)
    \sim N\bigl( \lambda_\al, \Sigma_\al \bigr),
\]
where $N$ denotes a multivariate normal distribution
with (suitably defined) expectation $\lambda_\al$ and
covariance $\Sigma_\al$. All four phases found for
probability histogram limits are also realized in 
this setting, so Gaussian histogram limits exist with
the same wide range of diffuse and point-like manifestations.
We argue that Gaussian histogram limits based on Green's
functions for the harmonic operator generalize the
well-known two-dimensional \emph{Gaussian Free Field}
(\cite{Werner20}) to higher dimensions, suggesting a potential
role in four-dimensional Euclidean quantum field theory.

To conclude we emphasize the constructive nature of the
existence theorems provided: random histogram systems
\emph{not only define but also approximate} random
measures. The approximative property has two large
advantages, one computational and one analytic: firstly
histogram systems consist of \emph{finite-dimensional}
probability distributions, which we can simulate: the
Dirichlet process, for example, derives much of its
immense popularity from its ease of numerical implementation
and use, and this considerable advantage extends to all
histogram methods. The second advantage lies in
mathematical accessibility: the analyses of example
histogram limits in
sections~\ref{sec:dirichletlimits},~\ref{sec:polya}
and~\ref{sec:gauss} are possible, \emph{only}
because calculations with finite-dimensional random
histograms are feasible, and limits of the results
correspond to properties of the infinite-dimensional
histogram limits. 

\subsubsection*{Acknowledgement} The author thanks Jan~van~Mill,
Harm~de~With and Georg~Meyl for numerous insightful discussions
at various stages in the development of this work.

%%%%%%%%%%%%%%%%%%%%%%%%%%%%%%%%%%%%%%%%%%%%%%%%%%%%%%%%%%%%%%%%%%%%%%%%%%%%%%%
\section{Limits of random histogram systems}
\label{sec:histograms}

The existence theorems that follow in sections~\ref{sec:weakhistogramlimits}
and~\ref{sec:prokhorov} require some (tedious but necessary) bookkeeping
of partitions and corresponding notation (subsection~\ref{sub:prelim}),
as well as a preparatory discussion of the relevant existence theorem,
the Bourbaki-Prokhorov-Schwartz theorem (subsection~\ref{sub:BPS}).

\subsection{Inverse systems of random histograms}
\label{sub:prelim}

We start by introducing directed sets of partitions and the associated
histogram systems.

\subsubsection{Measures, partitions and histograms}

Let $\scrX$ be a Hausdorff topological space (further specification,
\eg\ Polishness, compactness, \etc, follows below). For most purposes,
$\scrX$ can to be thought of
as a space like $\RR^d$: for statisticians, this space plays the
role of the \emph{sample space}, while for probabilists and
physicists $\scrX$ represents \emph{Euclidean space-time}. The
space $\scrX$ has a Borel $\sigma$-algebra we denote by $\scrB$.
We consider a collection $\scrA$ of partitions $\al$ of $\scrX$,
consisting of finite numbers of non-empty Borel sets. (The maximal
such collection, denoted $\scrA_0$, contains \emph{all} finite
partitions into non-empty Borel sets.)
We order $\scrA$ partially by refinement of partitions
(if $\al,\be\in\scrA$ and $\be $ refines $\al$, write
$\al\leq\be$) and we assume throughout that $\scrA$ forms a
directed set for ordering by refinement. Naturally, $\al\leq\be$
implies inclusion for the generated $\sigma$-algebras
$\sigma(\al)\subset\sigma(\be)\subset\scrB$. With the
notation ${|\al|}$ for the cardinality of $\al$, let $I(\al)$ denote
the index set $\{1,\ldots,{|\al|}\}$. Furthermore, we associate to
$\al$ a finite, discrete space $\scrX_\al=\{e_1,\ldots,e_{{|\al|}}\}$ and a
mapping $\varphi_\al:\scrX\to\scrX_\al$, such that $\varphi_\al(x)=e_i$
for all $x\in A_i$. For all $\al\leq\be$, we also define
$\varphi_{\al\be}:\scrX_\be\to\scrX_\al$ such that $\varphi_\al=
\varphi_{\al\be}\circ\varphi_\be$ (and $\varphi_{\al\al}$ as the
identity on $\scrX_\al$ for all $\al$), and define $J_{\al\be}(i)\subset
I(\be)$ to be such that $A_i=\cup_{j\in J_{\al\be}(i)}B_j$ for all
$i\in I(\al)$.

Consider $\scrX$ (or any of the discrete spaces $\scrX_\al$): define
$C_b(\scrX)$ (or $C(\scrX_\al)$) to be the linear space of all bounded,
continuous $f:\scrX\to\RR$ (or $f_\al:\scrX_\al\to\RR$), and let
$M_b(\scrX)$ (or $M(\scrX_\al)$) denote the space of all bounded,
signed Radon measures
$\mu$ on $\scrX$ (or $\mu_\al$ on $\scrX_\al$). For $\mu,\nu\in
M_b(\scrX)$, we say that \emph{$\mu$ dominates $\nu$} (notation
$\nu\ll\mu$), if $\mu(B)=0$ implies $\nu(B)=0$, for all $B\in\scrB$.
Define the bilinear form $\langle \mu,f\rangle = \int f\dd\mu$ (or
$\langle \mu_\al,f_\al\rangle_\al=\int f_\al\dd\mu_\al$). For any
$\mu\in M_b(\scrX)$, let $\mu_+,\mu_-$ denote the (unique pair of)
positive measures such that $\mu=\mu_+ - \mu_-$, and define
$|\mu|=\mu_+ + \mu_-$, noting that the total variational norm
$\|\mu\|$ equals $\langle|\mu|,1\rangle$. The bilinear form
$\langle \cdot,\cdot\rangle$ places
$C_b(\scrX)$ in dual correspondence with $M_b(\scrX)$
(or $C(\scrX_\al)$ with $M(\scrX_\al)$). We refer to the resulting
topology on $M_b(\scrX)$ as the \emph{tight topology}, denoted
$\scrT_C$. Let
$M_+(\scrX)$ denote the positive cone in $M_b(\scrX)$ and let
$M^1(\scrX)\subset M_+(\scrX)$ denote the space of all Radon probability
measures on $\scrX$ (or equivalent with $\scrX_\al$). If $\scrX$ is
a Polish space, then so are $M_+(\scrX)$ and $M^1(\scrX)$ (see
\III{IX}{5}{4}{prop.~10}). Alternatively, we view $C_b(\scrX)$ and
$M_b(\scrX)$ as normed spaces, with the supremum norm
$\|\cdot\|_{\infty,\scrX}$ (or $\|\cdot\|_{\infty,\al}$)
and total-variational norm $\|\cdot\|_{1,\scrX}$ (or
$\|\cdot\|_{1,\al}$)) respectively. We refer to the corresponding
norm topology on $M_b(\scrX)$ as the \emph{total-variational topology},
denoted $\scrT_{TV}$.
Below we also consider $M_b(\scrX)$ and $M^1(\scrX)$ in duality with
the space of all bounded Borel-measurable $h:\scrX\to\RR$, based on the
same bilinear form $\langle \mu,h\rangle$. We refer to
the corresponding topology on $M_b(\scrX)$ as the \emph{weak
topology}, denoted $\scrT_1$. Clearly the total-variational topology
refines the weak topology and the weak topology refines
the tight topology.

We specialize to probability measures in
most of this section and exclusively in
sections~\ref{sec:weakhistogramlimits}--\ref{sec:polya},
and generalize to signed measures only in \cref{sec:gauss}.
Summarizing the most basic requirements for $\scrX$, $\scrA$
and $M^1(\scrX)$, we give the following definition.
\begin{definition}
\label{def:minimalconditions}
We say that $\scrX$, $\scrA$ and $M^1(\scrX)$ satisfy
the \emph{minimal conditions}, if,
\begin{itemize}
\item[(i.)] $\scrX$ and $M^1(\scrX)$ are Hausdorff topological spaces;
\item[(ii.)] $\scrA$ is a directed set of finite partitions of $\scrX$
  in terms of non-empty, Borel measurable sets; and,
\item[(iii.)] for any $\al\in\scrA$ and all $A\in\al$,
  $M^1(\scrX)\to[0,1]:P\mapsto P(A)$ is Borel measurable.
\end{itemize}
\end{definition}
% [BOOK ONLY]
% To illustrate the above, consider the following example which
% plays a central role in section~\ref{sec:prokhorov}.
% \begin{example}
% \label{eq:AtoPAismsb}
% Let $\scrX$ be a Polish space and consider $M^1(\scrX)$ as
% a topological space with the tight topology (or finer).
% Then $\scrX$ is a Hausdorff completely regular space, so for any
% open $G\subset\scrX$, the indicator function
% $1_G:\scrX\to\RR$ equals the upper envelope function for
% the family of all continuous functions it dominates
% (\cf\ \GTII{XI}{1}{6}{prop.~5}):
% $1_G(x)=\sup\{f(x):f\in C_b(\scrX),f\leq1_G\}$ for all
% $x\in\scrX$. By separability of $\scrX$, there exists a
% monotone increasing sequence $(f_n)$ in $C_b(\scrX)$ with
% $1_G$ as its pointwise limit.
% Since, for any $f\in C_b(\scrX)$, $P\mapsto \int f\dd P$ is
% tightly continuous, $P\mapsto P(G)=\int 1_G\dd P$ is the monotone
% limit of a sequence of Borel measurable mappings and, as
% such, is itself Borel measurable. For any closed subset of
% $F$, the indicator $P\mapsto P(F)$ (equal to $1-P(G)$,
% for some open $G$) is also measurable, and the same holds
% for subsets $C$ in the Borel $\sigma$-algebra on $\scrX$,
% by $\sigma$-additivity.
% \end{example}
With regard to the third requirement it is noted that, if the
space $\scrX$ is a Polish space, the mappings $P\mapsto P(G)$,
($G\in\scrB$), are measurable with respect to the Borel
$\sigma$-algebra for the tight topology on $M^1(\scrX)$.

For any Borel probability measure $P\in M^1(\scrX)$ on $\scrX$,
there exists a mapping on $\scrA$, $\al\mapsto
P_\al = \bigl(P(A_1),\ldots,P(A_{{|\al|}})\bigr)$,
that takes a finite, measurable partition $\al$ of
$\scrX$ into the \emph{($\al$-)histogram} associated to $P$.
Note that $P(A_1)+\ldots+P(A_{{|\al|}})=P(\scrX)=1$, so
any $P_\al\in M^1(\scrX_\al)$ can be represented by an
element of the simplex $S_{{|\al|}}=\{p\in\RR^{{|\al|}}:
\Sigma_ip_i=1\}$ (and we shall interchange these two perspectives
freely in what follows). Consider $\al,\be\in\scrA$ such that
$\be$ refines $\al$.
%, denote $\al=\{A_1,\ldots,A_{{|\al|}}\}$,
%$\be=\{B_1,\ldots,B_{{|\be|}}\}$.
By finite additivity of the measure $P$, we have, for every
$A\in\al$, % $1\leq i\leq {|\al|}$,
\begin{equation}
  \label{eq:histogramadditive}
  P(A)=P\bigl(\,\cup_{B\subset A}B\bigr)
    =\sum\{P(B):B\in\be, B\subset A\},
  %P(A_i)=P\bigl(\,\cup_{j\in J_{\al\be}(i)}B_j\bigr)
  % =\sum_{j\in J_{\al\be}(i)}P(B_j),
\end{equation}
so the histograms $P_\al$ and $P_\be$ are related through
summation of probabilities for components that are unified when
partitions coarsen. Clearly, any probability measure $P$
defines a collection of probability histograms
related through (\ref{eq:histogramadditive}) which, conversely,
are enough to reconstruct $P$ if $\scrA$ is rich enough (as per
Carath\'eodory's extension). To give these observations regarding
histograms formal expression, we make the following definitions.
For every $\al\in\scrA$, there exists a \emph{projection mapping}
$\varphi_{*\,\al}:M^1(\scrX)\to M^1(\scrX_\al)$,
\begin{equation}
  \label{eq:projectionmapping}
  \varphi_{*\,\al}(P)=\bigl(P(A_1),\ldots,P(A_{{|\al|}})\bigr),
\end{equation}
that maps a probability distribution to its $\al$-histogram.
Based on (\ref{eq:histogramadditive}), for all $\al,\be\in\scrA$
such that $\al\leq\be$, there is a \emph{transition mapping}
$\varphi_{*\,\al\be}:M^1(\scrX_\be)\to M^1(\scrX_\al)$,
\begin{equation}
  \label{eq:transitionmapping}
  \varphi_{*\,\al\be}(P_\be) =
    \biggl( \sum_{B\subset A_1} P_\be(B)\,,\,\,\ldots\,\,,
    \sum_{B\subset A_{|\al|}} P_\be(B) \biggr),
\end{equation}
that maps $\be$-histograms to $\al$-histograms. Then
$\varphi_{*\,\al\al}$ is the identity for any $\al\in\scrA$,
and for any $\al\leq\be\leq\gamma$, we have,
\[
  \varphi_{*\,\al\gamma}= \varphi_{*\,\al\be}\circ\varphi_{*\,\be\gamma}.
\]
% because, for all $\al\leq\be\leq\gamma$ and all $i\in I(\al)$,
% $J_{\al\gamma}(i)=\cup\{J_{\be\gamma}(j):{j\in J_{\al\be}(i)}\}$.
and, for all $\al\leq\be$,
\begin{equation}
  \label{eq:coherentfamilyfunctions}
  \varphi_{*\,\al} = \varphi_{*\,\al\be}\circ\varphi_{*\,\be}.
\end{equation}

Together with the fact that $\scrA$ forms a directed set,
the following property implies that $\scrA$ is rich enough
for histograms to fix measures on all of the Borel
$\sigma$-algebra. 
\begin{definition}
\label{def:resolveX}
A set $\scrA$ of partitions of a Hausdorff topological
space $\scrX$ is said to \emph{resolve} $\scrX$, if
the $\sigma$-algebra generated by the union of all sets $A$
in all partitions in $\scrA$, is the Borel $\sigma$-algebra,
\ie\ if $\sigma(\{A:\al\in\scrA, A\in\al\})=\scrB$.
%for some basis $\scrU$ for the topology, and for any
%$U\in\scrU$, there exists an $\al\in\scrA$ such that
%$U$ is a union of partition sets $A\in\al$.
%[or maybe?] and an $A\in\al$ such that $A\subset U$.
\end{definition}
%It is clear that if $\scrA$ resolves $\scrX$, then for
%any $x\in\scrX$ the filter basis of those partition sets
%$A$ (in any $\al\in\scrA$) that contain $x$, converges
%to $x$.
% \begin{definition}
% \label{def:resolve} {(ALTERNATIVE)}\\
% A collection of partitions $\scrA$ of a topological space
% $\scrX$ is said to is said to be \emph{separating}, if for any
% pair $x,y\in\scrX$ such that $x\neq y$, there exists an
% $\al\in\scrA$ such that $\varphi_\al'(x)\neq\varphi_\al'(y)$.
% We say that $\scrA$ \emph{resolve} $\scrX$ if, for every
% open $U\subset\scrX$, there exists a $\al\in\scrA$ and an element
% $A\in\al$, such that $A\subset U$.
% \end{definition}
% \begin{remark} {(ALTERNATIVE)}\\
% If $\scrX$ is Hausdorff, any $\scrA$ that resolves
% $\scrX$ also separates $\scrX$. A collection of partitions
% $\scrA$ that resolves $\scrX$ makes it possible for every
% $x\in\scrX$, to find a filter $\scrF$ that converges to
% $x$ and consists of only sets that contain some $A\in\al$,
% for some $\al\in\scrA$.
% \end{remark}
To formulate necessary conditions below, we also need a construction of
partitions in terms of a topological basis for $\scrX$.
\begin{definition}
\label{def:generatedpartition}
Let $\scrU$ be a topological basis
% \varinindex{topological basis}{basis!topological}
for $\scrX$. We say that a partition $\al$ (or collection of
partitions $\scrA$) is \emph{generated by the basis}
%\varinindex{generated by the basis}%
%{partition!generated by basis}
$\scrU$, if, for (any $\al\in \scrA$ and) any $A\in\al$, $A$ is
the union of a finite number of subsets obtained through a
finite number of intersections of $\scrX$ with $U$ or with
$\scrX\setminus U$, $(U\in\scrU)$.
\end{definition}
\begin{example}
\label{ex:partitionscountablebasis}
In a topological space $\scrX$ with a countable basis
$\scrU$, we may construct a sequence of refining partitions based
on an enumeration of $\scrU$: start with $\al_0=\{\scrX\}$;
for all $n\geq1$, intersect all sets in $\al_{n-1}$ with $U_n$
and $\scrX\setminus U_n$, and then define $\al_n$ to
consist of all non-empty such intersections. The resulting
$\scrA=\{\al_n:n\geq1\}$ is a fully ordered set and $\scrA$
resolves $\scrX$.
\end{example}

\subsubsection{Domination, histogram densities and total variation}

In dominated families of probability measures, convergence of
histogram systems coincides with martingale convergence of
Radon-Nikodym densities (see also, \cite[appendix~A1.6]{daley07}).
Due to the monotony of the relation $\al\mapsto\sigma(\al)$,
$\scrF=\{\sigma(\al):\al\in\scrA\}$ is a directed filtration.
Furthermore, if $\scrA$ resolves $\scrX$ the limit of the
filtration $\scrF$ (which has the union of all $\sigma(\al)$,
$(\al\in\scrA)$, as a generating ring), is equal to the
Borel $\sigma$-algebra $\scrB$ on $\scrX$. 

Let $P,Q\in M^1(\scrX)$ be given and assume that $P\ll Q$, so that
$P$ has a ($Q$-almost-everywhere unique) Radon-Nikodym density
$p:\scrX\to[0,\infty)$ with respect to $Q$. Consider the
$\sigma(\al)$-measurable functions $p_\al:\scrX\to[0,\infty)$,
defined by,
\begin{equation}
  \label{eq:histogramdensity}
  p_\al(x)=\sum_{A\in\al} \biggl(\frac{1}{Q(A)}
    \int_{A}p(y)\dd Q(y) \biggr)\,1_{A}(x),
\end{equation}
for $Q$-almost-all $x\in\scrX$ (particularly, if $Q(A)=0$,
for some $A\in\al$, then the corresponding term proportional
to $1_{A}$ is ($Q$-almost-everywhere equal to $0$ and therefore)
not included in the sum). We may define, for every $\al\in\scrA$,
the $Q$-dominated probability measure $P_{Q,\al}:\scrB\to[0,1]$:
\begin{equation}
  \label{eq:Pprimealpha}
  P_{Q,\al}(B) = \int_{B}p_\al(x)\dd Q(x),
\end{equation}
for all $B\in\scrB$, where it is noted that, for all $A\in
\al$, $P_{Q,\al}(A)=P(A)$ (``$=P_\al(A)$'', in a
slightly abusive but natural notation that we introduce
in \cref{rem:Palphanotation}). 
\begin{lemma}
\label{lem:histogramconvergenceTV}
Let $\scrX$, $\scrA$ and $M^1(\scrX)$
satisfy the minimal conditions and assume that $\scrA$
resolves $\scrX$. Then, for any $P\in M^1(\scrX)$ and any
dominating $Q\in M^1(\scrX)$,
$P_{Q,\al}$ converges to $P$ in total variation.
% [BOOK ONLY]
% \[
%   \lim_{\al\in\scrA}\bigl\| P - P_{Q,\al} \bigr\|_{1,\scrX}=0,
% \]
% \ie\ $P_{Q,\al}$ converges (as a net) to $P$ in total variation.
\end{lemma}
\begin{proof}
The Radon-Nikodym density function $p_\al$ is
$Q$-almost-everywhere equal to the $\sigma(\al)$-measurable
conditional expectation
${\mathrm E}[p|\sigma(\al)]:\scrX\to[0,\infty)$, and as
such, the $p_\al$ form a non-negative, uniformly integrable
Doob martingale relative to the filtration $\scrF$.
Since $\scrA$ resolves $\scrX$, Doob's martingale convergence
guarantees that $\lim_\al p_\al=p$ in $L^1(\scrX,\scrB,Q)$.
The assertion now follows from the fact
that, for $Q$-dominated probability measures $P$ and $P_{Q,\al}$,
the total variational norm of their difference is proportional
to the $L^1(Q)$-norm of the difference between densities:
\begin{equation}
  \label{eq:TVLOneisometry}
  \|P-P_{Q,\al}\|_{1,\scrX}=\ft12\int\bigl|p(x)-p_\al(x)\bigr|\dd Q(x),
\end{equation}
for all $\al\in\scrA$.
\end{proof}
The above martingale convergence of densities has implications for
the total-variational norm that we shall appeal to in
sections~\ref{sec:weakhistogramlimits} and~\ref{sec:gauss}.
\begin{proposition}
\label{prop:alphatotalvariation}
Let $\mu$ be a bounded, signed Borel measure on $\scrX$. The mapping
$\al\mapsto\|\mu_\al\|_{1,\scrX_\al}$ is monotone increasing.
If $\scrA$ resolves $\scrX$, then the
total-variational norm for $\mu$ equals,
\[
  \|\mu\|_{1,\scrX} =
  \sup_{\al\in\scrA}\|\mu_\al\|_{1,\scrX_\al}
  = \sup_{\al\in\scrA} \sum_{A\in\al}|\mu_\al(A)|
\]
\end{proposition}
\begin{proof}
If $\al,\be\in\scrA$, and $\be$ refines $\al$, then,
\[
  \|\mu_\al\|_{1,\scrX_\al} = \sum_{A\in\al}|\mu_\al(A)|
    \leq \sum_{B\in\be}|\mu_\be(B)| = \|\mu_\be\|_{1,\scrX_\be}.
\]
Let a signed measure $\mu:\scrB\to\RR$ and $\ep>0$ be given.
According to the Hahn-Jordan decomposition, there exists a
$A_+\in\scrB$ such that, for any $A\in\scrB$, $A\subset A_+$,
we have $\mu(A)\geq0$; and for any $A\in\scrB$,
$A\subset\scrX\setminus A_+$, we have $\mu(A)<0$. Moreover,
$\|\mu\|_{1,\scrX} = \mu(A_+)-\mu(\scrX\setminus A_+)$.
Since $\scrA$ is directed and the union of all
$\sigma(\al)\subset\scrB$, ($\al\in\scrA$)
generates a generating ring for $\scrB$, there exist
an $\al\in\scrA$ and a $A_{\al,+}\in\sigma(\al)$ with
$\mu( (A_+ \setminus A_{\al,+})\cup(A_{\al,+} \setminus A_+) )
< \ft12\ep$, so that,
\[
  \|\mu\|_{1,\scrX} \leq \mu(A_{\al,+})-\mu(\scrX\setminus A_{\al,+}) + \ep
    \leq \|\mu_\al\|_{1,\scrX_\al} +\ep,
\]
proving the assertion.
\end{proof}
The quantities $\|P-P\wedge LQ\|$, used to control weak
compactness in section~\ref{sec:weakhistogramlimits}, are also
suprema of their histogram versions.
\begin{lemma}
\label{lem:DPGmatch}
Assume that $\scrA$ resolves $\scrX$. For any $P,Q\in M^1(\scrX)$
such that $P\ll Q$ and any $L>0$, we have,
\[
    \| P - P \wedge LQ\|_{1,\scrX}
   = \sup_{\al\in\scrA}
     \| P_\al - P_\al \wedge LQ_\al \|_{1,\scrX_\al}.
\]
\end{lemma}
\begin{proof}
Write the Radon-Nikodym derivative of $P$ with respect to $Q$ as $p=dP/dQ$.
Let $B\in\scrB$ be given. If $Q(B)=0$ then $P(B)=0$ and,
$(P - P\wedge LQ)(B)=0$ for any $L>0$. If $Q(B)>0$,
it follows from the convexity of $x\mapsto(x)_+$ and Jensen's
inequality that,
\[
  \begin{split}
  (P - P\wedge LQ)(B) &= Q(B)
    \Bigl( \frac1{Q(B)}\int 1_B(x) (p(x) - L)_+\,dQ(x) \Bigr)\\
    &\geq \Bigl( \int 1_B(x) (p(x) - L)\,dQ(x) \Bigr)_+
    = P(B) - P(B)\wedge LQ(B)
  \end{split}
\]
which implies that, for any $\al\in\scrA$,
\[
  \| P_\al - (P \wedge LQ)_\al\|_{1,\scrX_\al}
    \geq \| P_\al - P_\al \wedge LQ_\al \|_{1,\scrX_\al},
\]
and that the mapping $\al\mapsto\|P_\al-P_\al\wedge L Q_\al\|$
is monotone increasing.
Based on proposition~\ref{prop:alphatotalvariation}, we then find
\[
  \| P - P \wedge LQ\|_{1,\scrX}
   = \sup_{\al\in\scrA}\|P_\al - (P \wedge LQ)_\al\|_{1,\scrX_\al}
   \geq \sup_{\al\in\scrA}
     \| P_\al - P_\al \wedge LQ_\al \|_{1,\scrX_\al}.
\]
Note that we have, for every $\al\in\scrA$,
\[
  \begin{split}
  \biggl|\,\,
    \bigl\|P - & P\wedge LQ\bigr\|_{1,\scrX}
      -\bigl\|P_{\al}-P_{\al}\wedge L Q_\al)\bigr\|_{1,\scrX_\al}\biggr|
  \leq \bigl\|(P - P\wedge LQ)-(P_{Q,\al}
    -P_{Q,\al}\wedge L Q)\bigr\|_{1,\scrX}\\[0.3em]
  & = \ft12\int| (p(x) - L)_+ - (p_\al(x) - L)_+ |\,dQ(x)\\ 
  & = \ft12\int| (p(x)-L)1_{\{p(x)>L\}}
    - (p_\al(x)-L)1_{\{p_\al(x)>L\}}|\,dQ(x)\\
  & \leq \ft12\int 
      1_{\{p(x)>L,p_\al(x)>L\}} |p(x) - p_\al(x)|\\
  &\qquad \qquad   + 1_{\{p(x)>L,p_\al(x)\leq L\}} |p(x) - L|
      - 1_{\{p(x)\leq L,p_\al(x)>L\}} |p_\al(x) - L|
      \,dQ(x)\\[0.3em]
  & \leq \ft12\int|p(x)- p_\al(x)|\,dQ(x)
  = \| P-P_{Q,\al}\|_{1,\scrX}
  \end{split}
\]
An appeal to lemma~\ref{lem:histogramconvergenceTV} then proves
the assertion.
\end{proof}

\subsubsection{Random histogram systems and coherence}

Regarding \emph{random} elements $P\in M^1(\scrX)$ (\eg\ random
elements of a Bayesian statistical model), we can project $P$
onto its random histograms, as formalized in the following
proposition, which introduces the notion of \emph{coherence}.
\begin{proposition}
\label{thm:randomhistogram}
Let $\scrX$, $\scrA$ and $M^1(\scrX)$ satisfy the minimal
conditions and let $\Pi$ denote a Borel probability
distribution on $M^1(\scrX)$ describing a random element $P$.
Then for every $\al\in\scrA_0$,
\begin{equation}
  \label{eq:randomhistogram}
  P_\al:=\varphi_{*\,\al}(P)
  =\bigl(P(A_1),\ldots,P(A_{{|\al|}})\bigr)\sim\Pi_\al,
\end{equation}
induces a random histogram $P_\al$ with probability distribution
$\Pi_\al$ on $M^1(\scrX_\al)$. If
$\al\leq\be$, then $P_\al$ and $P_\be$ are \emph{coherent}, \ie,
the distribution of $P_\al$ follows from that of $P_\be$
through summation as in \cref{eq:histogramadditive}.
\end{proposition}
\begin{proof} By  assumption, for every $\al\in\scrA$
and every $A\in\al$, $M^1(\scrX)\to[0,1]:P\mapsto P(A)$
is Borel measurable. Accordingly,
$\Pi_\al=\Pi\circ\varphi_{*\,\al}^{-1}$ is a
Borel probability distribution on $M^1(\scrX_\al)$. Coherence
(\cref{eq:histogramadditive}) is a consequence of
\cref{eq:coherentfamilyfunctions}.
\end{proof}

Our main question may be paraphrased as the converse of the
above proposition: suppose that we provide distributions
$\Pi_\al$ for random histograms $P_\al$, for all
$\al\in\scrA$. {\it Under which conditions does a
collection of (probability) histogram distributions define
a random (probability) measure (uniquely)?} According
to proposition~\ref{thm:randomhistogram}, coherence is
necessary.
\begin{definition}
\label{def:inverselimits}
Let $\scrX$, $\scrA$ and $M^1(\scrX)$ satisfy the minimal
conditions. For every $\al\in\scrA$, let $\Pi_\al$ be a
distribution for a random histogram $P_\al$ as in
\cref{eq:randomhistogram}. Assume that the resulting
system of random histograms has the following property:
if $\al\leq\be$, then the distribution $\Pi_\al$
follows from $\Pi_\be$ through summation as in
\cref{eq:transitionmapping}, \ie,
\begin{equation}
  \label{eq:measuresyscoherence}
  \Pi_\al = \Pi_\be\circ\varphi_{*\,\al\be}^{-1}
\end{equation}
Then we refer to $(\Pi_\al,\varphi_{*\,\al\be})$ as a
\emph{coherent (inverse) system of random histogram
distributions}. If there exists a unique Radon probability
distribution $\Pi$ on $M^1(\scrX)$ with projections
$\Pi_\al$ for all $\al\in\scrA$, then $\Pi$ is called
its \emph{random histogram limit}.
\end{definition}
For later reference, we define mean measures for Borel
probability distributions on $M^1(\scrX)$. 
\begin{definition}
\label{def:meanmeasure}
Let $\scrX$, $\scrA$ and $M^1(\scrX)$ satisfy the minimal
conditions. Consider $M^1(\scrX)$ with a Borel probability
measure $\Pi$. The \emph{mean measure} $G$ under $\Pi$
is defined pointwise,
\[
  G(A) = \int_{M^1(\scrX)} P(A)\dd \Pi(P),
\]
for every Borel set $A$ in $\scrX$. Its restrictions to
the sub-$\sigma$-algebras $\sigma(\al)$ are denoted $G_\al$.
\end{definition}
To see that $G$ is a well-defined probability measure,
note that $\sigma$-additivity of $G$ is guaranteed by
monotone convergence. Also note that the restrictions
$G_\al:=G\bigr|_{\sigma(\al)}$ are mean
measures for the distributions $\Pi_\al$ on $M^1(\scrX_\al)$:
for any $A\in\al$,
\[
  G(A)=\int_{M^1(\scrX)}
    P(A)\dd \Pi(P)=\int_{M^1(\scrX)} P_\al(A)\dd \Pi(P)
    =\int_{M^1(\scrX_\al)} P_\al(A)\dd \Pi_\al(P_\al)
    = G_\al(A).
\]
\begin{remark}
\label{rem:Palphanotation} In the above we have abused notation
slightly: for any probability measure in $M^1(\scrX_\al)$, the
domain is $\sigma(\scrX_\al)$ rather than $\sigma(\al)\subset\scrB$.
So when we mean to refer to $\varphi_{*\,\al}(P)(\{e_i\})$, we shall
often use the more natural notation $P_\al(A_i)$ instead. 
% or, using \cref{eq:Pprimealpha}, $P'_\al(A_i)$.
\end{remark}

\subsection{The Bourbaki-Prokhorov-Schwartz theorem}
\label{sub:BPS}

The conditions we derive in subsequent sections are based on
a theorem from \citep{Schwartz73} (referred to as Prokhorov's
theorem in \citep{BourIII}), which says that the existence of
a limiting positive Radon measure in inverse systems of positive
measures is equivalent to a form of inner regularity that holds
for all projections simultaneously. This leads to characterization
of those inverse systems $(\Pi_\al,\varphi_{*,\al\be})$ that
consistently define Radon probability measures $\Pi$ on
$M^1(\scrX)$ with various topologies.

To discuss the Bourbaki-Prokhorov-Schwartz theorem, we
first have to generalize somewhat: let $\scrA$ be a directed set,
assume that $\scrY_\al$, $\al\in\scrA$, are Hausdorff topological
spaces and that for any $\al\leq\be$, there exist continuous,
surjective transition mappings $\psi_{\al\be}:\scrY_\be\to\scrY_\al$.
Together, they form an inverse system of Hausdorff spaces
(see \GTIthree{I}{4}{4}, \GTI{I}{2}{3}{Prop.~4}), denoted
$(\scrY_\al,\psi_{\al\be})$. If $T$ denotes a Hausdorff
topological space, a family of (projection) mappings
$\psi_\al:T\to\scrY_\al$, $\al\in\scrA$, is said
to be \emph{coherent} if, for all $\al\leq\be$,
$\psi_\al=\psi_{\al\be}\circ\psi_\be$,
and it is said to be \emph{separating}
if, for all $x,y\in T$, $x\neq y$, there exists an
$\al\in\scrA$ such that $\psi_\al(x)\neq \psi_\al(y)$.

\begin{theorem}
\label{thm:prokhorov}(Bourbaki-Prokhorov-Schwartz)\\
Let $(\scrY_\al,\psi_{\al\be})$ be an inverse system of
Hausdorff topological spaces indexed by $\al\in I$, $T$ a Hausdorff
topological space and $\psi_\al:T\to\scrY_\al$ a coherent and
separating family of continuous mappings. Let
$(\mu_\al,\psi_{\al\be})$ be a coherent inverse system of
positive measures on $(\scrY_\al,\psi_{\al\be})$. There exists
a bounded positive Radon measure $\mu$ on $T$ projecting to
$\mu_\al$ for all $\al\in I$, if and only if the following
property is satisfied:
\begin{itemize}
  \item[] for every $\ep>0$, there is a compact $H\subset T$
  such that for all $\al\in I$,
  \[
    \tag{P}\label{property:P}
    \mu_\al\bigl(\scrY_\al\setminus \psi_\al(H)\bigr)\leq\ep.
  \]
\end{itemize}\vspace*{-0.75em}
When $(P)$ holds, the measure $\mu$ is uniquely determined and
$\mu(L) = \inf\bigl\{\mu_\al(\psi_\al(L)):\al\in I\bigr\}$
for every compact set $L$ in $T$.
\end{theorem}
\begin{proof} See theorem~1 of \IIIthree{IX}{4}{2}.
\end{proof}
If all conditions of theorem~\ref{thm:prokhorov} are met, except
the system of functions $\psi_\al$ is not separating, then a
measure $\mu$ exists but may not be unique. 

\cite{BourIII} continues with application to a proof of
existence of Wiener measure and Kolmogorov's perspective, and the
definition of so-called \emph{promeasures} (also commonly
known as \emph{cylinder set measures}) which can be compared
with the coherent histogram systems we define below:
for a locally convex space $E$, \cite{BourIII}, Ch.\,{IX},
\S\,{6} considers the collection of all linear subspaces
$V$ of finite co-dimension in $E$ with continuous projections
$p_V:E\to E/V$ (and canonical $P_{VW}:E/W\to E/V$ for
$W\subset V$), to introduce $(E/V,p_{VW})$ as the
\emph{inverse system of finite-dimensional quotients}. 
A coherent system of positive measures $\Pi_V$ on the
finite-dimensional spaces $E/V$, $(\Pi_V,p_{VW})$, is
called a \emph{promeasure} on $E$. It is noted that
\IIIthree{IX}{6}{8--10} formulates a sufficient condition,
(\emph{Minlos's theorem}, \III{IX}{6}{10}{Theorem 2}
(based on \cite{Minlos63}), but it appears difficult to
apply unless $E$ is a (barrelled) \emph{nuclear space}.

In subsequent sections, we apply theorem~\ref{thm:prokhorov}
directly to spaces of (bounded/signed/{\-}positive/probability)
measures with various topologies, limiting the inverse
system of finite-dimensional quotients and promeasures,
to inverse systems of partitions and random histograms.
Let us prepare the discussion with some specifications
pertaining to the situation where $\scrX$, $\scrA$ and
$M^1(\scrX)$ satisfy the minimal conditions and
$I=\scrA$, $T=M^1(\scrX)$, $\scrY_\al=M^1(\scrX_\al)$
and $\psi_{\al\be}=\varphi_{*\,\al\be}$, in the form
of the following proposition.
\begin{proposition}
\label{prop:invsysmsr}
Let $\scrX$, $\scrA$ and $M^1(\scrX)$ satisfy the minimal
conditions. For all $\al\leq\be$, the mappings
$\varphi_{*\,\al\be}:M^1(\scrX_\be)\to M^1(\scrX_\al)$,
are continuous and surjective, and
$(M^1(\scrX_\al),\varphi_{*\,\al\be})$ forms an
inverse system of compact Hausdorff topological spaces, with
non-empty, compact Hausdorff inverse limit $N$.
\end{proposition}
\begin{proof}
Let $\al\leq\be$ be given. For any $g\in C(\scrX_\al)$, the mapping
$g\circ\varphi_{\al\be}:\scrX_\be\to\RR$ is an element of
$C(\scrX_\be)$. Because $\varphi_{\al\be}$ is surjective, the
induced mapping $\varphi^*_{\al\be}:C(\scrX_\al)\to C(\scrX_\be)$
is a bounded linear operator (with norm equal to one). The transpose
mapping $\varphi_{*\,\al\be}:M(\scrX_\be)\to M(\scrX_\al)$ is
defined by,
\begin{equation}
  \label{eq:defofvarphistaralbe}
  \langle \varphi_{*\,\al\be}(\mu_\be) , g \rangle_\al
    = \langle \mu_\be , \varphi^*_{\al\be}(g) \rangle_\be
    = \langle \mu_\be , g\circ\varphi_{\al\be} \rangle_\be,
\end{equation}
for all $\mu_\be\in M(\scrX_\be)$ and $g\in C(\scrX_\al)$. The linear
mapping $\varphi_{*\,\al\be}$ is bounded (with norm less than or equal
to one) and surjective. Note that if we express $\mu_\be$ as a
vector $(\mu_{\be,1},\ldots,\mu_{\be,{|\be|}})$ in $\RR^{{|\be|}}$,
\[
  \langle \mu_\be , g\circ\varphi_{\al\be} \rangle_\be
    = \sum_{j\in I(\be)} \mu_{\be,j}\, g(\varphi_{\al\be}(e_j))
    %= \sum_{j\in I(\be)} \mu_{\be,j}\, g(e_{i_{\al\be}(j)})
    = \sum_{i\in I(\al)}\Bigl( \sum_{j\in J_{\al\be}(i)}\mu_{\be,j}
      \Bigr) g(e_i)
    = \sum_{i\in I(\al)}\varphi_{*\,\al\be}(\mu_\be)_i g(e_i).
\]
in accordance with (\ref{eq:transitionmapping}). Finally, it is
noted that inverse limits of non-empty, compact spaces are non-empty
and compact (see \GTI{I}{9}{6}{Prop.~8}).
\end{proof}
The space $N$ consists of \emph{finitely additive} probability
set-functions on the $\sigma$-algebra generated by the partitions
in $\scrA$. Existence theorems for inverse limit probability measures
on associated inverse limit spaces like $N$ have been studied
extensively: Bochner's theorem \citep{Bochner55} and Choksi's
theorem \citep{Choksi58} give relatively mild sufficient conditions
for the existence of a limiting probability measure $\Pi$ on $N$
for inverse systems of Radon probability spaces (see also,
\cite{Rao71,Rao81,Beznea14}). But although $M^1(\scrX)$ (with
the weak topology of section~\ref{sec:weakhistogramlimits}) is
homeomorphic to a subspace of $N$, it has proven difficult
to formulate an additional condition to specify that $\Pi$
is concentrated on the image of $M^1(\scrX)$ in $N$ (see,
however, \cite{Ferguson73,Ghosh03,Ghosal17} and the correct
proof of the \emph{mean-measure condition} in \cite{Orbanz11}).
One of the strengths of the Bourbaki-Prokhorov-Schwartz theorem
is that $T=M^1(\scrX)$ is projected directly onto the spaces
$\scrY_\al=M^1(\scrX_\al)$, without detour via the inverse
limit $N$. In that way, theorem~\ref{thm:prokhorov} avoids
the (attractive but misleading) suggestion that a probability
distribution $\Pi$ on $N$ is an easy way to get `close to'
the desired distribution on $M^1(\scrX)$. By insisting only
on continuous projections $T\to\scrY_\al$,
theorem~\ref{thm:prokhorov} focusses on inner regularity as
the central issue. (Compare \citep{Schwartz73}, Theorem~21
and \III{IX}{4}{2}{Theorem~1}.)

%%%%%%%%%%%%%%%%%%%%%%%%%%%%%%%%%%%%%%%%%%%%%%%%%%%%%%%%%%%%%%%%%%%%%%%%%%%%%%%
\section{Random histogram limits with the weak topology}
\label{sec:weakhistogramlimits} % was {sec:lecamschwartz}

Again, let $\scrX$, $\scrA$ and $M^1(\scrX)$ satisfy the minimal
conditions, and fix the topology on $M^1(\scrX)$ to be the
\emph{weak topology}, defined as the
subspace topology that $M^1(\scrX)$ inherits from $M_b(\scrX)$
with the weak topology. Compactness of a subset of $M^1(\scrX)$
is characterized by the \emph{Dunford-Pettis-Grothendieck theorem}
(as presented, for example, in \citep{LeCam86}, Theorem~6 of
Appendix~8): a subset $H$ of $M^1(\scrX)$ is relatively compact
in the weak topology, if and only if, for some $Q\in M^1(\scrX)$,
\[
  \sup_{P\in H}\,\bigl\| P-P\wedge L\,Q \bigr\|_{1,\scrX} \to 0,
\]
as $L\to\infty$. The more classical characterization of weak
compactness, is that of the \emph{Dunford-Pettis theorem}
(see \cite{Diestel91}), which says that a $Q$-dominated subset
$H$ of $M^1(\scrX)$ is relatively compact in the weak topology,
if and only if, the subset of Radon-Nikodym densities
$\{dP/dQ:P\in H\}$ in $L^1(\scrX,\scrB,Q)$ is \emph{uniformly
integrable}, \ie,
\[
  \sup_{P\in H}\,\int_{\{x\in\scrX:dP/dQ(x)>L\}}
    \frac{dP}{dQ}(x)\dd Q(x) \to 0,
\]
as $L\to\infty$. 
And finally, we may characterize relative compactness of a
subset $H$ of $M^1(\scrX)$ in the weak topology, by the
condition that there exists a $Q\in M^1(\scrX)$ such that for
every $\ep>0$, there exists a $\delta>0$ such that,
\begin{equation}
  \label{eq:uniformintegrabilitythree}
  Q(A)<\delta \,\,\,\,\Rightarrow\,\,\,\, \sup_{P\in H} P(A)<\ep,
\end{equation}
for every measurable $A\subset\scrX$. These three
characterizations are equivalent.

\subsection[Support and approximation]{Support and approximation
of weak histogram limits}
\label{sub:WeakSupport}

Before we apply theorem~\ref{thm:prokhorov} to define Radon
probability measures on $M^1(\scrX)$ with the weak
topology (and total-variational topology), let us consider some
consequences, that is, necessary conditions for the existence of a
random histogram limit. First, we characterize the support of Borel
probability measures, next we consider approximations of random $P$
by random $P_\al$.

\subsubsection{Support and domination}

The following lemma is immediate but central enough to emphasize. 
\begin{lemma}
\label{lem:WeakSupport}
Let $\scrX$, $\scrA$ and $M^1(\scrX)$ satisfy the minimal
conditions. Consider $M^1(\scrX)$ with a Borel probability measure
$\Pi$. For any Borel set $A$ in $\scrX$, $G(A)=0$ implies that
$\Pi(\{P\in M^1(\scrX):P(A)>0\})=0$.
\end{lemma}
\begin{proof}
Let a Borel set $A$ in $\scrX$ be given and assume that $G(A)=0$.
If the Borel set $B=\{P\in M^1(\scrX):P(A)>0\}$ in $M^1(\scrX)$ has
probability $\Pi(B)>0$, then by $\sigma$-additivity, for some $\ep>0$
the Borel set $B'=\{P\in M^1(\scrX):P(A)>\ep\}$ has probability
$\Pi(B')>0$. That would imply that
$G(A)=\int P(A)\dd\Pi\geq\int_{B'}P(A)\dd \Pi(P)\geq \ep \Pi(B') > 0$,
contradicting the assumption.
\end{proof}
Domination by the mean measure plays a role in the following
proposition concerning the support of weakly-Borel probability
measures on $M^1(\scrX)$.
\begin{proposition}
\label{prop:WeakSupport}
Let $\scrX$, $\scrA$ and $M^1(\scrX)$ satisfy the minimal
conditions. Consider $M^1(\scrX)$ with the weak topology
and a Borel probability distribution $\Pi$. Let $G$ be the mean
measure under $\Pi$. Then $\{P\in M^1(\scrX):P\ll G\}$ is closed in
$M^1(\scrX)$ and,
\[
  \supp_{\scrT_1}(\Pi) \subset \{P\in M^1(\scrX):P\ll G\}.
\]
Moreover, if $P\in M^1(\scrX)$ is such that for all
measurable partitions $\al\in\scrA$, $P_\al$ lies in
the support of $\Pi_\al$ in $M^1(\scrX_\al)$, then $P$
lies in the weak support of $\Pi$.
\end{proposition}
\begin{proof}
If $P\in M^1(\scrX)$ is not dominated by $G$, then there exists a
Borel set $A$ such that $P(A)>0=G(A)$. Consequently for small
enough $\ep'>0$, the  weakly open neighbourhood
$U=\{Q\in M^1(\scrX):|Q(A)-P(A)|<\ep'\}$ does not meet
$\{Q\in M^1(\scrX):Q\ll G\}$, so $\{Q\in M^1(\scrX):Q\ll G\}$
is weakly closed. According to lemma~\ref{lem:WeakSupport},
$\Pi(U)\leq\Pi(\{Q\in M^1(\scrX):Q(A)>0\})=0$, so $U$
receives $\Pi$-mass zero, implying that
$P\not\in\supp_{\scrT_1}(\Pi)$. 

Regarding the last assertion, it is noted that, since
$M^1(\scrX)$ with the weak topology is homeomorphic to
a subspace of the inverse limit $N$ of
proposition~\ref{prop:invsysmsr}, the collection of sets
$\{\varphi_{*\,\al}^{-1}(V):\al\in\scrA,V\in\scrU_\al\}$
(where $\scrU_\al$ is any topological basis for $M^1(\scrX_\al)$)
in $M^1(\scrX)$ forms a basis for the weak topology.
Consequently, for any weak neighbourhood $U$ of
$P\in M^1(\scrX)$ there exists an $\al\in\scrA$ and
a $V\in\scrU_\al$ such that $\varphi_{*\,\al}^{-1}(V)\subset U$,
and,
\[
  \Pi(U)\geq\Pi(\varphi_{*\,\al}^{-1}(V))=\Pi_\al(V)>0,
\]
by assumption.
% BOOK ONLY
% let $P\in M^1(\scrX)$ be given and
% assume that $P_\al\in\supp(\Pi_\al)$ for all measurable partitions
% $\al\in\scrA_0$. If $U$ is an open neighbourhood of $P$, there
% exist $\ep>0$, $k\geq1$ and measurable $\phi_l:\scrX\to[0,1]$,
% $1\leq l\leq k$, such that,
% \[
%   \Biggl\{\, Q\in M^1(\scrX):
%     \max_{1\leq l\leq k}\biggl|\int\phi_l\dd (P-Q)\biggr|<\ep
%     \,\Biggr\}\subset U.
% \]
% For every $1\leq l\leq k$, there exists a finite measurable
% partition $\al_l$ and $0\leq f_{l,1},\ldots,f_{l,N(\al_l)}\leq 1$
% such that the simple function $f_l:\scrX\to[0,1]$,
% $f_l(x)=\sum_{m=1}^{N(\al_l)} f_{l,m} 1_{A_l}(x)$ approximates
% $\phi_l$ uniformly over $\scrX$ to within $\ep/4$. On the (finite,
% measurable) partition $\al\in\scrA_0$ that consists of the intersections
% of all elements of the partitions $\al_1,\ldots,\al_k$, all simple
% approximations can be realized simultaneously. Therefore, for every
% $1\leq l\leq k$,
% \[
%   \biggl|\int\phi_l\dd (P-Q)\biggr|
%     \leq \biggl|\int f_l\dd (P-Q)\biggr|+\frac{\ep}2
%     \leq \max_{1\leq i\leq {|\al|}}\bigl|(P_\al-Q_\al)(A_i)\bigr|+\frac{\ep}2
% \]
% Noting that $\max_{1\leq i\leq {|\al|}}\bigl|(P_\al-Q_\al)(A_i)\bigr|$
% can serve as the norm on $M(\scrX_\al)$ such that $M^1(\scrX_\al)$
% has the subspace topology, we see that for the open norm ball $B_\al$
% of radius $\ep/2$ around $P_\al$ in $M^1(\scrX_\al)$,
% we have $\Pi(U)\geq \Pi_\al( B_\al )$, and $\Pi_\al(B_\al)>0$
% because $P_\al$ lies in the support of $\Pi_\al$. Conclude that
% $P\in\supp_{\scrT_1}(\Pi)$. 
\end{proof}
Recall the open problem regarding construction of random
elements from specific dominated families, as posed in
\cite{daley07}, on p.~42:
\begin{quotation}
  \noindent\emph{Indeed, it appears to be an open problem to find
    simple sufficient conditions, analogous to Corollary~9.3.VI, for
    the realizations of a random measure to be [almost-surely]
    absolutely continuous with respect to a given measure.}
\end{quotation}
(\cite{daley07}, Corollary~9.3.VI formulates a condition
for a system of random histograms with limit, to be almost-surely
non-atomic.) Proposition~\ref{prop:WeakSupport} says that, given
some probability measure $G$, we should look for coherent systems
of random histograms with the projections $G_\al$ as their mean
histograms, and a limit $\Pi$ that is a weakly-Borel probability
measure on $M^1(\scrX)$.  In subsection~\ref{sub:WeakExistence},
we provide a relatively simple necessary and sufficient condition
for a coherent random histogram system to have a unique
weakly-Radon random histogram limit $\Pi$.

\subsubsection{Approximation by weakly convergent histograms}

Next, we consider the way in which histogram systems with
a weak limit approximate (random) probabilities $P(A)$.
Let $\Pi$ be a \emph{Radon} probability measure on $M^1(\scrX)$ with
the weak topology, with mean measure $G$. For
every $\eta>0$ and every $\al\in\scrA$, let $\scrC_{G,\al}(\eta)$
denote the collection of all Borel sets $B$ in $\scrX$
that are approximated by elements of the $\sigma$-algebra
$\sigma(\al)$ to within $G$-measure $\eta$:
\[
  \scrC_{G,\al}(\eta) = \bigl\{ B\in\scrB:
    \inf\{G((B\setminus C)\cup(C\setminus B))
    :C\in\sigma(\al)\}<\eta \bigr\}.
\]
Note that for any $\eta>0$ and any $B$, there exists an $\al$
such that $B\in\scrC_{G,\al}(\eta)$ (see theorem~4.4 in
\citep{Kingman66}).
\begin{proposition}
Let $\scrX$, $\scrA$ and $M^1(\scrX)$ satisfy the minimal
conditions. If $\Pi$ is a weakly-Radon probability measure
on $M^1(\scrX)$ with mean measure $G$, then for every
$\delta,\ep>0$, there exists a partition $\al\in\scrA$ and
an $\eta>0$, such that for all $B\in\scrC_{G,\al}(\eta)$,
\[
  \Pi\Bigl( \bigl\{P\in M^1(\scrX):
    \inf\{P((B\setminus C)\cup(C\setminus B))
    :C\in\sigma(\al)\}>\delta \bigr\} \Bigr) < \ep.
\]
\end{proposition}
\begin{proof}
Let $\ep>0$ be given. By inner regularity, there exists a
weakly compact $H$ in
$M^1(\scrX)$ such that $\Pi(H)>1-\ep$. For every $\delta>0$,
there exists an $\eta>0$, such that for all Borel sets $A$
in $\scrX$, $G(A)<\eta$ implies that $P(A)<\delta$ for all
$P\in H$, \cf\ (\ref{eq:uniformintegrabilitythree}). In
particular, if $B\in\scrC_{G,\al}(\eta)$,
then for some $C\in\sigma(\al)$,
$G((B\setminus C)\cup(C\setminus B))<\eta$, implying that
for all $P\in H$, $P((B\setminus C)\cup(C\setminus B))<\delta$.
\end{proof}
This observation is important from a computational perspective:
the practitioner chooses an approximating partition $\al$ to
perform computations with histograms and would like to
be able to control accuracy of his approximations for the
$P$ in terms of their restrictions to $\sigma(\al)$. He has
control over the probability measures $\Pi_\al$, and as a
result, control over the mean measures $G_\al$. Accordingly,
he can choose a level of refinement (expressed by a choice
for some partition $\al$), making approximations in $G$-measure
by $\al$-histogram. The Radon property ensures
that the approximation in $G$-measure carries over to
approximation in $P$-measure, uniformly in $P$, with
arbitrarily high $\Pi$-probability, depending on the
degree of approximation in the level $\al$ that is chosen
for actual computations. Such a guarantee concerning degrees
of approximation is not automatic if $\Pi$ is a Radon measure
for the tight topology of section~\ref{sec:prokhorov}.

\subsection{Existence of weak histogram limits}
\label{sub:WeakExistence}

Let $\scrA$ denote a set of finite Borel measurable
partitions of $\scrX$, directed for ordering by refinement.
If we equip $T=M^1(\scrX)$ with the weak topology,
theorem~\ref{thm:prokhorov} takes the following form.
% (Below we abbreviate $\varphi_{*\,\al}(Q)=Q_\al$,
% for $Q\in M^1(\scrX)$ and $\al\in\scrA$.)
\begin{theorem}
\label{thm:WeakExistence} % was {thm:ToneExistence}
Let $\scrX$, $\scrA$ and $M^1(\scrX)$ satisfy the minimal
conditions. Assume that $\scrA$ resolves $\scrX$ and
consider $M^1(\scrX)$ with the weak topology. Let
$(\Pi_\al,\varphi_{*\,\al\be})$ be a coherent system
of Borel probability measures on the inverse system
$(M^1(\scrX_\al),\varphi_{*\,\al\be})$. There exists a
unique weakly-Radon probability measure $\Pi$ on $M^1(\scrX)$
projecting to $\Pi_\al$ for all $\al\in\scrA$, if and only if,
\begin{itemize}
\item[] there is a $Q\in M^1(\scrX)$ such that, for
  every $\ep,\delta>0$ there is a $L>0$ such that,
  \[
    \tag{P-weak}\label{property:Pws}
    \Pi_\al\bigl(\{ P_\al\in M^1(\scrX_\al):
      \|P_\al-P_\al\wedge LQ_\al\|_{1,\scrX_\al}
      > \delta \}\bigr) < \ep,
  \]
  for all $\al\in\scrA$.
\end{itemize}
% When condition~(\ref{property:Pws}) holds,
% $\Pi(H) = \inf\bigl\{\Pi_\al(\varphi_{*\,\al}(H)):\al\in\scrA\bigr\}$,
% for every compact set $H$ in $M^1(\scrX)$.
\end{theorem}
\begin{proof}
According to \cref{prop:invsysmsr}, 
$(M^1(\scrX_\al),\varphi_{*\,\al\be})$ forms an
inverse system of Hausdorff topological spaces.
For all $\al\in\scrA$, $\varphi_{*\,\al}:M^1(\scrX)\to M^1(\scrX_{\al})$ 
is continuous with respect to the weak topology. 
If $P,Q\in M^1(\scrX)$ and $P\neq Q$, then there exists a set
$B$ in the $\sigma$-algebra generated by the $\sigma(\al)$ such that
$P(B)\neq Q(B)$, which cannot be the case unless, for some
$\al\in\scrA$, the histogram projections $\varphi_{*\,\al}(P)$ and
$\varphi_{*\,\al}(Q)$ differ. Combining this with
\cref{eq:coherentfamilyfunctions}, we conclude that
$(\varphi_{*\,\al},\varphi_{*\,\al\be})$ forms a coherent and
separating family of continuous mappings $M^1(\scrX)\to M^1(\scrX_\al)$.

The assertion now follows from theorem~\ref{thm:prokhorov} if we can show
that condition~(\ref{property:P})
holds. To that end,
let $\ep>0$ be given and define $\ep_n=2^{-n}\ep$. Given some
monotone decreasing sequence $(\delta_n)$ such that $\delta_n>0$,
$\delta_n\to0$, let $L_n$ be positive constants such that,
\[
  \Pi_\al\Bigl(\bigl\{ P_\al\in M^1(\scrX_\al)\,:\,
    \|P_\al-P_\al\wedge L_n\,Q_\al\|_{1,\scrX_\al} > \delta_n \bigr\}\Bigr) < \ep_n,
\]
for every $\al\in\scrA$. Define,
\[
  H = \bigcap\bigl\{ P\in M^1(\scrX)\,:\,
    \|P_\al-P_\al\wedge L_n\,Q_\al\|_{1,\scrX_\al} \leq \delta_n,\, n\geq1,\, \al\in\scrA \bigr\}.
\]
Let $\delta>0$ be given,
choose $n\geq1$ such that $\delta_n<\delta$ and define $L=L_n$.
Since $\scrA$ resolves $\scrX$, proposition~\ref{prop:alphatotalvariation}
and lemma~\ref{lem:DPGmatch} say that $\|P-P\wedge L\,Q\|_{1,\scrX}=
\sup\{\|P_\al-P_\al\wedge L\,Q_\al\|_{1,\scrX_\al}:\al\in\scrA\}$
for all $P$, and hence,
\[
  \sup\{\|P-P\wedge L\,Q\|_{1,\scrX}:P\in H\}
    =\sup\{\|P_\al-P_\al\wedge L\,Q_\al\|_{1,\scrX_\al}:P\in H,\al\in\scrA\}
    \leq\delta,
\]
Conclude that $H$ is relatively compact with respect to the
weak topology, \cf\ the
Dunford-Pettis-Grothendieck condition.  For the compact closure
$\overline{H}$ of $H$ in $M^1(\scrX)$ and any $\al$, we have (by
monotony of $\al\mapsto\|P_\al-P_\al\wedge L\,Q_\al\|_{1,\scrX_\al}$
for any $L>0$),
\[
  \begin{split}
  \Pi_\al\bigl(
    M^1(\scrX_\al&)\setminus\varphi_{*\,\al}(\overline{H}) \bigr)
  \leq \Pi_\al\bigl( M^1(\scrX_\al)\setminus\varphi_{*\,\al}(H)
    \bigr)\phantom{\biggm|}\\
  &\leq \Pi_\al\biggl( M^1(\scrX_\al)\setminus
    \bigcap_{n\geq1} \varphi_{*\,\al}\Bigl(\bigl\{ P\in M^1(\scrX)\,:\,
    \|P_\al-P_\al\wedge L_n\,Q_\al\|_{1,\scrX_\al} \leq
    \delta_n \bigr\}\Bigr) \biggr)\\
  &= \Pi_\al\biggl( \bigcup_{n\geq1}\{P_\al\in M^1(\scrX_\al)\,:\,
    \|P_\al-P_\al\wedge L_n\,Q_\al\|_{1,\scrX_\al} > \delta_n \}\biggr)\\  
  &\leq \sum_{n\geq1} \Pi_\al\bigl( \{P_\al\in M^1(\scrX_\al)\,:\,
    \|P_\al-P_\al\wedge L_n\,Q_\al\|_{1,\scrX_\al} > \delta_n \}\bigr) < \ep,
  \end{split}
\]
which shows that condition~(\ref{property:P}) of theorem~\ref{thm:prokhorov}
is satisfied. Conclude that there exists a unique Radon probability measure
$\Pi$ on $M^1(\scrX)$ that projects to $\Pi_\al$ for all
$\al\in\scrA$.

Conversely, let $\Pi$ be a weakly-Radon probability measure $\Pi$
on $M^1(\scrX)$. According to proposition~\ref{prop:WeakSupport},
the weak support of $\Pi$ is dominated by the mean measure $G$.
Again appealing to lemma~\ref{lem:DPGmatch}, we see that
that for every $\delta>0$, all $L>0$ and all $\al\in\scrA$,
\[
  \begin{split}
    \Pi\Bigl(\bigl\{ P\in & M^1(\scrX):
      \| P-P\wedge L\,G \|_{1,\scrX} > \delta \bigr\}\Bigr)\\[2mm]
    & = \Pi\Bigl(\bigl\{ P\in M^1(\scrX):
      \sup_{\be\in\scrA}\| P_\be-P_\be\wedge L\,G_\be\|_{1,\scrX_\be}
        > \delta \bigr\}\Bigr)\\
    &\geq \Pi\Bigl(\bigl\{ P\in M^1(\scrX):
      \| P_\al-P_\al\wedge L\,G_\al \|_{1,\scrX} > \delta \bigr\}\Bigr)\\[2mm]
    & = \Pi_\al\Bigl(\bigl\{ P_\al\in M^1(\scrX_\al):
      \| P_\al-P_\al\wedge L\,G_\al \|_{1,\scrX_\al} > \delta \bigr\}\Bigr).
  \end{split}
\]
Since $\Pi$ is weakly-Radon, for every $\delta,\ep>0$ there
exists a constant $L>0$ such that,
\[
  \Pi\Bigl(\bigl\{ P\in M^1(\scrX):
    \| P-P\wedge L\,G \|_{1,\scrX} > \delta \bigr\}\Bigr) < \ep,
\]
verifying that condition~(\ref{property:Pws}) holds.
\end{proof}

On first sight, condition~(\ref{property:Pws}) may appear technical
and inaccessible. It is noted, however, that in practice one has
considerable control: one may choose $\scrA$ (large enough
to resolve $\scrX$ but otherwise) as small as possible, and the
histogram system $(\Pi_\al:\al\in\scrA)$ as simple as possible,
in order to enable verification of condition~(\ref{property:Pws})
in manageable form. Moreover, all subsequent calculations involve
only probability distributions on finite-dimensional simplices,
enhancing feasibility greatly.

Regarding the measure $Q$, we simplify by 
appeal to a necessary condition: it is clear that \emph{if}
theorem~\ref{thm:WeakExistence} holds, then the support of the
probability measure $\Pi$ is dominated by the mean measure $G$,
\cf\ proposition~\ref{prop:WeakSupport}. So, when looking for
a candidate dominating measure $Q$ to verify
condition~(\ref{property:Pws}), we can turn to the mean measures
$G_\al$ of the histogram distributions $\Pi_\al$: if we show,
either, that the $G_\al$ are the histograms associated to a mean
measure $G$ (the \emph{mean measure condition} \citep{Orbanz11}),
or, that condition~(\ref{property:Ptight}) below holds (\cf\
definition~\ref{def:meanmeasure}), then $G$ can play the role of
$Q$.

Based on those two remarks, condition~(\ref{property:Pws}) can
be rewritten as follows: 
\begin{itemize}
\item[] for some $G\in M^1(\scrX)$, $\al$-histograms
  $\varphi_{*\,\al}(G)$ equal the mean measures $G_\al$
  for all $\al\in\scrA$, and, for every $\ep,\delta>0$ there
  is a $L>0$ such that,
  \[
    \tag{P-weak'}\label{property:Pweakprime}
    \Pi_\al\biggl(\Bigl\{ P_\al\in M^1(\scrX_\al):
      \sum_{A\in\al}\bigl(P_\al(A)-LG_\al(A)\bigr)_+ > \delta
      \Bigr\}\biggr) < \ep,
  \]
  for all $\al\in\scrA$.
\end{itemize}
This form of the condition forms the starting point for
the examples of sections~\ref{sec:polya} and~\ref{sec:gauss}.

\subsection{Existence of total-variational histogram limits}
\label{sub:totalvariation}

Let $\scrX$ be a Hausdorff topological space and let $\scrP$
be a subset of $M^1(\scrX)$ dominated by a probability measure
$Q\in M^1(\scrX)$. In this subsection, we distinguish $\scrP$
from $\scrP_b$, represented by the same set, denoted $\scrP$
when equipped with the weak topology, and $\scrP_b$ when
equipped with the total-variational topology. The identity
mapping $i:\scrP_b\to\scrP$ is a continuous bijection. Write
$\scrB(\scrP)$ and $\scrB(\scrP_b)$ for the associated Borel
$\sigma$-algebras.
\begin{proposition}
If $\scrX$ is separable and $\scrP$ is a
dominated subset of $M^1(\scrX)$, then $\scrP_b$ is separable
and $\scrB(\scrP)=\scrB(\scrP_b)$.
\end{proposition}
\begin{proof}
Let $Q$ denote the probability measure that dominates $\scrP$.
By separability of $\scrX$, the Banach space $L^1(\scrX,\scrB,Q)$
of $Q$-integrable, real-valued functions on $\scrX$,
is separable, and so is its subspace of Radon-Nikodym densities
$\scrP'_b=\{dP/dQ:P\in\scrP\}$. Since $\scrP_b$ and $\scrP'_b$
are homeomorphic, $\scrP_b$ is separable too. It can be shown
\citep{pfanzagl69,Strasser85} that then the total-variational
norm is measurable with respect to the minimal $\sigma$-algebra
for measurability of the mappings $P\mapsto P(A)$, $A\in\scrB$,
which is contained in $\scrB(\scrP)$. Accordingly,
$\scrB(\scrP_b)\subset \scrB(\scrP)$. Since the total-variational
topology refines the weak topology, also $\scrB(\scrP)\subset
\scrB(\scrP_b)$.
\end{proof}
As a consequence, any Borel probability measure $\Pi$ on
$\scrP$ (viewed as a $\sigma$-additive set-function with
$\scrB(\scrP)$ as its domain) gives rise to a Borel probability
measure $\Pi_b$ on $\scrP_b$ (viewed as a $\sigma$-additive
set-function with $\scrB(\scrP_b)$ as its domain).
\begin{corollary}
\label{cor:TVExistence}
Let $\scrX$ be separable and, together with $\scrA$ and
$M^1(\scrX)$, satisfy the minimal conditions. Consider $M^1(\scrX)$
with the total-variational topology. Assume that $\scrA$ resolves
$\scrX$. Let $(\Pi_\al,\varphi_{*\,\al\be})$ be a coherent system
of Borel probability measures on the inverse system
$(M^1(\scrX_\al),\varphi_{*\,\al\be})$. If
condition~(\ref{property:Pws}) holds, there exists a unique
$\scrT_{TV}$-Radon probability measure $\Pi$ on $M^1(\scrX)$
projecting to $\Pi_\al$ for all $\al\in\scrA$.
\end{corollary}
\begin{proof}
Under stated conditions, theorem~\ref{thm:WeakExistence} asserts
the existence of a Radon probability measure $\Pi$ on
$M^1(\scrX)$ with the weak topology,
and \cref{prop:WeakSupport} guarantees that
$\scrP=\supp_{\scrT_1}(\Pi)$ is dominated by the mean measure
$G$. Because $\scrX$ is separable, $\scrB(\scrP)=\scrB(\scrP_b)$,
so $\Pi$ is a Borel probability measure on $\scrP_b$. Uniqueness
follows from the uniqueness of the weak histogram limit.
By the Radon-Nikodym theorem, $\scrP_b$ is homeomorphic
(isometrically, even, \cf\ (\ref{eq:TVLOneisometry})) to a
closed subset of the Polish space $L^1(\scrX,\scrB,G)$, and
therefore a Radon space, so that $\Pi$ is a Radon measure.
\end{proof}
So, remarkably, existence of a total-variational random
histogram limit does not impose stricter conditions than
existence of a weak random histogram limit; moreover, not
even inner regularity is lost in the transition from $\scrP$ to
$\scrP_b$. From the perspective of theorem~\ref{thm:phases}
this amplification is inconsequential, but events and statements
involving the total-variational norm are very common and
measurability of total-variational balls is crucial for many
applications (for example, in large-sample limits of posterior
distributions on metric spaces in non-parametric statistics
\citep{Ghosal17,Kleijn21}).

%%%%%%%%%%%%%%%%%%%%%%%%%%%%%%%%%%%%%%%%%%%%%%%%%%%%%%%%%%%%%%%%%%%%%%%%%%%%%%%
\section{Random histogram limits with the tight topology}
\label{sec:prokhorov} % was sec:prokhorov

The existence question of a limit for coherent random histogram
systems has been studied extensively with the \emph{tight topology}
for $M^1(\scrX)$: a rich literature has grown from Kingman's
original work on \emph{completely random measures}
\citep{Kingman67}, with an emphasis on limits with almost-surely
purely atomic realizations \citep{daley03,daley07}. Here we revisit the
existence problem without restricting to point-processes,
and derive a necessary and sufficient condition
in subsection~\ref{sub:TightExistence}, based on
theorem~\ref{thm:prokhorov}. In subsection~\ref{sub:suppweak} we
consider the support of tight random histogram limits.

\subsection{Existence of tight histogram limits}
\label{sub:TightExistence}

Let $\scrX$ be a Polish space with topology $\scrT$. We are
interested in the construction of Radon probability measures
on $M^1(\scrX)$ with the tight topology. In comparison
with the construction of subsection~\ref{sub:WeakExistence}, the
assertion is weaker since the weak topology refines the tight
topology. Accordingly, compactness as in
condition~(\ref{property:P}) constitutes a less stringent
restriction, while the continuity requirement of histogram
projections becomes harder to satisfy. 

Indeed, when one tries to reproduce the initial steps in the
proof of \cref{sub:WeakExistence} with the tight topology,
a disappointment awaits: if $\scrX=\RR^d$ with the standard
topology, for example, then for any partition $\al$ of $\scrX$
into two or more subsets, the projection mapping
$\varphi_{*\,\al}:M^1(\scrX)\to M^1(\scrX_\al)$ of
equation~(\ref{eq:projectionmapping}) is \emph{not continuous}:
superficially, it appears that \cref{thm:prokhorov} cannot
be applied.

In order to correct this, we refine to a
\emph{zero-dimensional} version $\hat{\scrX}$ of $\scrX$,
rendering projection mappings continuous for a
collection of partitions that is large enough to be
separating and resolving.
While this leaves the Borel $\sigma$-algebra unchanged, the
transition to $\hat{\scrX}$ does complicate the nature of tight
compactness in $M^1(\hat{\scrX})$. A counterexample at the end
of this subsection shows that this complication corresponds
directly to the precise way in which a coherent system of random
histograms can \emph{fail to have a tight limit}.

\subsubsection{Zero-dimensional refinements of Polish spaces}
\label{subsub:zerodim}

With a countable basis $\scrU$ for the topology $\scrT$, define
the topological %\varinindex{subbasis}{subbasis!topological},
sub-basis,
\begin{equation}
  \label{eq:zerodimsubbasis}
  \scrS = \{ U,\scrX\setminus U\,:\,U\in\scrU \},
\end{equation}
for a topology $\hat{\scrT}$ on the set $\scrX$ in which each
basis element $U\in\scrU$ is %\inindex{clopen}
clopen; denote the resulting topological
space by $\hat{\scrX}$.
% \begin{example}
% In the case of a
% %\varinindex{second countable}{topological space!second countable}
% second countable space $\scrX$, enumerate $\scrU=\{U_n:n\geq1\}$
% and define $\scrX_0=\scrX$ and $\scrX_n$, ($n\geq1$) to be the
% topological sum of $U_n$ and $\scrX_{n-1}\setminus U_n$.
% Note that there is a corresponding sequence of partitions
% $\scrA=\{\al_n:n\geq1\}$ (as in example~\ref{ex:defaultpartitionsR})
% that are generated by the basis $\scrU$,
% and that the $\scrX_n$ form an inverse system
% of topological spaces $(\scrX_n,\varphi_{mn})$ (where the
% set-theoretic identity on the set $\scrX$ acts as continuous
% projection mapping $\varphi_{mn}$ for all $m\leq n$) with inverse
% limit homeomorphic to the space $\hat{\scrX}$.
% \end{example}
\begin{proposition}
\label{prop:concretehatX} % was {thm:concretehatX}
The space $\hat{\scrX}$ is zero-dimensional and the identity
mapping $i:\hat{\scrX}\to\scrX$ is continuous. If $\scrU_1$ and
$\scrU_2$ are two bases for $\scrX$, the corresponding spaces
$\hat{\scrX}_1$ and $\hat{\scrX}_2$ are homeomorphic. If
$\scrX$ is Polish, then $\hat{\scrX}$ is also Polish.
\end{proposition}
\begin{proof}
The sub-basis $\scrS$ gives rise to a basis consisting
of clopen sets, so $\hat{\scrX}$ is zero-dimensional, and
the identity $i$ is continuous because $\hat{\scrT}$ refines
$\scrT$. Because any $U_2\in\scrU_2$ contains a
$U_1$ from the basis $\scrU_1$ and \vv, the identity
mapping on $\scrX$ is continuous from $\hat{\scrX}_1$ to
$\hat{\scrX}_2$ and also from $\hat{\scrX}_2$ to
$\hat{\scrX}_1$. Assuming $\scrX$ is Polish, the countable
product space $\scrX^{\NN}=\prod_{n\geq1}\scrX$ is Polish
(see \cite{Kechris94}, prop.~3.3) and has a diagonal
$\Delta=\{(x,x,\ldots)\in\prod_{n\geq1}\scrX:x\in\scrX\}$
that is a closed subspace, homeomorphic to $\scrX$. 
Enumerate the basis sets in $\scrU=\{U_i:i\geq1\}$ and
define $\scrY_n$ to be the refinement of $\scrX$ with $U_1,
\ldots,U_n$ made clopen (\eg\ $\scrY_1$ is the topological
sum of $U_1$ and $\scrX\setminus U_1$, \etc). The canonical
set-theoretic identification $i_n:\scrY_n\to\scrX$,
is continuous. The spaces $\scrY_n$ are all Polish, as the
topological sums of $U_n$ and $\scrX\setminus U_n$ (which
are Polish).
% and countable topological sums of Polish spaces are Polish spaces. 
The product space $\prod_{n\geq1}\scrY_n$ is Polish and the map
$j:\prod_{n\geq1}\scrY_n\to\scrX^{\NN}$ is a continuous
bijection. Then $j^{-1}(\Delta)$ is Polish and homeomorphic to
$\hat{\scrX}$.
\end{proof}
\begin{proposition}
\label{prop:msrstaysmsr}
Let $\scrX$ be a Polish space.
The Borel sets on $\scrX$ and $\hat{\scrX}$ are equal and
any set function $\mu$ is a (bounded/signed/{\-}positive/probability)
Borel measure on $\scrX$, if and only if, $\mu$ is a
(bounded/signed/{\-}positive/probability) Borel measure on $\hat{\scrX}$.
\end{proposition}
\begin{proof}
Note that the Borel $\sigma$-algebra on $\scrX$
generated by the basis $\scrU$ is identical to
the $\sigma$-algebra generated by $\scrU$ and its
complements, which form the sub-basis for $\hat{\scrX}$.
Conclude that $\scrX$ and $\hat{\scrX}$ have the same Borel
sets. Boundedness, signedness or positivity, being a probability
measure and countable additivity are then identical as properties
of set functions $\mu$ on the Borel $\sigma$-algebra.
\end{proof}
Proposition~\ref{prop:msrstaysmsr} implies the existence of
a bijective mapping $i_*$ with the following properties.
\begin{proposition}
\label{prop:M1hatXtoM1X}
The mapping $i_*:M_b(\hat{\scrX})\to M_b(\scrX)$ is a continuous
bijection. If $\scrX$ is Polish, $M^1(\scrX)$ and $M^1(\hat{\scrX})$
are Polish and $i_*^{-1}:M^1(\scrX)\to M^1(\hat{\scrX})$ is Borel
measurable.
\end{proposition}
\begin{proof}
Any bounded, continuous $f:\scrX\to\RR$ is also bounded,
continuous when viewed as $f:\hat{\scrX}\to\RR$, so there exists a
linear, injective mapping $j_*:C_b(\scrX)\to C_b(\hat{\scrX})$ of norm
one, and transpose to that, a bounded, injective, linear
$i_*=j_*^t:M_b(\hat{\scrX})\to M_b(\scrX)$ of norm one
(see \TVS{II}{6}{4}{proposition~5} and
\TVS{IV}{1}{3}{proposition~8}). As noted earlier, if $\scrX$ is
a Polish space, then so is $M^1(\scrX)$ (\III{IX}{5}{4}{prop.~10}),
so based on \cref{prop:concretehatX}, both
$M^1(\scrX)$ and $M^1(\hat{\scrX})$ are Polish spaces.
According to theorems~12.4 and~14.12 (\emph{Souslin's theorem}) in
\citep{Kechris94}, if $\scrX,\scrY$ are standard Borel
spaces and $f:\scrX\to\scrY$ is a Borel measurable injection,
then its inverse on $f(\scrX)$ is also Borel measurable. Applied
to $i_*$, this proves the last assertion.
\end{proof}
For all $\al\in\scrA$, define the mappings
$\hat\varphi_{*\,\al}:M_b(\hat{\scrX})\to M(\scrX_\al)$,
\begin{equation}
  \label{eq:msrproj}
  \hat\varphi_{*\,\al}(\mu)
    = \bigl( \mu(A_1),\ldots,\mu(A_{{|\al|}}) \bigr),
\end{equation}
that takes any bounded, signed Radon measure $\mu$ on
$\hat{\scrX}$ into its $\al$-histogram.
\begin{proposition}
\label{prop:clopencont}
Let $\al$ be a partition of $\scrX$ generated by a
basis $\scrU$, and let $\hat{\scrX}$ denote the associated
zero-dimensional version of $\scrX$. The mapping
$\hat{\varphi}_{*\,\al}:M_b(\hat{\scrX})\to M(\scrX_\al)$
is continuous for the tight topology.
\end{proposition}
\begin{proof}
For any partition $\al$ generated by the basis $\scrU$
(\cf\ definition~\ref{def:generatedpartition}), any
$A\in\al$ is clopen in $\hat{\scrX}$. For any clopen $A$,
% both $A$ and $\hat{\scrX}\setminus A$ are open, so 
$1_A$ is a bounded, continuous function on $\hat{\scrX}$.
Therefore $M^1(\hat{\scrX})\to\RR:P\mapsto P(A)$ is
continuous with respect to the tight topology and so
is $\hat{\varphi}_{*\,\al}$.
\end{proof}
Compactness in $\hat{\scrX}$ has a more stringent meaning
than in the original space $\scrX$: indeed, according
to \emph{Brouwer's theorem}, any compact $\hat{K}\subset\hat{\scrX}$
is a union of a subspace homeomorphic to the Cantor set with
isolated points. For example, with $\scrX=\RR$ in its standard
topology, $[0,1]$ is not compact in the space $\hat{\scrX}$.

\subsubsection{Tight histogram limits with zero-dimensional compacta}

Tight compactness in $M_b(\scrX)$ is characterized by
\emph{Prokhorov's theorem} (see \cite{Prokhorov56}), which
says that a subset $H$ of $M_b(\scrX)$ is relatively tightly
compact, if and only if,
\begin{itemize}
  \item[(i.)] $\sup\{\|\Phi\|_{1,\scrX}:\Phi\in H\}<\infty$,
  \item[(ii.)] for every $\ep>0$, there exists a compact
    $K\subset\scrX$ such that,
    \begin{equation}
      \label{eq:prokhorovscondition}
      \sup\{ |\Phi|(\scrX\setminus K):\Phi\in H\}<\ep.
    \end{equation}
\end{itemize}
We are now in a position to apply \cref{thm:prokhorov} to
$M^1(\scrX)$ (or rather, $M^1(\hat{\scrX})$) with the tight
topology. (Note: mention of the Radon property in the statement
of theorem~\ref{thm:tightexistence}, is accurate
but strictly speaking redundant, since $M^1(\scrX)$ is a
Radon space.)
\begin{theorem}
\label{thm:tightexistence} % was {thm:weakexistence}
Let $\scrX$ be a Polish space with a directed set $\scrA$ of
partitions  that resolves $\scrX$, generated by a basis that
gives rise to a zero-dimensional $\hat{\scrX}$. Consider
$M^1(\scrX)$ with the tight topology. Let
$(\Pi_\al,\varphi_{*\,\al\be})$ be a coherent system of
Borel probability measures on the inverse system
$(M^1(\scrX_\al),\varphi_{*\,\al\be})$. There exists a
unique Radon probability measure $\Pi$ on $M^1(\scrX)$
projecting to $\Pi_\al$ for all $\al\in\scrA$, if and only if,
\begin{itemize}
\item[] for every $\ep,\delta>0$ there is a compact
  $\hat{K}\subset\hat{\scrX}$ such that,
  \[
    \tag{P-tight}\label{property:Ptight}
    \Pi_\al\bigl(\bigr\{ P_\al\in M^1(\scrX_\al)\,:\,
      P_\al(\varphi_\al(\hat{K}))<1-\delta \bigr\}\bigr) < \ep,
  \]
  for all $\al\in\scrA$.
\end{itemize}
% When condition~(\ref{property:Ptight}) holds,
% $\Pi(H) = \inf\bigl\{\Pi_\al(\varphi_{*\,\al}(H)):\al\in\scrA\bigr\}$,
% for every tightly compact set $H$ in $M^1(\scrX)$.
\end{theorem}
\begin{proof}
Under the condition of the theorem, $\hat{\scrX}$, $\scrA$ and
$M^1(\hat{\scrX})$ satisfy the minimal conditions. Like before
(see \cref{prop:invsysmsr}), $(M^1(\scrX_\al),\varphi_{*\,\al\be})$
forms an inverse system of compact Hausdorff topological spaces.
As in the proof of \cref{thm:prokhorov}, $(\hat\varphi_{*\,\al},
\varphi_{*\,\al\be})$ is a coherent and separating family of
mappings on $M^1(\hat{\scrX})$, and from \cref{prop:clopencont}
we conclude that the $\hat\varphi_{*\,\al}$ are also continuous.

To show that condition~(\ref{property:P}) holds, let $\ep>0$ be given and
define $\ep_n=2^{-n}\ep$. Given some decreasing sequence $(\delta_n)$
such that $\delta_n>0$, $\delta_n\to0$, let $\hat{K}_n$ be compact
subsets of $\hat{\scrX}$ such that,
\[
  \Pi_\al\Bigl(\bigr\{ P_\al\in M^1(\scrX_\al)\,:\,
    P_\al(\varphi_\al(\hat{K}_n))<1-\delta_n \bigr\}\Bigr) < \ep_n,
\]
for every $\al\in\scrA$ and every $n\geq1$. Define,
\[
  H = \bigcap\bigl\{ P\in M^1(\scrX)\,:\,
    P_\al(\varphi_\al(\hat{K}_n))\geq 1-\delta_n,\, 
    n\geq1,\, \al\in\scrA \bigr\}.
\]
Let $\delta>0$ be given, choose $n\geq1$ such that $\delta_n<\delta$
and choose $\hat{L}=\hat{K}_n$. Since the Borel sets
$L_\al=(\varphi_\al^{-1}\circ\varphi_\al)(\hat{L})$ decrease as the
level of refinement of the partition $\al$ increases, and since
$\scrA$ resolves $\scrX$, $\hat{L}=\cap_{\al\in\scrA} L_\al$, and
$P(\hat{L})=\inf\{P_\al(\varphi_\al(\hat{L})):\al\in\scrA\}$ by
monotony. Conclude that $H$ is relatively compact with respect to the
tight topology, according to (\ref{eq:prokhorovscondition}). For
the compact closure $\overline{H}$ of $H$ in $M^1(\hat{\scrX})$
and any $\al$, we have (by monotony of $\al\mapsto
P_\al(\varphi_\al(K))$ for any $K$),
\[
  \begin{split}
  \Pi_\al\bigl(
    M^1(\scrX_\al&)\setminus\hat\varphi_{*\,\al}(\overline{H}) \bigr)
  \leq \Pi_\al\bigl(
    M^1(\scrX_\al)\setminus\hat\varphi_{*\,\al}(H) \bigr)\phantom{\biggm|}\\
  &\leq \Pi_\al\biggl( M^1(\scrX_\al)\setminus
    \bigcap_{n\geq 1} \hat\varphi_{*\,\al}\Bigl(\bigl\{ P\in M^1(\scrX)\,:\,
    P_\al(\varphi_\al(\hat{K}_n))\geq 1-\delta_n \bigr\}\Bigr) \biggr)\\
  &\leq \Pi_\al\biggl( \bigcup_{n\geq1} \{P_\al\in M^1(\scrX_\al)\,:\,
    P_\al(\varphi_\al(\hat{K}_n))<1-\delta_n \}\biggr)\\  
  &\leq \sum_{n\geq1} \Pi_\al\Bigl( \{P_\al\in M^1(\scrX_\al)\,:\,
    P_\al(\varphi_\al(\hat{K}_n))<1-\delta_n\}\Bigr) < \ep,\phantom{\biggm|}
  \end{split}
\]
which shows that condition~(\ref{property:P}) of theorem~\ref{thm:prokhorov}
is satisfied. Conclude that there exists a unique Radon probability measure
$\hat\Pi$ on $M^1(\hat{\scrX})$ that projects to $\Pi_\al$ for all
$\al\in\scrA$. The continuous mapping $i_*:M^1(\hat{\scrX})\to M^1(\scrX)$
of \cref{prop:M1hatXtoM1X} serves to define $\Pi=\hat\Pi\circ i_*^{-1}$,
a Radon probability measure on $M^1(\scrX)$, and $\Pi$ still projects
to $\Pi_\al$ for all $\al\in\scrA$.

Conversely, since $\scrX$ is Polish, \cf\
proposition~\ref{prop:M1hatXtoM1X}, $\hat{\scrX}$, $M^1(\hat{\scrX})$
and $M^1(\scrX)$ are Polish spaces, and the mapping\
$i_*^{-1}:M_b(\scrX)\to M_b(\hat{\scrX})$
is Borel measurable. Therefore, the mapping $\hat\Pi=\Pi\circ i_*$
defines a Borel probability measure on $M^1(\hat{\scrX})$, which is Radon
because $M^1(\hat{\scrX})$ is a Radon space. So, according to
Prokhorov's theorem, for every $\delta,\ep>0$, there exists a compact
$\hat{K}$ in $\hat{\scrX}$ such that,
\[
  \Pi\Bigl(\bigr\{P\in M^1(\hat{\scrX}):
    P(\hat{K})<1-\delta \bigr\}\Bigr)<\ep.
\]
With $\hat{K}_\al=(\varphi_\al^{-1}\circ\varphi_\al)(\hat{K})
\subset\hat{\scrX}$, for every $\al\in\scrA$, we have
$\hat{K}\subset\hat{K}_\al$ and accordingly,
\[
  \begin{split}
  \Pi\Bigl(\bigr\{P\in &M^1(\hat{\scrX}):
    P(\hat{K})<1-\delta \bigr\}\Bigr)
  \geq
  \Pi\Bigl(\bigr\{P\in M^1(\hat{\scrX}):
    P(\hat{K}_\al)<1-\delta \bigr\}\Bigr)\\
  &=
  \Pi_\al\Bigl(\bigr\{ P_\al\in M^1(\scrX_\al)\,:\,
    P_\al(\varphi_\al(\hat{K}))<1-\delta \bigr\}\Bigr),
  \end{split}
\]
for any $\delta>0$, which implies the converse.
\end{proof}
Let us paraphrase: to have a coherent inverse system of probability
measures for histograms define a limit $\Pi$ that is a Radon probability
measure on $M^1(\scrX)$ for the tight topology, we look for compacta
$\hat{K}$ in a zero-dimensional version of $\scrX$ that capture most
of the mass of the projected measures $P_\al$ with high
$\Pi_\al$-probability, uniformly in $\al\in\scrA$.

\subsubsection{Tight histogram limits with ordinary compacta}
\label{subsub:tightwithnormalcompacta}

In certain histogram systems (like those that define
Dirichlet process distributions), there is an easy way to prove
the \emph{mean measure condition} \citep{Orbanz11} (see
the proof of theorems~\ref{thm:existencediri}). In histogram
systems where this
condition is less or not accessible (like those that define
the P\'olya-tree distributions), zero-dimensional compacta in
the space $\hat{\scrX}$ are unwieldy, so we also provide a
re-formulation of theorem~\ref{thm:tightexistence} that relies
only on compacta in $\scrX$.

To avoid mention of the zero-dimensional space $\hat{\scrX}$,
we re-construct compacta $\hat{K}$ from decreasing sequences
of compacta in $\scrX$. Let $\scrA$ be a directed set of
partitions generated by a basis. For every $\al\in\scrA$
consider the topological space $\scrY_\al$ obtained from
$\scrX$ by declaring all sets $A\in\al$ clopen, \ie\
$\scrY_\al$ is the topological sum of all the partition sets
$A\in\al$ with their subspace topologies. Note that the
set-theoretic identity mapping on $\scrX$ is continuous as a
mapping $\pi_\al:\hat{\scrX}\to\scrY_\al$. (See also the proof of
proposition~\ref{prop:concretehatX}.)
\begin{lemma}
\label{lem:hatcompact}
A subset $\hat{K}$ is compact in $\hat{\scrX}$, if and
only if, $K_\al=\pi_\al(\hat{K})$ is compact in $\scrY_\al$,
for all $\al\in\scrA$. Conversely, given compact subsets
$K_\al$ of $\scrY_\al$ for all $\al\in\scrA$, the
subset $\cap_\al K_\al$ is compact in $\hat{\scrX}$.
\end{lemma}
\begin{proof}
Consider the product space $\scrY=\prod_{\al\in\scrA}\scrY_\al$. The
diagonal $\Delta=\{(x,x,x,\ldots)\in\scrY:x\in\scrX\}$
is a closed subset of $\scrY$, homeomorphic to $\hat{\scrX}$
and the mappings $\pi_\al$ are the canonical projections
$\scrY\to\scrY_\al$, applied after the homeomorphism
$\hat{\scrX}\to\Delta$. A compact $\hat{K}$ in $\hat{\scrX}$
has compact images $\pi_\al(\hat{K})$ in all $\scrY_\al$,
$\al\in\scrA$. Conversely, if $H$ is a subset of $\hat{\scrX}$
such that $\pi_\al(H)$ is compact in $\scrY_\al$ for all
$\al\in\scrA$, then $B=\prod_{\al\in\scrA} \pi_\al(H)$ is
compact in $\scrY$ by \emph{Tychonov's theorem}, and so is
the closed subspace $\Delta\cap B$. Set-theoretically,
$\pi_\al(H)=H$ for all
$\al\in\scrA$, implying that $\Delta\cap B$ is homeomorphic
to $H$, so $H$ is compact. Given compact subsets
$K_\al$ of $\scrY_\al$ for all $\al\in\scrA$, the
subset $\hat{K}=\cap_\al K_\al$ is compact as a subset of
any $\scrY_\al$, ($\al\in\scrA$), so $\hat{K}$ is compact
as a subset of $\hat{\scrX}$.
\end{proof}
% For the following corollary, we assume that the partitions
% $\al\in\scrA$ form an ordered sequence $(\al_n)$. This is not
% unreasonable in most applications, in light of
% example~\ref{ex:partitionscountablebasis}.
\begin{corollary}
\label{cor:tightexistence}{} % was {thm:weakexistence}
Let $\scrX$ be a Polish space with a countable basis $\scrU$
and a well-ordered sequence $\scrA$ of partitions
generated by the basis, that resolves $\scrX$. Consider
$M^1(\scrX)$ with the tight topology. Let
$(\Pi_\al,\varphi_{*\,\al\be})$ be a coherent system of
Borel probability measures on the inverse system
$(M^1(\scrX_\al),\varphi_{*\,\al\be})$. If,
\begin{itemize}
\item[] for all $\al\in\scrA$, all $A\in\al$ and all $\ep,\delta>0$,
  there is a $K\subset A$, compact in $\scrX$, such that,
  \[
    \tag{P-tight'}\label{property:Ptightprime}
    \Pi_\be\bigl(\bigr\{ P_\be\in M^1(\scrX_\be)\,:\,
      P_\be(\varphi_\be(K))< P_\be(\varphi_\be(A))-\delta
        \bigr\}\bigr) < \ep,
  \]
  for all $\be\in\scrA$ such that $\al\leq\be$,
\end{itemize}
then there exists a unique Radon probability measure $\Pi$ on
$M^1(\scrX)$ projecting to $\Pi_\al$ for all $\al\in\scrA$.
\end{corollary}
\begin{proof}
Enumerate the partitions in $\scrA$, $\scrA=\{\al_n:n\geq0\}$,
and let $\delta,\ep>0$ be given.
To find a compact subset $\hat{K}$ of $\hat{\scrX}$
to satisfy property~(\ref{property:Ptight}), we
construct a decreasing sequence of non-empty, compact
sets in the spaces $\scrY_{\al_n}$, ($n\geq0$) by
induction, and take the intersection. For now, assume
that $\al_0=\{\scrX\}$. According to
condition~(\ref{property:Ptightprime}), there exists 
a compact set $K_0$ in $\scrY_{\al_0}=\scrX$ such that,
\[
  \Pi_{\al_n}\Bigl( \bigl\{ P_{\al_n}\in M^1(\scrX_{\al_n}):
    P_{\al_n}(\varphi_{\al_n}(K_0)) < 1-\ft12\delta \bigr\} \Bigr)
    < \ft12\ep,
\]
for all $n\geq0$. Make the induction assumption that,
for given $n\geq0$, there is a compact $K_n$ in
$\scrY_{\al_n}$ with,
\begin{equation}
  \label{eq:KnboundOne}
  \Pi_{\al_m}\biggl( \Bigl\{ P_{\al_{m}}\in M^1(\scrX_{\al_{m}}):
          P_{\al_m}(\varphi_{\al_m}(K_n))
            < 1-\frac{2^{n+1}-1}{2^{n+1}}\delta \Bigr\} \biggr)
    < \frac{2^{n+1}-1}{2^{n+1}}\ep,
\end{equation}
for all $m\geq n$. Fix $m\geq n+1$. To combine masses back
at a later stage, choose $0<\lambda_i<1$ for all
$1\leq i\leq {|\al_{n+1}|}$, such that $\sum_i\lambda_i=1$.
For any $1\leq i\leq {|\al_{n+1}|}$,
there exists a compact $K_i\subset A_i$, such that,
\begin{equation}
  \label{eq:KnboundTwo}
  \Pi_{\al_{m}}\Bigl( \bigl\{ P_{\al_{m}}\in M^1(\scrX_{\al_{m}}):
    P_{\al_{m}}(\varphi_{\al_{m}}(A_i)) - P_{\al_{m}}(\varphi_{\al_{m}}(K_i)) >
        2^{-(n+2)}\lambda_i\delta \bigr\} \Bigr)
        < 2^{-(n+2)}\lambda_i\ep.
\end{equation}
for all $m\geq n+1$.
The intersection $K_{n+1}=K_n\cap (\cup_i K_i)$ is not only
compact in $\scrY_{\al_n}$, but also in $\scrY_{\al_{n+1}}$. Then, for
any $m\geq n$, if $P_{\al_m}$ does not lie in any of the
$M^1(\scrX_{\al_m})$-subsets on the left-hand sides of
inequalities~(\ref{eq:KnboundOne}) and~(\ref{eq:KnboundTwo}),
\[
  P_{\al_m}\bigl(\varphi_{\al_m}(K_{n+1})\bigr)
    %\geq P_{\al_m}\bigl(\varphi_{\al_m}(K_{n})\bigr)
    %  - \sum_i P_{\al_m}\bigl(\varphi_{\al_m}(A_i\setminus K_i)\bigr)
    \geq 1-\frac{2^{n+2}-1}{2^{n+2}}\delta,
\]
and the $\Pi_{\al_m}$-probability of that event is lower-bounded by,
\[
  \begin{split}
  \Pi_{\al_{m}}\biggl( \Bigl\{ P_{\al_{m}}\in & M^1(\scrX_{\al_{m}}):
    P_{\al_m}\bigl(\varphi_{\al_m}(K_{n+1})\bigr)
      \geq 1-\frac{2^{n+2}-1}{2^{n+2}}\delta \Bigr\} \biggr)\\
    &\geq 1 - \frac{2^{n+1}-1}{2^{n+1}}\ep - \sum_i 2^{-(n+2)}\lambda_i\ep
    = 1 - \frac{2^{n+2}-1}{2^{n+2}}\ep,
  \end{split}
\]
completing the induction step. Define $\hat{K}=\cap_{n\geq1} K_n$,
which is compact in $\hat{\scrX}$ by lemma~\ref{lem:hatcompact}, and,
\[
  \Pi_{\al_n}\Bigl( \bigl\{ P_{\al_{n}}\in M^1(\scrX_{\al_{n}}):
    P_{\al_n}(\hat\varphi_{\al_n}(\hat{K})) < 1-\delta
      \bigr\} \Bigr) < \ep,
\]
for all $n\geq0$, showing that condition~(\ref{property:Ptight}) is
satisfied, and the assertion follows from theorem~\ref{thm:tightexistence}.
Coming back to the assumption that $\al_0=\{\scrX\}$,
if $\al_0$ consists of more than one set, then the induction argument
is started from a (finite) partition $\al_0$ that coincides with some 
$\al_n$-stage in the proof as provided above.
\end{proof}
The requirements that corollary~\ref{cor:tightexistence} places
on $\scrX$ and $\scrA$ are more specific than those of
theorem~\ref{thm:tightexistence}, but not necessarily more
restrictive: all Polish spaces have countable bases and
well-ordered partition systems $\scrA$, generated by some
countable basis $\scrU$, can all be derived as subsequences
of the generic situation, \cf\
example~\ref{ex:partitionscountablebasis}.

\subsubsection{Coherent random histogram systems without limit}
\label{subsub:nolimit}

To conclude, we consider the cases in which condition~(\ref{property:Ptight})
does not hold. We start with a counterexample that illustrates
concretely how failure of condition~(\ref{property:Ptight}) is
related to the `leaking away' of probability mass in the limit of
refining $\al$.

\begin{example}
\label{ex:coherentbutnoPi}
Consider $\scrX=\RR$ with a basis $\scrU$ defined by all open
intervals with rational midpoints and rational radii. Consider a
triangular array defined by $\{q_{n,m}: n\geq1, 1\leq m\leq M_n=2^n-1\}$
of values in $\QQ$, such that for every $n\geq1$, we have
$q_{n+1,1}\leq q_{n,1}$, $q_{n+1, M_{n+1}}\geq q_{n+1,M_n}$; for every
$1\leq m\leq M_n$, $q_{n+1,2m}=q_{n,m}$; and for every $m\leq M_n-1$,
$q_{n,m}<q_{n+1,2m-1}<q_{n,m+1}$. Defining, $\al_n$ to be of the form,
\[
  \al_n=\bigl\{ A_{n,1}=(-\infty,q_{n,1}], A_{n,2}=(q_{n,1},q_{n,2}],
    \ldots, A_{n,{|\al_n|}}=(q_{n,M_n},\infty) \bigr\},
\]
one verifies that the $\al_n$ are generated by the basis $\scrU$ and
$\al_{n+1}$ refines $\al_n$ for any $n\geq1$. Assuming that the set
$\cup\{q_{n,m}:n\geq1,1\leq m\leq M_n\}$ is dense in $\RR$, the
resulting partitions collectively generate the Borel $\sigma$-algebra.
(For later reference, we indicate the possibility to choose
$q_{n,1}=0$, $q_{n,M_n}=1$, to define partitions on $(0,1]$.)

The simplest example of a coherent histogram system that does
not satisfy condition~(\ref{property:Ptight}), is constructed as
follows. Choose some $\delta>0$, $N\geq1$ and define histogram
distributions $\Pi_{\al_n}$ for all $n\geq N$, for the probability
vectors $(P_{\al_n}(A_{n,m}):1\leq m\leq {|\al_n|})$ satisfying,
\[
  \Pi_{\al_n}\Bigl( P_{\al_n}(A_{n,1})+
    P_{\al_n}(A_{n,{|\al_n|}})=\delta\Bigr)=1,
\]
that is, some non-zero fraction of the total probability mass
in the $n$-th histogram is concentrated in the `outside' sets
$A_{n,1},A_{n,{|\al_n|}}$ with $\Pi_{\al_n}$-probability one. As
$A_{n+1,1}\subset A_{n,1}$ and $A_{n+1,{|\al_{n+1}|}}\subset
A_{n,{|\al_n|}}$, coherence of the histogram system
$(\Pi_{\al_n},\varphi_{*\,\al_n\al_m})$ is maintained.
Assuming that $-q_{n,1},q_{n,M_n}\to\infty$, any compact $K$
in $\scrX=\RR$, fails to meet the `outside' sets,
$K\cap A_{n,1}=\emptyset$ and $K\cap A_{n,{|\al_n|}}=\emptyset$,
for large enough $n$, which invalidates condition~(\ref{property:Ptight}).
\end{example}
To summarize the above example, the problem occurs because the
$\Pi_{\al_n}$ shift a non-zero amount of mass towards
$\pm\infty$ without limitations as $n$ grows. For any presumed
limit measure $\Pi$ on $M^1(\scrX)$, this would mean that
for all compact sets $K$ in $\RR$, $\Pi(P(K)\leq 1-\delta)=1$.
This shows that property~(\ref{property:Ptight}) cannot be
satisfied, and no such limit $\Pi$ exists as a probability
distribution on $M^1(\scrX)$.

Non-compactness of $\scrX$ appears essential in the above
example; however the next example shows that the situation is
more complicated: mass can `leak away' not just to points
at infinity, but at any boundary between partition sets.
\begin{example}
\label{ex:coherentbutnoPigeneric}
In example~\ref{ex:coherentbutnoPi}, take $\scrX$ equal to the
compact subset $[0,1]$, define the points
$s_{n,m}=\ft12+\frac{1}{\pi}\arctan(q_{n,m})$ for all $n\geq1$,
$1\leq m\leq M_n$ and consider partitions,
\[
  \be_n=\bigl\{ B_{n,0}=\{0\},B_{n,1}=(0,s_{n,0}], B_{n,2}=(s_{n,0},s_{n,1}],
    \ldots, B_{n,{|\al_n|}-1}=(s_{n,M_n},1), B_{n,{|\be_n|}}=\{1\}
    \bigr\},
\]
so that ${|\be_n|}={|\al_n|}+2$. Now define the histogram distributions
$\Pi_{\be_n}$ for the probability vectors
$(P_{\be_n}(B_{n,m}):1\leq m\leq {|\be_n|})$, by
\[
  \Pi_{\be_n}\bigl(P_{\be_n}(B_{n,0})=0\bigr)
    =\Pi_{\be_n}\bigl(P_{\be_n}(B_{n,{|\be_n|}})=0\bigr)=1,
\]
and the distribution of $(P_{\be_n}(B_{n,m}):2\leq m\leq {|\be_n|}-1)$
equal to that of $(P_{\al_n}(A_{n,m}):1\leq m\leq {|\al_n|})$. Again, we
have a coherent system of histogram distributions, and like in
example~\ref{ex:coherentbutnoPi}, probability mass is shifted up against
the boundary points at $0$ and $1$, but no limiting distribution $\Pi$
on $M^1([0,1])$ exists with the $\Pi_{\be_n}$ as histogram projections.

In fact, probability mass does not even have to disappear at points of
the boundary of $\scrX$: if we make this example on $[0,1]$ part of a
refining system of partitions of $\RR$, with some fraction of the
total probability mass in (a random histogram system on) the complements
$(-\infty,0)\cup(1,\infty)$, the construction on $\{0\}\cup(0,1)\cup\{1\}$
will continue to make some non-zero fraction of the mass `leak away'
across the boundary of $(0,1)$, which lies in the interior of $\scrX$.
\end{example}

The concluding remark in example~\ref{ex:coherentbutnoPigeneric} is close
to the generic situation: if we partition $\RR$ into intervals, boundaries
between partition sets create the potential for coherent random histogram
systems that make probability mass disappear there in the limit. If we
generalize to higher dimensions, it becomes clear that mass does not
necessarily disappear at specific points, it may be concentrated in
any decreasing sequence of partition sets with empty limit; this shows
in which way a (histogram-specific) form of $\sigma$-additivity makes
a re-appearance.

These counterexamples highlight the significance of
condition~(\ref{property:Ptight}): requiring the existence of a compact
$K$ in $\scrX$ would prevent counterexample~\ref{ex:coherentbutnoPi} but
not~\ref{ex:coherentbutnoPigeneric}. In order to prevent `leakage'
of the latter type, we have to impose the stronger requirement of
a existence of a compact $\hat{K}$ in $\hat{\scrX}$, which keeps
probability mass away from all potential points of `leakage'
simultaneously.

In case condition~(\ref{property:Ptight})
cannot be satisfied, as in examples~\ref{ex:coherentbutnoPi}
and~\ref{ex:coherentbutnoPigeneric}, it is possible to
consider compactification of $\hat{\scrX}$, for example, the
Stone-\v{C}ech compactification $\beta\hat{\scrX}$. With canonical
extension of partitions of $\scrX$ to the space $\beta\hat{\scrX}$
% (note that the indicator for $\beta\hat{\scrX}\setminus\hat{\scrX}$,
% $\hat{x} \mapsto 1-\Sigma\{\hat{1}_A(\hat{x}):A\in\al\}$, is
% continuous (with $\hat{1}_A$ denoting the (unique) continuous
% extension of the indicator $1_A$ to $\beta\hat{\scrX}$)),
condition~(\ref{property:Ptight})
is satisfied trivially. The limiting probability measure $\Pi$
on $M^1(\beta\hat{\scrX})$ may not be unique (because the
projections onto the spaces $M^1(\scrX_\al)$ are not necessarily
separating). Moreover, in applications, the added points in the
closed subset $\beta\hat{\scrX}\setminus\hat{\scrX}$ lack
interpretation.
% we continuously embed $\hat{\scrX}$ as a
% dense subset of a zero-dimensional compact space $\hat{\scrX}'$
% (for example inside a \emph{Cantor cube}, a product
% space of the form $\{0,1\}^I$). In that case
% condition~(\ref{property:Ptight}) is satisfied and a limiting
% Radon probability measure $\hat\Pi'$ exists on $M^1(\hat{\scrX}')$.

\subsection{Support in the tight topology}
\label{sub:suppweak}

Below, theorem~\ref{thm:tightexistence} is used to characterize
the support of histogram limit measures $\Pi$ on $M^1(\scrX)$
with the tight topology, for a Polish space $\scrX$. As it
turns out, the appropriate relation to the mean measure is
inclusion of supports. This assertion was already known
in the literature (see, for example, theorem~4.15 in
\cite{Ghosal17}), but the proof given here was not. In the
formulation of the following theorem, let $G$ denote the
mean measure of definition~\ref{def:meanmeasure}.
\begin{proposition}
\label{prop:supptight} % was {prop:diriweaksupport}
Let $\scrX$ be a Polish space. Consider $M^1(\scrX)$
with the tight topology and a Borel probability
distribution $\Pi$. Let $G$ be the mean measure under $\Pi$.
Then $\{P\in M^1(\scrX):\supp(P)\subset\supp(G)\}$ is
closed in $M^1(\scrX)$ and,
\[
  \supp_{\scrT_C}(\Pi) \subset \{P\in M^1(\scrX):\supp(P)\subset\supp(G)\}.
\]
Moreover, if $P\in M^1(\scrX)$ is such that for all
partitions $\al\in\scrA$, $P_\al$ lies in
the support of $\Pi_\al$ in $M^1(\scrX_\al)$, then $P$
lies in the tight support of $\Pi$.
\end{proposition}
\begin{proof}
If $P$ is such that $\supp(P)\not\subset\supp(G)$, there
exist an $x\in\supp(P)\setminus\supp(G)$ and, by complete
regularity of $\scrX$, a continuous $f:\scrX\to[0,1]$ with $f=0$
on $\supp(G)$ and $f(x)=1$.
% (see definition~\ref{def:completelyregularspace}).
While $\langle G, f\rangle=0$, the open neighbourhood of $x$ for which
$f>\ft12$ receives non-zero $P$-probability, and we see that
$\langle P, f\rangle>0$. So if $Q$ lies in the tight neighbourhood
$\{Q\in M^1(\scrX):|\langle(P-Q),f\rangle|<\ft12\ep\}$ of $P$
(for some $0<\ep<\langle P, f\rangle$), $\langle Q, f\rangle>0$ and
accordingly, $\supp(Q)\not\subset\supp(G)$, from which it follows
that $\{P\in M^1(\scrX):\supp(P)\subset\supp(G)\}$ is closed.
Moreover, by Markov's inequality and Fubini's theorem,
\[
  \begin{split}
  \Pi\bigl( \{Q\in &M^1(\scrX):|\langle(P-Q),f\rangle|<\ft12\ep
      \}\bigr)\\
    &\leq \Pi\bigl( \{Q\in M^1(\scrX): \langle Q, f\rangle>\ft12\ep
      \}\bigr)
    \leq \frac{2}{\ep}\int \langle Q, f\rangle \dd\Pi(Q)
      = \frac{2\langle G,f\rangle}{\ep} = 0.
  \end{split}
\]
Conclude that $P$ has a tight neighbourhood of $\Pi$-mass zero,
which means that $P$ does not lie in the tight support of $\Pi$.

Regarding the last assertion, it is noted that, since
$M^1(\scrX)$ with the tight topology is the continuous
image of a subset of the inverse limit $N$ of
proposition~\ref{prop:invsysmsr}, the collection of sets
$\{\varphi_{*\,\al}^{-1}(V):\al\in\scrA,V\in\scrU_\al\}$
(where $\scrU_\al$ is any basis for $M^1(\scrX_\al)$,
\eg\ total-variational balls) in $M^1(\scrX)$ forms a
basis for the tight topology. Consequently, for any
tight neighbourhood $U$ of $P\in M^1(\scrX)$ there exists
an $\al\in\scrA$ and a $V\in\scrU_\al$ such that
$\varphi_{*\,\al}^{-1}(V)\subset U$, and,
\[
  \Pi(U)\geq\Pi(\varphi_{*\,\al}^{-1}(V))=\Pi_\al(V)>0,
\]
by assumption.
\end{proof}

%%%%%%%%%%%%%%%%%%%%%%%%%%%%%%%%%%%%%%%%%%%%%%%%%%%%%%%%%%%%%%%%%%%%%%%%%%%%%%%
\section{Phase structure of probability histogram limits}
\label{sec:phasestructure}

In this section we combine the two main existence theorems of preceding
sections with the general theory of completely random measures
\citep{Kingman67}, to describe the various ways in which random
histogram limits manifest. In subsection~\ref{sub:completelyrandom}
we review completely random point processes \citep{daley07},
and show in subsection~\ref{sub:phases} how combination leads to the
conclusion that random histogram limits occur in one of four distinct
\emph{phases}: continuous-singular or dominated, each either purely
atomic or not (see theorem~\ref{thm:phases} below). The phase of a
random histogram limit depends on the topology on $M^1(\scrX)$ and
on independence within random histogram distributions. In
sections~\ref{sec:polya} and~\ref{sec:gauss}, we demonstrate that,
both in the P\'olya-tree family and in the Gaussian family of
histogram limits, changes in the defining parameters of their histogram
systems can cause the limit to transition from one phase to another.

\subsection{Completely random measures}
\label{sub:completelyrandom}

In \citep{Kingman67,Kingman75} so-called \emph{completely random measures}
are defined as positive random measures $\nu\sim\Pi'$ that assign
stochastically independent random masses to disjoint measurable subsets
of the underlying space $\scrX$, and it is shown that (the
random part of) a completely random measure is a \emph{purely atomic}
measure with $\Pi'$-probability one. (Note, we say that a positive
measure $\nu$ is \emph{purely atomic}, if the collection $D$ of points
$x\in\scrX$ for which $\nu(\{x\})>0$ (so-called \emph{atoms}),
contains all $\nu$-mass, \ie\
$\nu(D)=\nu(\scrX)$; we say that $\nu$ is \emph{non-atomic},
if $D=\emptyset$.) Below we give the briefest of introductions to
completely random measures (following \cite[chapters~9 and~10]{daley07}),
and relate the results to the existence theorems of
sections~\ref{sec:weakhistogramlimits} and~\ref{sec:prokhorov}.

\begin{definition}
\label{def:completelyrandommeasure}
Let $(\scrX,\scrB)$ be a Polish space. A random positive Radon
measure $\nu$ on $\scrX$, distributed according to $\Pi'$,
is called a
%\varinindex{completely random measure}{measure!completely random}
\emph{completely random measure}, if, for any finite collection
of disjoint measurable sets $A_1,\ldots,A_n\in\scrB$, the measures
$\nu(A_1),\ldots,\nu(A_n)$ are independent.
%If, for some
%completely random measure $\nu\sim\Pi'$, $\Pi'(\nu(\scrX)<\infty)=1$,
%then $\nu/\nu(\scrX)$ is a random probability measure
%called a \varinindex{normalized completely random measure}%
%{measure!normalized completely random}.
\end{definition}
Any (random) positive Radon measure $\nu\sim\Pi'$ decomposes
as a sum of a (random) purely atomic measure $\nu_d$ and a
(random) non-atomic measure $\nu_n$ in a unique way
\citep[Proposition~9.3.IV]{daley07},
\begin{equation}
  \label{eq:atomicnonatomic}
  \Pi'\bigl( \nu = \nu_n +\nu_d \bigr) =1,
\end{equation}
and for a random positive Radon measure to be almost-surely
non-atomic, it is necessary and sufficient that for any
$\ep,\delta>0$, there is a finite Borel
measurable partition $\al$ of $\scrX$ such that for all
finer finite Borel measurable partitions $\be$,
$\Pi'\bigl( \max\{\nu(B):B\in\be \}>\delta \bigr)<\ep$.
In the case of a \emph{completely} random measure $\nu$,
this implies that $\nu$ is $\Pi'$-almost-surely equal to
some fixed (that is, non-random) non-atomic measure $\nu_n$
\citep[Proposition~10.1.II]{daley07}. As a consequence,
the atomic part of any completely random measure
can be fixed or random, while the non-atomic part is
always non-random.
\begin{definition}
\label{def:cumulant}
For any random positive Radon measure $\nu\sim\Pi'$ on
$(\scrX,\scrB)$ and any $t>0$, define the
%\varinindex{cumulant}{measure!cumulant}
\emph{cumulant} $\lambda_t:\scrB\to[0,\infty]$, by,
\[
  \lambda_t(A) = \log\int e^{t\nu(A)}\dd\Pi'(\nu).
\]
\end{definition}
\emph{Fubini's theorem} implies that, if $\nu\sim\Pi'$ is a
completely random measure, for any $t>0$, the cumulant
$\lambda_t$ is a positive Borel measure.
The theorem below says that the atomic part
of a completely random measure decomposes
into a sum of random atoms at \emph{fixed points}
in $\scrX$, and a sum $\nu_r$ of random atoms at
\emph{random points} in $\scrX$.
\begin{theorem}
\label{thm:Kingman} \citep{Kingman67}\\
Let $\nu\sim\Pi'$ be a completely random measure with cumulant
measures $\lambda_t$ for $t>0$. If all $\lambda_t$ are
$\sigma$-finite, then with $\Pi'$-probability one, $\nu$
satisfies the decomposition,
\begin{equation}
  \label{eq:KingmanCRM}
  \nu = \nu_n + \nu_f + \nu_{r},
\end{equation}
where $\nu_n$ is a non-random, non-atomic, $\sigma$-finite
measure on $(\scrX,\scrB)$; $\nu_f$ is a purely atomic measure
supported on a fixed, countable subset $D\subset\scrX$
where $\nu_f(\{x\})$ and $\nu_f(\{x'\})$ are independent
if $x,x'\in D$, $x\neq x'$; and $\nu_{r}$ is a random
purely atomic measure that is independent of $\nu_f$.
\end{theorem}
As it turns out, $\sigma$-finiteness of $\lambda_t$ is equivalent
to the existence of countable cover $C_1,C_2,\ldots\in\scrB$
of $\scrX$ such that $\Pi'(\nu(C_i)<\infty)>0$, for every $i\geq1$.
Furthermore, the set of fixed atoms $D$ is the set of atoms of
$\lambda_t$ and $\sigma$-additivity of $\lambda_t$ implies
countability of $D$. The random purely atomic measure
$\nu_{r}$ is realized with the help of a Poisson point-process
$N$ on $\scrX\times(0,\infty)$ (\cf\
\cite[Proposition~9.1.III-(v)]{daley07}), as follows:
\[
  \nu_r(A)=\int y\,N(A\times\dd y),
\]
with an intensity measure $\mu$ that may be unbounded on sets
of the form $A\times(0,\ep)$, $(\ep>0)$, but satisfies,
\[
  \int \min\{1,y\}\dd \mu(A\times\dd y)<\infty. 
\]
The measure $\nu_n$ appears in $\lambda_t$ as the $t$-linear
contribution: $\lambda_t(A)=t\nu_n(A)+\ldots$. 
So completely random measures with cumulant measures without
$t\nu_n$-terms ($\nu_n=0$) and without fixed atoms ($\nu_f=0$)
are characterized as purely atomic with \emph{random}
locations, purely $\nu_r$; and similarly, completely random
measures with $\nu_n=0$ and Poisson intensity measure $\mu=0$,
are characterized as purely atomic with \emph{fixed}
locations, purely $\nu_f$.

Complete randomness imposes a purely atomic nature on random
probability measures too, after normalization: if a given positive
random measure $\nu\sim\Pi'$ satisfies $0<\nu(\scrX)<\infty$ with
$\Pi'$-probability one, then $P=\nu/\nu(\scrX)\sim\Pi$ defines
a random probability measure called a \emph{normalized} completely
random measure, and $P$ inherits the purely atomic nature of $\nu$.
The histogram distributions $\Pi_\al$ for $P_\al$ follow from
the distributions $\Pi'_\al$ through,
\begin{equation}
  \label{eq:NCRM}
  \bigl(P_\al(A_{\al,1}),\ldots,P_\al(A_{\al,{|\al|}})\bigr)
    =\frac1{\nu_\al(\scrX)}\bigl( \nu_\al(A_{\al,1}),
      \ldots,\nu_\al(A_{\al,{|\al|}}) \bigr)\sim\Pi_\al,
\end{equation}
where $\nu_\al(\scrX)=\nu_\al(A_{\al,1})+\ldots
+\nu_\al(A_{\al,{|\al|}})$. We say that the random histograms
$P_\al$ are independent up to normalization.

\subsection{Phases of probability histogram limits}
\label{sub:phases}

Combining the conclusions of sections~\ref{sec:weakhistogramlimits}
and~\ref{sec:prokhorov} with the presence or absence of complete
randomness, we arrive at the following theorem. (Given a random
probability measure, let $G$ denote the associated mean
measure.)
\begin{theorem}
\label{thm:phases}
{\it (Phases of random histogram limits)}\\
Let $\scrX$, $\scrA$ and $M^1(\scrX)$ satisfy the minimal
conditions. Let $\scrA$ be a directed set of finite, Borel
measurable partitions that resolves $\scrX$, with a coherent
system of Borel histogram probability measures
$(\Pi_\al,\varphi_{*\,\al\be})$ on the inverse system
$(M^1(\scrX_\al),\varphi_{*\,\al\be})$.
\begin{itemize}
  \item[(i.)] \emph{(absolutely-continuous)}\\%{\it (absolutely-continuous)}\\
  If condition~(\ref{property:Pws}) is satisfied, the histogram
  limit describes a random element $P$ of $M^1(\scrX)$,
  distributed according to a weakly-Radon probability
  measure $\Pi$, such that $\Pi$-almost-surely, for all
  measurable $B\in\scrB$, $G(B)=0$ implies $P(B)=0$:
  \[
    \Pi\bigl(\{P\in M^1(\scrX):P\ll G\}\bigr)=1.
  \]
  The random element $P$ can be identified isometrically with
  a random positive Radon-Nikodym density function $p$ in
  $L^1(\scrX,\scrB,G)$ of norm one, and we can write, for all
  $B\in\scrB$,
  \[
      P(B) = \int_B p(x) \dd G(x);
  \]
  \item[(ii.)] \emph{(fixed-atomic)}\\%{\it (fixed-discrete)}\\
  if condition~(\ref{property:Pws}) is satisfied and the $\Pi_\al$
  describe normalized completely random histograms, \cf\
  (\ref{eq:NCRM}), the histogram limit $P\sim\Pi$ is a normalized
  version of the sum $\nu$ of a fixed non-atomic measure $\nu_n\ll G$
  and a random purely atomic measure $\nu_f$ supported on
  the fixed, countable set $D=\{x\in\scrX:G(\{x\})>0\}$. For all
  $B\in\scrB$,
  \[
    P(B) = \frac1{\nu(\scrX)}\bigl(\nu_n(B)+\nu_f(B)\bigr).
  \]
\end{itemize}
Assume, in addition, that $\scrX$ is a Polish space and that $\scrA$
is a directed set of finite partitions generated by a basis that
resolves $\scrX$.
\begin{itemize}
  \item[(iii.)] \emph{(continuous-singular)}\\%{\it (continuous-singular)}\\
  If condition~(\ref{property:Ptight}) is satisfied, the histogram
  limit describes a random element $P$ of $M^1(\scrX)$,
  distributed according to a tightly-Radon probability
  measure $\Pi$, such that $\Pi$-almost-surely, for all open
  $U\subset\scrX$, $G(U)=0$ implies $P(U)=0$:
  \[
    \Pi\bigl(\{P\in M^1(\scrX):\supp(P)\subset\supp(G)\}\bigr)=1;
  \]
\end{itemize}
\begin{itemize}
  \item[(iv.)] \emph{(random-atomic)}\\%{\it (random-discrete)}\\
  if condition~(\ref{property:Ptight}) is satisfied and the $\Pi_\al$
  describe normalized completely random histograms, \cf\
  (\ref{eq:NCRM}), the histogram limit $P\sim\Pi$ is a normalized
  version of the sum $\nu$ of a fixed non-atomic measure $\nu_n\ll G$,
  a random purely atomic measure $\nu_f$ supported on
  the fixed, countable set $D=\{x\in\scrX:G(\{x\})>0\}$, and
  a random purely atomic measure $\nu_r$. For all $B\in\scrB$,
  \[
    P(B) = \frac1{\nu(\scrX)}\bigl(\nu_n(B)+\nu_f(B)+\nu_r(B)\bigr).
  \]
\end{itemize}
\end{theorem}
\begin{proof}
In cases {\it (i.)} and {\it (iii.)}, the theorem states the assertions
of theorems~\ref{thm:WeakExistence} and~\ref{thm:tightexistence};
in cases {\it (ii.)} and {\it (iv.)}, these assertions are combined
with those of theorem~\ref{thm:Kingman}, specific to normalized
completely random measures (where it is observed that 
the set of atoms of $\lambda_t$ (for any $t>0$), is equal to the
set of atoms of $G$).
\end{proof}
In qualitative terms, we may describe the phase structure of
random histogram limits as follows: the most general, least
constrained type of limit above, is that of the \emph{continuous-singular
phase}. According to (\ref{eq:atomicnonatomic}), any
continuous-singular random $P$ decomposes into a random atomic
component and a random non-atomic component.
The random component of any completely random case manifests as
purely atomic (the \emph{random-atomic phase}), with independent,
randomly-sized point masses at fixed locations and at independent
random locations. Many examples of (normalized)
completely random families are known, including the
well-known Dirichlet family (which is discussed in
section~\ref{sec:dirichletlimits}), and a
sub-family of Gaussian histogram systems (see
section~\ref{sec:gauss}).

The random non-atomic component of a histogram limit in the 
continuous-singular phase is novel and more interesting: it
is implied by the above that \emph{dependence in histogram
distributions is required to induce a random non-atomic
continuous-singular component}. To illustrate the nature
of such a component, we may think, for example, of
$\scrX=[0,1]$ and a $\Pi$ with $G$ equal to Lebesgue
measure, describing a random Stieltjes function
$F:[0,1]\to[0,1]$, from a class that is \emph{everywhere
continuous but not everywhere (or even nowhere)
differentiable} (\eg\ the so-called \emph{Cantor
distribution}). Such distributions are non-atomic
but cannot be identified with random Radon-Nikodym density
functions. The Gaussian histogram systems of section~\ref{sec:gauss}
are in the non-atomic continuous-singular phase generically,
and only Gaussian histogram systems with diagonal covariance
matrices are in the random-atomic phase.

In the \emph{absolutely-continuous phase} the histogram
distributions are such that the histogram probabilities
$P_\al(A)$ may be larger than their means $G_\al(A)$, but
not to such a degree that ($G_\al$-averages of) proportions
between $P_\al$ and $G_\al$ grow unbounded in the limit.
This is borne out by the formulation of
property~(\ref{property:Pweakprime}),
and also serves to interpret later bounds (\eg\
(\ref{eq:boundforws})). The upper bound on
the proportions between $P_\al(A)$ and $G_\al(A)$ induces
domination $P\ll G$ with $\Pi$-probability one. Extending
the above example with $\scrX=[0,1]$, the absolutely
continuous phase describes a random
Stieltjes function $F:[0,1]\to[0,1]$ which is
\emph{everywhere differentiable}, and can be identified
with a random Radon-Nikodym density function with respect
to $G$. In section~\ref{sec:gauss} we discuss
Gaussian random histogram limits in the 
absolutely-continuous phase. If we specify that an
absolutely-continuous random histogram limit is also
normalized completely random, then the limit is in the
\emph{fixed-atomic phase}: combining the resulting purely
atomic character of the random component with domination
by $G$, we find only random point masses at the
\emph{fixed} locations of the atoms of $G$. The sub-family
of Dirichlet process distributions with countably supported
base measures are in the fixed-atomic phase.

The distinction between the random-atomic and fixed-atomic phases
provides an alternative explanation for the decomposition
$\nu_f+\nu_r$ of the random purely atomic component in Kingman's
theorem: based on the above and the Radon-Nikodym theorem,
we explain this by the fact that any random probability measure
decomposes uniquely into a random component dominated by its
mean measure $G$, and a random component that is mutually
singular with respect to $G$ (but still with support inside
the support of $G$).

\section{Existence and phases of Dirichlet histogram limits}
\label{sec:dirichletlimits} % was {sub:dirichlettightlimits} 

The best-known family of histogram limits is the Dirichlet
family; its definition is based most conveniently on the
observation that, if $Z_1,\ldots,Z_n$ are independent and
distributed according to Gamma-distributions
$\Gamma(\nu_1,1),\ldots,\Gamma(\nu_n,1)$,
then $Z_1+\ldots+Z_n$ is distributed according to
$\Gamma(\nu_1+\ldots+\nu_n,1)$. (Below, we use convention
that $\Gamma(0,1)$ is a single atom of mass one located at zero.)
\begin{definition}
Let $\nu$ be a non-zero, bounded, positive Borel measure on
a Polish space $\scrX$ and define, for every Borel measurable
partition $\al$,
\begin{equation}
  \label{eq:dirichlethistograms}
  (P(A):A\in\al) \sim \mathrm{Dir}_{\nu_\al},
\end{equation}
where $\nu_\al=(\nu(A):A\in\al)$. The histogram distributions
$\mathrm{Dir}_{\nu_\al}$ on $M^1(\scrX_\al)$ are those of
the normalized positive random elements
$(Z_1/S,\ldots,Z_{|\al|}/S)$, where $Z_i\sim\Gamma(\nu(A_i),1)$,
($i\in I(\al)$) and $S=\sum_iZ_i$. 
Together, the distributions $(\mathrm{Dir}_{\nu_\al},\al\in\scrA_0)$
are coherent and form the
\emph{Dirichlet histogram system with base measure $\nu$}.
\end{definition}
It is clear that the \emph{Gamma-process}, defined by the positive
random vectors $(Z_1,\ldots,Z_n)\sim\prod_i\Gamma(\nu(A_i),1)$
is completely random and that Dirichlet histogram systems are
\emph{normalized completely random}. Limits of Dirichlet
histogram systems therefore describe random probability measures
in one of the two atomic phases.

A second immediate observation,
is that coherence of the histogram system could have been
guaranteed based on parametrization in terms of a
\emph{finitely additive} base measure $\nu$.
The well-known \emph{mean measure condition}
\citep{Orbanz11} requires $\nu$ to be \emph{countably additive}
to guarantee existence of a unique histogram limit
with respect to the tight topology on $M^1(\scrX)$. We come back
to the mean measure condition below.

\subsection{Tight limits of Dirichlet histogram systems}
\label{sub:tightdirichletlimits}

The following theorem is
the (by now, classical) existence result for Dirichlet histogram
limits, with a new proof in terms of
condition~(\ref{property:Ptight}).
\begin{theorem}
\label{thm:existencediri}
%\varinindex{(Existence of Dirichlet process distributions)}%
%{existence!Dirichlet process distribution}\\
Let $\scrX$ be a Polish space, endow $M^1(\scrX)$ with
the tight topology and let $\nu$ be a non-zero, bounded,
positive Borel measure on $\scrX$. There exists a unique
Radon probability measure $\mathrm{Dir}_{\nu}$ on
$M^1(\scrX)$ projecting to the Dirichlet histogram
distributions (\ref{eq:dirichlethistograms}), describing
a random probability measure in the random atomic phase.
\end{theorem}
\begin{proof}
Let $\scrU$ be a countable basis for $\scrX$ and let $\scrA$ be
a refining sequence of partitions, generated by $\scrU$, that
resolves $\scrX$. By assumption there exist distributions
$\mathrm{Dir}_{\nu_\al}$ for the random histograms
$P_\al\in M^1(\scrX_\al)$, ($\al\in\scrA$). As said, coherence
of the inverse system $(\mathrm{Dir}_{\nu_\al},\varphi_{*\,\al\be})$
follows from finite additivity of the measure $\nu$.

To prove condition~(\ref{property:Ptight}), let $\ep>0, \delta>0$
be given. According to proposition~\ref{prop:msrstaysmsr}, $\nu$
defines a bounded positive Borel measure on $\hat{\scrX}$, and
according to proposition~\ref{prop:concretehatX}, $\hat{\scrX}$
is Polish, so $\nu$ is a Radon measure on $\hat{\scrX}$. Hence
there exists a compact $\hat{K}$ in $\hat{\scrX}$ such that,
\[
  \nu(\hat{\scrX}\setminus \hat{K}) < \delta\ep\,\nu(\hat{\scrX}).
\]
Let $\al$ be given. By Markov's inequality and the fact that
under $\mathrm{Dir}_{\nu_\al}$,
\begin{equation}
  \label{eq:dirichletisbeta}
  P_\al(A)\sim\mathrm{Beta}\bigl(\nu(A),\nu(\hat{\scrX})-\nu(A)\bigr),
\end{equation}
for any $A\in\sigma(\al)$,
%(see example~\ref{ex:dirfortwo})
we have,
\[
  \begin{split}
  \mathrm{Dir}_{\nu_\al}&\Bigl( \bigl\{P_\al\in M^1(\scrX_\al):
    P_\al(\hat{\varphi}_\al(\hat{K})) < 1 - \delta \bigr\}
      \Bigr)\\[.5em]
  &= \mathrm{Dir}_{\nu_\al}\Bigl( \bigl\{P_\al\in M^1(\scrX_\al):
    P_\al(\scrX_\al) - P_\al(\hat{\varphi}_\al(\hat{K})) > \delta \bigr\}
      \Bigr)\\[.5em]
  &\leq \frac1\delta \int_{M^1(\scrX_\al)}
    P_\al (\scrX_\al\setminus\hat{\varphi}_\al(\hat{K}) )\dd
      \mathrm{Dir}_{\nu_\al}(P_\al)\\[.5em]
  &= \frac1\delta \frac{\nu(\hat{\scrX} \setminus
    (\hat{\varphi}_\al^{-1}\circ\hat{\varphi}_\al)(\hat{K}))}{\nu(\hat{\scrX})}
  \leq \frac1\delta \frac{\nu(\hat{\scrX}\setminus\hat{K})}{\nu(\hat{\scrX})}
    < \ep,
  \end{split}
\]
by Markov's inequality, the fact that the $G_\al$
are proportional to $\nu_\al$, and the fact that
$\hat{K}\subset
(\hat{\varphi}_\al^{-1}\circ\hat{\varphi}_\al)(\hat{K})$.
Conclude that there exists a unique histogram limit
$\mathrm{Dir}_{\nu}$, a Radon probability measure on $M^1(\scrX)$
with the tight topology. Because the histogram system is
normalized completely random, the limiting random element $P$
is in the random-atomic phase.
\end{proof}
To conclude, two remarks are in order: firstly coming back
to the \emph{mean measure condition}, it is noted
that the above proof relies on $\nu$ being not just finitely,
but countably additive, to imply the Radon property.
Secondly, we note that restriction to $\scrA$ with partitions
generated by the basis may be confusing, since the most common
definition of the Dirichlet histogram system involves
all Borel measurable partitions, $\scrA_0$. We argue this
distinction expresses the difference between the roles that
$\scrA$ plays in theorem~\ref{thm:existencediri} and
proposition~\ref{thm:randomhistogram}: to
define $\mathrm{Dir}_{\nu}$, we restrict to directed sets
$\scrA$ of a special form, while after proving existence,
we may use histograms associated to all $\al\in\scrA_0$.

\subsection{Weak limits of Dirichlet histogram systems}
\label{sub:weakdirichletlimits}

Whether $\mathrm{Dir}_\nu$ is a Radon measure with respect to
the weak topology as well, depends on the base measure $\nu$. To
make a preliminary assessment, note that, given $\al\in\scrA$
and $L>0$, for any $P_\al,Q_\al\in M^1(\scrX_\al)$,
\begin{equation}
  \label{eq:boundforws}
  \begin{split}
  \|P_{\al}-&P_{\al}\wedge L\,Q_{\al}\|_{1,\scrX_{\al}}
    =\sum \bigl\{ P_{\al}(A_i)
      \,:\, i\in I(\al),\,P_{\al}(A_i)>L\,Q_{\al}(A_i)\bigr\}\\[.5em]
    &\leq \frac{1}{L} \sum \biggl\{
      \biggl(\frac{P_{\al}(A_i)}{Q_{\al}(A_i)}\biggr)^2\,Q_{\al}(A_i)
        \,:\, i\in I(\al),\,P_{\al}(A_i)>L\,Q_{\al}(A_i)\biggr\}\\[.5em]
    &\leq \frac{1}{L} \sum_{i\in I(\al)} \frac{P_{\al}(A_i)^2}{Q_{\al}(A_i)}.
  \end{split}
\end{equation}
Based on (\ref{eq:dirichletisbeta}), we see that, for every $A\in\al$,
\[
  \int_{M^1(\scrX_\al)} P_\al (A)^2 \dd\mathrm{Dir}_{\nu_\al}(P_\al)
    %= \frac{\nu(A_i)(\nu(\scrX)-\nu(A_i))}{\nu(\scrX)^2(\nu(\scrX)+1)}
    %  + \frac{\nu(A_i)^2}{\nu(\scrX)^2}
    = \frac{\nu(A)^2+\nu(A)}{\nu(\scrX)^2+\nu(\scrX)}.
\]
Now, let $\delta>0$ be given. Due to the bound (\ref{eq:boundforws}),
for any $L>0$ and any $\al\in\scrA$, Markov's inequality gives,
\[
  \begin{split}
  \Pi_{\al}\bigl(\{ &P_{\al}\in M^1(\scrX_{\al}):
    \|P_{\al}-P_{\al}\wedge LG_{\al}\|_{1,\scrX_{\al}}
    > \delta \}\bigr)\phantom{\biggm|}\\
  &\leq
  \Pi_{\al}\biggl(\Bigl\{ P_{\al}\in M^1(\scrX_{\al}):
    \frac1{L}\sum_{i\in I(\al)} \frac{P_{\al}(A_i)^2}{G_{\al}(A_i)}
    > \delta \Bigr\}\biggr)\\
  &\leq \frac1{L\delta} \sum_{i\in I(\al)} \frac1{G_{\al}(A_i)}
    \int_{M^1(\scrX_{\al})}P_{\al}(A_i)^2\dd\Pi_{\al}(P_{\al})\\
  &= \frac1{L\delta} \sum_{i\in I(\al)}
    \frac{\nu(A_i)+1}{\nu(\scrX)+1}=\frac1{L\delta}
    \frac{\nu(\scrX)+{|\al|}}{\nu(\scrX)+1}
  \end{split}
\]
for all $\al\in\scrA$. Since ${|\al|}\to\infty$, as $\al\in\scrA$ refines
(unless $\scrX$ is finite), this shows that the most obvious upper
bound to imply uniform integrability does not lead to a useful argument.
However, we show the following.
\begin{theorem}
\label{thm:weakdirichlet}
Let $\scrX$, $\scrA=\scrA_0$ and $M^1(\scrX)$ satisfy the minimal
conditions and consider $M^1(\scrX)$ with the weak topology. 
Let $\nu$ be a non-zero, bounded, positive, purely atomic
measure on $\scrX$. Then there exists a unique
Radon probability measure $\mathrm{Dir}_{\nu}$ on $M^1(\scrX)$
with the weak topology, projecting to $\mathrm{Dir}_{\nu_\al}$
for all $\al\in\scrA_0$. In that case, $\mathrm{Dir}_{\nu}$
describes a normalized completely random measure in the
fixed-atomic phase. 
\end{theorem}
\begin{proof}
First consider a \emph{countable set} $D$ with the
\emph{discrete topology} (which is a Polish space), with a
bounded, positive Borel measure $\nu_D$ on $D$. According
to theorem~\ref{thm:existencediri}, the Dirichlet
histogram system with base measure $\nu_D$ has a Radon
histogram limit $\mathrm{Dir}_{\nu_D}$ on $M^1(D)$ with
the tight topology. Since any bounded $f:D\to\RR$ is continuous,
the tight and weak topologies are equal. Therefore
$\mathrm{Dir}_{\nu_D}$ is also Radon with respect to the
weak topology on $M^1(D)$ by default.

Now, let $\scrX$ be Polish, and
let $D$ denote the set $\{x\in\scrX:\nu(\{x\})>0\}$.
Let $\scrA_D$ denote the set of all finite partitions of $D$,
and let $\scrA_{\scrX\setminus D}$ denote the set of all finite,
Borel measurable partitions of $\scrX\setminus D$. Define $\scrA$
to contain all partitions $\al$ that combine a partition $\al_D$
from $\scrA_D$ and a partition $\al_{\scrX\setminus D}$ from
$\scrA_{\scrX\setminus D}$, to partition the whole space $\scrX$.
Note that $\scrA$ resolves $\scrX$, and $\scrA$ is directed
and co-final in $\scrA_0$. For any $\al
=(\al_D,\al_{\scrX\setminus D})\in\scrA$, the Dirichlet
histogram distribution $\mathrm{Dir}_{\nu_\al}$ is such that,
for the ($\sigma(\al)$-measurable) subset $\scrX\setminus D$,
\[
  \mathrm{Dir}_{\nu_\al}\bigl(
  \{P_\al\in M^1(\scrX_\al):P_\al(\scrX\setminus D)=0\} \bigr) = 1.
\]
so $P_\al(D)=1$ with $\mathrm{Dir}_{\nu_\al}$-probability one.
The projections of $(P(A):A\in\al)$ onto $(P(A):A\in\al_D)$
give rise to a Dirichlet histogram system with base measure
$\nu_D$, the restriction of $\nu$ to subsets of $D$. As argued
above, the limit $\mathrm{Dir}_{\nu_D}$ is a Radon probability
measure on $M^1(D)$ with the weak topology. The space $M^1(D)$
is weak-to-weak homeomorphic to the weakly closed subspace
$M$ of all $P\in M^1(\scrX)$ such that $P\ll \nu$, through
the mapping $\phi:M^1(D)\to M^1(\scrX)$, $\phi(P)(B)=P(B\cap D)$
for all $B\in\scrB$. Conclude that the histogram system based
on partitions in $\scrA$ has a histogram limit
$\mathrm{Dir}_{\nu}$ that is Radon on $M^1(\scrX)$ with the
weak topology.
%[DO WORK HERE] Conversely ...
\end{proof}

%%%%%%%%%%%%%%%%%%%%%%%%%%%%%%%%%%%%%%%%%%%%%%%%%%%%%%%%%%%%%%%%%%%%%%%%%%%%%%%

\section{Existence and phases of P\'olya-tree histogram limits}
\label{sec:polya}

Here we give only a very brief introduction to
P\'olya-tree distributions, for much
more see \citep{Kraft64,Mauldin92,Lavine92,Lavine94} and the
overviews in \citep{Ghosh03,Ghosal17}.

The P\'olya-tree distribution is defined through a sequence
of refining partitions of a Polish space $\scrX$ (usually $\RR$
or the interval $[0,1]$), where in each step, every set in the preceding
partition is split in two subsets. To describe the resulting
tree of refinements, we define the following: for every $m\geq0$,
we denote by $\scrE_m$ the set of all binary sequences
$\vep$ of length $m$ (and we denote the empty binary
sequence formally as $\vep_\emptyset$, forming the only
element of the set denoted $\scrE_0$). We also define the
set $\scrE=\cup_{m\geq0}\scrE_m$ of all \emph{finite binary
sequences} (including the empty one). For any two binary
sequences $\vep\in\scrE_m$, $\vep'\in\scrE_{m'}$, we
write $\vep\vep'$ for the concatenation in $\scrE_{m+m'}$.
In particular, for any $\vep\in\scrE_m$, $\vep0$
($\vep1$) in $\scrE_{m+1}$ appends a zero (one) to $\vep$.
Also note that $\vep_\emptyset\vep=\vep\vep_\emptyset=\vep$
for all $\vep\in\scrE$. We write out $\vep\in\scrE_m$
as $\vep=e_1\ldots e_m$, and use the notation
$\vep_l:=e_1\ldots e_l\in\scrE_{l}$ for the projections
onto the first $1\leq l\leq m$ binary digits.
We also define, for any $\vep\in\scrE_m$ with $e_m\in\{0,1\}$,
$\vep$ with the last digit flipped: $\hat\vep=\vep_{m-1}(\neg e_m)$.

We use $\scrE$ to organise a refining sequence
$\scrA=\{\al_n:n\geq0\}$ of partitions,
$\al_0=\{\scrX\}$, $\al_1=\{A_0,A_1\}$,
$\al_2=\{A_{00},A_{01},A_{10},A_{11}\}$, \etc, into
a \emph{dyadic tree}, defining $\al_n=\{ A_\vep:\vep\in\scrE_n\}$
and for all $\vep\in\scrE$,
\begin{equation}
  \label{eq:splitA}
  A_\vep = A_{\vep0}\cup A_{\vep1}.
\end{equation}
Mostly we shall
look at refinement through intersection with basis sets and
their complements, \ie\ for every $\vep\in\scrE$ either
$A_{\vep0}$ or $A_{\vep1}$ equals $A_\vep \cap U$ for some
element $U$ in a basis $\scrU$ for $\scrX$. Note that in
the case of a countable basis $\scrU$, iterative
application of the above construction gives rise to a
countable $\scrA=\{\al_m:m\geq1\}$ that resolves $\scrX$.
% Indeed, \emph{if} a random probability measure $P\sim\Pi$
% as in theorem~\ref{thm:randomhistogram} \emph{exists},
% then random histograms as in (\ref{eq:randomhistogram}) are
% defined for all Borel measurable partitions. However, to
% define an inverse system and prove existence, smaller families of
% partitions are possible and, as in this case of the P\'olya
% tree random histogram systems, more practical. % was {ex:halfhalf}
\begin{example} 
\label{ex:dyadictree}
A typical example is a \emph{dyadic tree} of partitions of
$\scrX=(0,1]$ (or $[0,1]$), constructed by iteratively
bisecting every interval at the mid-point.
This leads to a sequence of refining partitions $\al_m$,
$m\geq0$, consisting of $2^m$ intervals of the forms $(l,u]$
where $l=u-2^{-m}$ and $u=2^{-m}k$, $k=1,2\ldots,2^m$, which
is generated by a basis and which resolves $(0,1]$. (In case
$\scrX=[0,1]$ we add to every partition the singleton $\{0\}$.)
\end{example}
% \begin{example}
% \label{ex:dyadictreeRR}
% We also specify a dyadic tree of partitions of $\RR$.
% Let $\scrE=\cup_{m\geq0}\scrE_m$ denote the set of all
% \emph{finite binary sequences} and define a refining
% sequence of partitions $\scrA=\{\al_m:m\geq0\}$ into
% intervals based on a strictly increasing positive sequence
% $0=a_0<a_1<a_2<\ldots$, $a_m\to\infty$, as follows:
% $A_\emptyset=\RR$, $A_0=(-\infty,-a_0)$,
% $A_1=[a_0,\infty)$, $A_{00}=(-\infty,-a_1)$, $A_{01}=[-a_1,a_0)$,
% $A_{10}=[a_0,a_1]$, $A_{11}=(a_1,\infty)$, and we continue
% splitting the outer-most intervals like this: for $m\geq2$,
% $\vep=0\ldots0\in\scrE_m$, $\iota_m=1\ldots1\in\scrE_m$ the
% elements $A_{\omicron_m0}, A_{\omicron_m1}, A_{\iota_m0},
% A_{\iota_m1}\in\al_m$ are given by $A_{\omicron_m0}=(-\infty,-a_m)$,
% $A_{\omicron_m1}=[-a_m,-a_{m-1})$, $A_{\iota_m0}=[a_{m-1},a_m]$,
% $A_{\iota_m1}=(a_m,\infty)$ (and, of course, suitable
% dyadic refinement into intervals of the intervening sets
% in the partitions $\al_m$). Such partitions are
% generated by a basis and resolve $\RR$.
% \end{example}

To arrive at random histogram distributions for the P\'olya-tree,
we define, for every $\vep\in\scrE$, a so-called \emph{splitting
variable} $V_{\vep0}$ (and $V_{\vep1}=1-V_{\vep0}$) taking values
in $[0,1]$ such that,
\begin{itemize}%[nosep]
\item[(i.)] for any $\vep,\vep'\in\scrE$ such that $\vep\neq\vep'$,
  $V_{\vep0}$ is independent of $V_{\vep'0}$;
\item[(ii.)] for every $\vep\in\scrE$, there exist
  $\be_{\vep0},\be_{\vep1}>0$ such that $V_{\vep0}$
  has a $\mathrm{Beta}(\be_{\vep0},\be_{\vep1})$
  distribution.
\end{itemize}
(In case $\scrX=[0,1]$ we assign a separately chosen, fixed
probability $0\leq p_0\leq1$ to $\{0\}$ with $\Pi_\al$-probability
one, for all $\al\in\scrA$. As a default, we choose $p_0=0$.)
\begin{remark}
Here and below, we extend the usual family of
Beta-distributions somewhat: we consider $\be_{\vep0}=\infty$
and $\be_{\vep1}=\infty$ and define
$\mathrm{Beta}(\infty,\be_{\vep1})=\delta_1$
for all $0<\be_{\vep1}<\infty$;
$\mathrm{Beta}(\be_{\vep0},\infty)=\delta_0$
for all $0<\be_{\vep0}<\infty$; and
$\mathrm{Beta}(\infty,\infty)=\delta_{1/2}$.
\end{remark}
The splitting variables $V_{\vep0}$ are interpreted as random
fractions that determine how much of the probability mass of
$A_\vep$ goes to $A_{\vep0}$ and how much remains for
$A_{\vep1}$, in accordance with (\ref{eq:splitA}):
\[
  V_{\vep0} = P( A_{\vep0}| A_\vep),
  \quad V_{\vep1} = P( A_{\vep1}| A_\vep).
\] 
Consequently, for every $m\geq1$, $\vep=e_1\ldots e_m\in\scrE_m$,
the random probability for $A_\vep$ can be written as a product
of independent fractions,
\[
  P(A_\vep) = V_{e_1}V_{e_1e_2}\ldots V_{e_1\ldots e_m}
    = \prod_{l=1}^m V_{e_1\ldots e_{l}},
\]
which fixes the histogram probability
measures $\Pi_{\al_m}$ on $\scrX_{\al_m}$ for all $m\geq1$,
\begin{equation}
  \label{eq:Polyainvsystem}
  \bigl( P(A_\vep):\vep\in\scrE_m \bigr) \sim \Pi_{\al_m}.
\end{equation}
By construction, the $\Pi_{\al_m}$ are such that
refinement and coarsening of partitions (corresponding
to relations of the type (\ref{eq:histogramadditive}))
are accommodated coherently.

For later reference, we note the first two moments of the
random variables $P(A_\vep)$: for every $m\geq1$ and every
$\vep\in\scrE_m$, the mean measure equals,
\begin{equation}
  \label{eq:polyameanmeasure}
  G(A_\vep):=\int_{M^1(\scrX_{\al_m})} P_{\al_m}
    (A_\vep)\dd\Pi_{\al_m}(P_{\al_m})
    =\int 
    \prod_{l=1}^m V_{e_1\ldots e_l}\dd\Pi_{\al_m}
    = \prod_{l=1}^m
    \frac{\be_{\vep_{l-1}e_{l}}}
    {\be_{\vep_{l-1}0} + \be_{\vep_{l-1}1}},
\end{equation}
by independence of the variables $V_{\vep0}$ and
expectations of the $\mathrm{Beta}$-distributions. Expressed
in terms of the parameters $\be$, the second moment of
$P(A_\vep)$ takes the form,
\[
  \begin{split}
  \int_{M^1(\scrX_{\al_m})}
    P_{\al_m}(A_{\vep})^2&\dd\Pi_{\al_m}(P_{\al_m})
  =\int
    \prod_{l=1}^m V_{e_1\ldots e_l}^2\dd\Pi_{\al_m}\\
  &=\prod_{l=1}^m
     \Biggl(\frac{\be_{\vep_{l-1}0}\be_{\vep_{l-1}1}}
     {(\be_{\vep_{l-1}0} + \be_{\vep_{l-1}1})^2
      (\be_{\vep_{l-1}0} + \be_{\vep_{l-1}1}+1)}
     +\frac{\be_{\vep_{l-1}e_l}^2}
     {(\be_{\vep_{l-1}0} + \be_{\vep_{l-1}1})^2}\Biggr),
  \end{split}
\]
based on independence of the $V_{\vep0}$, the variances of
the corresponding $\mathrm{Beta}$-distributions and
\cref{eq:polyameanmeasure}.

To have a sub-class of relatively simple examples,
we define so-called \emph{homogeneous P\'olya-tree systems}.
\begin{definition}
\label{def:homogeneouspolyatreesystem}
Let $\scrA$ denote the a \emph{dyadic tree} of partitions
of $\scrX=(0,1]$ (or $[0,1]$) as in example \ref{ex:dyadictree}.
A P\'olya-tree system is called \emph{homogeneous}, if we choose
$\be_m>0$ for all $m\geq1$, and set 
$\be_\vep=\be_m$, for all $\vep\in\scrE_m$.
\end{definition}
Accordingly, in a homogeneous P\'olya-tree system,
splitting variables are distributed symmetrically around
$v=\ft12$ and the mean measure $G$ for any homogeneous
P\'olya-tree system with a limit, is Lebesgue measure.

\subsection{Tight limits of P\'olya-tree histogram systems}
\label{sub:polyatightlimits}

First, the general case of the P\'olya-tree histogram system 
is analysed with corollary~\ref{cor:tightexistence}: here, the
particulars of the partition play a role in the formulation of
the condition, so we have to be specific regarding $\scrX$ and its
partitioning. In this subsection, we specify that $\scrX=(0,1]$
(or $\RR$), with a dyadic tree of partitions. We use the
following notation: for all $m\geq0$,
$\omicron_m=0\ldots0\in\scrE_m$ and $\iota_m=1\ldots 1\in\scrE_m$,
\begin{theorem}
\label{thm:existencepolya}
Let $\scrX=(0,1]$ and let $\scrA=\{\al_m:m\geq0\}$ be the
dyadic tree of example~\ref{ex:dyadictree}.
Let $(\Pi_{\al_m},\varphi_{*,\al_m\al_n})$ be a coherent inverse
system of P\'olya-tree measures (with parameter
$\be=\{\be_\vep:\vep\in\scrE\}$) on the inverse system
$(M^1(\scrX_{\al_m}), \varphi_{*\,\al_m\al_n})$.
Then there exists a unique probability measure on $M^1(\scrX)$
that is Radon with respect to the tight topology and projects
to the P\'olya-tree histograms parametrized by
$\{(\be_{\vep0},\be_{\vep1}):\vep\in\scrE\}$, if and only if,
\begin{equation}
  \label{eq:Polyaexistence}
  \prod_{m\geq0}
    \frac{\be_{\vep\omicron_m 0}}
      {\be_{\vep\omicron_m 0} + \be_{\vep\omicron_m 1}}=0,
\end{equation}
for every $\vep\in\scrE$, and the resulting random element
$P$ of $M^1(\scrX)$ is in the \emph{continuous-singular phase}.
\end{theorem}
\begin{proof}
Given $m\geq1$, the partition $\al_m$ consists of $2^m$ intervals
of the forms $(a,b]$ where $a=b-2^{-m}$ and $b=2^{-m}k$,
$k=1,2\ldots,2^m$, which is generated by a basis for the standard
topology on $(0,1]$. The well-ordered set of partitions
$\scrA=\{\al_m:m\geq1\}$ resolves $(0,1]$.
For given $\vep\in\scrE_m$ we consider $A_\vep\in\al_m$. Let also
$\delta,\eta>0$ be given. If $G_{\al_m}(A_\vep)=0$,
$P_{\al_m}(A_\vep)=0$, with $\Pi_{\al_m}$-probability one and
any compact $K\subset A_\vep$ satisfies
property~(\ref{property:Ptightprime}). Assuming that
$G_{\al_m}(A_\vep)>0$, we write $A_\vep=(a,b]$ for certain fixed
$a,b$ like above, and consider the
sequence of half-open intervals $(I_{\vep,l})_{l\geq m}$ in
$\scrX$, defined by,
\[
  I_{\vep,l}=A_\vep\setminus A_{\vep\omicron_{m-l}}=(a+2^{-l},b].
\]
% Note that $I_{\vep,l}$ increases to $A_\vep$ as $l\to\infty$. 
Assuming that (\ref{eq:Polyaexistence}) holds, choose $l\geq m$
large enough such that,
\[
  \prod_{k=0}^{l-1}
    \frac{\be_{\vep\omicron_k 0}}
      {\be_{\vep\omicron_k 0} + \be_{\vep\omicron_k 1}}<\frac{\delta\eta}{G_{\al_m}(A_\vep)}.
\]
Note that for all $m\leq k\leq l$,
$\varphi_{\al_k}(I_{\vep,l})=\varphi_{\al_k}(A_\vep)$, while
for any $l'\geq l$, $\varphi_{\al_l'}(I_{\vep,l})
=\varphi_{\al_l'}(A_\vep)\setminus\varphi_{\al_l'}(A_{\vep\omicron_{m-l}})$.
Defining $K$ to be the closure of $I_{\vep,l}$, by Markov's inequality,
we have,
\[
  \begin{split}
  \Pi_{\al_{l'}}\bigl( \{ & P_{\al_{l'}}\in M^1(\scrX_{\al_{l'}}):
    P_{\al_{l'}}(K)<P_{\al_{l'}}(\varphi_{\al_{l'}}(A_\vep))
    -\delta \}\bigr)\phantom{\prod_{k=0}^{l-1}}\\[.25em]
  &\leq \Pi_{\al_{l'}}\bigl( \{  P_{\al_{l'}}\in M^1(\scrX_{\al_{l'}}):
    P_{\al_{l'}}(I_{\vep,l})<P_{\al_{l'}}(\varphi_{\al_{l'}}(A_\vep))
    -\delta \}\bigr)\phantom{\prod_{k=0}}\\
  &= \Pi_{\al_l}\bigl( \{ P_{\al_l}\in M^1(\scrX_{\al_l}):
    P_{\al_l}(A_{\vep\omicron_{m-l}})>\delta \}\bigr)
  \leq \frac1\delta\int_{M^1(\scrX_{\al_l})}
    P_{\al_l}(A_{\vep\omicron_{m-l}})\dd\Pi_{\al_l}(P_{\al_l})\\
  &=\frac{G_{\al_m}(A_\vep)}\delta\prod_{k=0}^{l-1}
    \frac{\be_{\vep \omicron_{k}0}}
    {\be_{\vep \omicron_{k}0} + \be_{\vep \omicron_{k}1}}<\eta,
  \end{split}
\]
which shows that property~(\ref{property:Ptightprime}) holds.

Conversely, suppose that there exists a $\vep\in\scrE$, such that,
\[
  \prod_{m\geq0}
    \frac{\be_{\vep\omicron_m 0}}
      {\be_{\vep\omicron_m 0} + \be_{\vep\omicron_m 1}}>0,
\]
Then, 
\[
  \lim_{m\to\infty}\int P_{\al_m}(A_{\vep\omicron_m})
    \dd\Pi_{\al_m}(P_{\al_m})>0,
\]
while the sequence $(A_{\vep\omicron_m})_{m\geq0}$ decreases to
$\emptyset$. Hence, the mean measures $G_\al$ do not define a measure
(on the ring that is formed by the union of all $\sigma(\al)$,
$\al\in\scrA$), which precludes the existence of a Borel probability
measure $\Pi$ on $M^1(\scrX)$ with the tight topology (if $\Pi$
would exist, $B\mapsto\int P(B)\dd\Pi$ would define a Borel mean
measure).
\end{proof}
\begin{remark}
The above applies to examples with $\scrX=\RR$ as well, but in that case,
we have to require, in addition to (\ref{eq:Polyaexistence}), that,
\begin{equation}
  \label{eq:P1betas}
  \prod_{m\geq0}
    \frac{\be_{\iota_m 1}}
    {\be_{\iota_m 0} + \be_{\iota_m 1}}=0,
\end{equation}
because aside from the open, left-sided boundaries of half-open
intervals $A\in\al$, there are directions towards $\pm\infty$ where mass
can `leak away' in the limit.
\end{remark}
\begin{example}
\label{rem:polyadirichlet}
It is well known \citep{Lavine92}, that a P\'olya-tree histogram
system with defining parameters
$\{(\be_{\vep0},\be_{\vep1}):\vep\in\scrE\}$ that satisfy,
\begin{equation}
  \label{eq:polyadirichlet}
  \be_{\vep} = \be_{\vep0} + \be_{\vep1},
\end{equation}
for all $\vep\in\scrE$, coincides with a Dirichlet histogram
system (not on all of $\scrA_0$, but on a smaller set of dyadic
partitions $\scrA$ that resolves $\scrX$, generated by a basis).
Accordingly, such Dirichlet-P\'olya-tree histogram systems have limits
that are Radon probability measures on $M^1(\scrX)$ with the
tight topology, and the resulting random element
$P$ of $M^1(\scrX)$ is in the \emph{random-atomic phase}.
\end{example}

In the example below, we make a choice for the parameters
$\{(\be_{\vep0},\be_{\vep1}):\vep\in\scrE\}$ that gives
rise to a coherent histogram system without a tight limit.
This choice is not singular by construction, in the sense
that parameters either grow very large or vanish in the
limit: for all $\vep\in\scrE$, we have
$\be_{\vep0}^2+\be_{\vep1}^2=1$. To introduce the example,
we define the following function on $\scrE$.
\begin{definition}
\label{def:BaireCantor} % was {ex:BaireCantor}
In the standard construction of Cantor space as a subspace
$\scrC$ of $[0,1]$ by successive deletions of open
mid-sections of intervals, we define the %\inindex{Cantor mid-point function}
\emph{Cantor mid-point function} $x$ that parametrizes
the set of all mid-points of deleted intervals in terms
of finite binary sequences: 
$x:\scrE\to[0,1]$ maps $\vep\in\scrE_m$ to the midpoint
of the interval that is deleted in the $m$-th transition in the
construction of the set $\scrC$: for example, $x(\vep_\emptyset)=1/2$
in $\scrE_0$, $x(0)=1/6$, $x(1)=5/6$ in $\scrE_1$,
$x(00)=1/18$, $x(01)=5/18$, $x(10)=13/18$, $x(11)=17/18$ in
$\scrE_2$, \etc.
\end{definition}
\begin{example}
\label{ex:Polyatreethatdoesntexist}
Take $\scrX=\RR$ with a \emph{dyadic tree} of
partitions as defined in example~\ref{ex:coherentbutnoPi}, and,
for all $m\geq0$, $\vep\in\scrE_m$, 
\begin{equation}
  \label{eq:cossinpolya}
    \be_{\vep0}=\cos\bigl(\ft12\pi x(\vep)\bigr),\quad
    \be_{\vep1}=\sin\bigl(\ft12\pi x(\vep)\bigr).
\end{equation}
Note that,
\[
  \begin{split}
    \prod_{m\geq0}&\frac{\be_{\omicron_m 0}}
      {\be_{\omicron_m 0} + \be_{\omicron_m 1}}
    =\prod_{m\geq0}\frac{\cos\bigl(\ft12\pi x(\omicron_m)\bigr)}
      {\cos\bigl(\ft12\pi x(\omicron_m)\bigr)
        + \sin\bigl(\ft12\pi x(\omicron_m)\bigr)}\\
    &=\prod_{m\geq0}\Bigl(1+\tan\bigl(\ft12\pi x(\omicron_m)\bigr)\Bigr)^{-1}
    =\exp\biggl( -\sum_{m\geq0}
      \log\Bigl(1+\tan\bigl(\ft12\pi x(\omicron_m)\bigr)\Bigr)\biggr).
  \end{split}
\]
It is noted that $x(\omicron_m)=1/2(1/3)^m$ and,
\[
  \sum_{m\geq0}\log\Bigl(1+\tan\bigl(\ft12\pi x(\omicron_m)\bigr)\Bigr)
    \approx \sum_{m\geq0}\tan\bigl(\ft12\pi x(\omicron_m)\bigr)
    \approx \frac\pi{2} \sum_{m\geq0} x(\omicron_m)
    =\frac\pi{4}\sum_{m\geq0}\Bigl(\frac13\Bigr)^m=\frac{3\pi}{8}<\infty,
\]
Similarly,
\[
  \begin{split}
    \prod_{m\geq0}&\frac{\be_{\iota_m 1}}
      {\be_{\iota_m 0} + \be_{\iota_m 1}}
    =\prod_{m\geq0}\frac{\sin\bigl(\ft12\pi x(\iota_m)\bigr)}
      {\cos\bigl(\ft12\pi x(\iota_m)\bigr)
        + \sin\bigl(\ft12\pi x(\iota_m)\bigr)}\\
    &=\prod_{m\geq0}\Bigl(1+1/\tan\bigl(\ft12\pi x(\iota_m)\bigr)\Bigr)^{-1}
    =\exp\biggl( -\sum_{m\geq0}
      \log\Bigl(1+1/\tan\bigl(\ft12\pi x(\iota_m)\bigr)\Bigr)\biggr)
  \end{split}
\]
Since $x(\iota_m)=1-x(\omicron_m)$ for all $m\geq0$,
\[
  1/\tan\bigl(\ft12\pi x(\iota_m)\bigr)
    =1/\tan\bigl(\ft12\pi(1- x(\omicron_m))\bigr)
    =\tan\bigl(\ft12\pi x(\omicron_m)\bigr).
\]
Conclude that,
\[
  \prod_{m\geq0}
    \frac{\be_{\omicron_m 0}}
    {\be_{\omicron_m 0} + \be_{\omicron_m 1}}=
  \prod_{m\geq0}
    \frac{\be_{\iota_m 1}}
    {\be_{\iota_m 0} + \be_{\iota_m 1}}>0,
\]
which implies that the P\'olya-tree random histograms defined in 
(\ref{eq:cossinpolya}) form a coherent
system that \emph{does not} lead to a limiting probability measure
on $M^1(\RR)$ with the tight topology.
\end{example}

\subsection{Weak limits of P\'olya-tree histogram systems}
\label{sub:weakpolyalimits}

Second, we formulate a sufficient condition for the parameters
$\{(\be_{\vep0},\be_{\vep1}):\vep\in\scrE\}$ such that
the corresponding P\'olya-tree histogram system has a limit
$\Pi$ that is a Radon probability measure on $M^1(\scrX)$
with the weak topology. Based on this condition, it is
demonstrated that homogeneous P\'olya-tree systems with
$\be_m^{-1}=O(m^{-1})$ give rise to such weak histogram
limits. This rate of growth is lower than required in the
sufficient condition of (\cite{Kraft64}), which is elaborated
upon in \citep{Ferguson74,Mauldin92,Lavine92} and re-visited
in \citep{Ghosal17}.
\begin{theorem}
\label{thm:WeakExistencePolya}
Let $\scrX$ be a second countable metrizable space with countable
basis $\scrU$, and corresponding dyadic tree $\scrA$ of 
partitions $\al_m$, $m\geq1$, generated by the basis. Let
$(\Pi_{\al_m},\varphi_{*,\al_m\al_n})$ be a coherent inverse
system of P\'olya-tree measures (with parameter
$\be=\{\be_\vep:\vep\in\scrE\}$) on the inverse system
$(M^1(\scrX_{\al_m}), \varphi_{*\,\al_m\al_n})$. Assume
also that condition~(\ref{property:Ptight}) holds. Then there
exists a unique Radon probability measure $\Pi$ on
$M^1(\scrX)$ with the weak topology, projecting to $\Pi_{\al_m}$
for all $m\geq1$, if,
\begin{equation}
  \label{eq:Ppolya}
  \sup_{m\geq1} \sum_{\vep\in\scrE_m}
   \prod_{l=1}^m \frac1{\be_{\vep_{l-1}0} + \be_{\vep_{l-1}1}}
   \Biggl( \frac{\be_{\hat\vep}}
    {\be_{\vep_{l-1}0} + \be_{\vep_{l-1}1}+1}
   +\be_{\vep}\Biggr)
  <\infty.
\end{equation}
The resulting random element $P$ of $M^1(\scrX)$ is in the
\emph{absolutely-continuous phase}.
\end{theorem}
\begin{proof}
Condition~(\ref{property:Ptight}) implies the
existence of a tightly-Borel probability measure $\Pi'$ on
$M^1(\scrX)$ and a corresponding mean measure
$G\in M^1(\scrX)$, which serves as our choice of $Q$ in
the proof for property~(\ref{property:Pws}).
Let $\delta>0$ be given. For any $L>0$ and every $m\geq1$, Markov's inequality gives,
\[
  \begin{split}
  \Pi_{\al_m}\bigl(\{ &P_{\al_m}\in M^1(\scrX_{\al_m}):
    \|P_{\al_m}-P_{\al_m}\wedge LG_{\al_m}\|_{1,\scrX_{\al_m}}
    > \delta \}\bigr)\phantom{\biggm|}\\
  &\leq
  \Pi_{\al_m}\biggl(\Bigl\{ P_{\al_m}\in M^1(\scrX_{\al_m}):
    \frac1{L}\sum_{i\in I(\al_m)} \frac{P_{\al_m}(A_i)^2}{G_{\al_m}(A_i)}
    > \delta \Bigr\}\biggr)\\
  &\leq \frac1{L\delta} \sum_{i\in I(\al_m)} \frac1{G_{\al_m}(A_i)}
    \int_{M^1(\scrX_{\al_m})}P_{\al_m}(A_i)^2\dd\Pi_{\al_m}(P_{\al_m})\\
  &= \frac1{L\delta} \sum_{\vep\in\scrE_m}
     \prod_{l=1}^m\Biggl(\frac{\be_{\vep_{l-1}0}\be_{\vep_{l-1}1}}
       {\be_{\vep_{l-1}e_{l}}}
     \frac1{(\be_{\vep_{l-1}0} + \be_{\vep_{l-1}1})
      (\be_{\vep_{l-1}0} + \be_{\vep_{l-1}1}+1)}
     +\frac{\be_{\vep_{l-1}e_l}}
     {\be_{\vep_{l-1}0} + \be_{\vep_{l-1}1}}\Biggr)\\
  &< \frac{K}{L\delta},
  \end{split}
\]
for all $m\geq1$, where $K$ denotes the value of the supremum in
condition~(\ref{eq:Ppolya}). Consequently,
condition~(\ref{property:Pws}) is satisfied
and theorem~\ref{thm:WeakExistence} asserts that there exists
a unique weakly-Radon probability measure $\Pi$ on $M^1(\scrX)$
that projects to $\Pi_{\al_m}$ for all $m\geq1$.
\end{proof}
\begin{corollary}
\label{cor:weakpolya}
Assume the conditions of \cref{thm:WeakExistencePolya} and let
a sequence $\be_m>0$, ($m\geq1$) be given. If the $\be_m$ grow
like $m$ or faster, $\be_m^{-1}=O(m^{-1})$, there
exists a unique Radon probability measure $\Pi$ on $M^1(\scrX)$
with the weak topology, projecting to the associated homogeneous
P\'olya-tree histogram system.
\end{corollary}
\begin{proof}
Substituting $\be_\vep=\be_m$ in condition~(\ref{eq:Ppolya}),
we find, for every $m\geq1$,
\[
  \sum_{\vep\in\scrE_m}
   \prod_{l=1}^m \frac1{\be_{\vep_{l-1}0} + \be_{\vep_{l-1}1}}
   \Biggl( \frac{\be_{\hat\vep}}
    {\be_{\vep_{l-1}0} + \be_{\vep_{l-1}1}+1}
   +\be_{\vep}\Biggr)
  = \Biggl( \frac1{2\be_m+1}+1\Biggr)^m,
\]
which behaves like $\exp(m/(2\be_m+1))$ in the limit $m\to\infty$.
Since $m/\be_m=O(1)$ by assumption, the right-hand side stays
bounded and property~(\ref{eq:Ppolya}) is satisfied.
\end{proof}
Note that the sufficient condition of (\cite{Kraft64}; see also
\cite{Ghosal17}) suggests that absolute continuity of homogeneous
P\'olya-tree limits sets in when $\be_m$ grows as $O(m^2)$ or
faster; here it is shown that absolute continuity already
obtains with $\be_m$ that grow more slowly, as $O(m)$ or
faster.

% [BOOK ONLY?]
% Combining the conclusion of corollary~\ref{cor:weakpolya} with
% remark~\ref{rem:polyadirichlet}, we observe that homogeneous
% solutions of condition~(\ref{eq:polyadirichlet}) (\ie\
% a sequence $(\be_m)$ such that $\be_{m}=2\be_{m+1}$
% for all $m\geq0$) cannot grow with $m$, let alone
% grow as $O(m)$ or faster. So there are no
% homogeneous Dirichlet-P\'olya histogram systems that give rise to
% a random element $P$ of $M^1(\scrX)$ in the \emph{fixed-atomic
% phase}.

%%%%%%%%%%%%%%%%%%%%%%%%%%%%%%%%%%%%%%%%%%%%%%%%%%%%%%%%%%%%%%%%%%%%%%%%%%%%%%%
\section{Existence and phases of Gaussian histogram limits}
\label{sec:gauss}

Most known examples of random histogram systems with a limit,
are of the (normalized) completely random type \citep{daley07}.
The reason for the preference for
systems with independent components is \emph{coherence},
\cf\ (\ref{eq:histogramadditive}) or (\ref{eq:measuresyscoherence}),
which is analysed most conveniently with infinite divisibility,
requiring independence between summands.
The consequence is that most known histogram limits
are in one of the atomic phases of theorem~\ref{thm:phases}. In
this section we introduce the family of \emph{Gaussian random
measures}, random \emph{signed} measures on the space $\scrX$
with components that display dependence generically, manifesting
in one of the non-atomic phases of theorem~\ref{thm:phases}.

\subsection{Random histogram limits with signed measures}
\label{sub:signedlimits}

To arrive at a proof of existence for Gaussian histogram
limits, we have to generalize the approaches of
subsections~\ref{sub:WeakExistence}
and~\ref{sub:TightExistence}. Consider the case of a
\emph{locally compact} Polish space $\scrX$. The most natural
generalization of our random histogram question calls for
construction of Radon probability measures on $M(\scrX)$,
the space of all signed (and potentially unbounded) Radon
measures on $\scrX$, with the \emph{vague topology}
(see \IIthree{III}{1}{9}), rather than $M^1(\scrX)$
with the tight topology as in section~\ref{sec:prokhorov}. 
However to make histogram projections continuous,
transition to a zero-dimensional refinement $\hat{\scrX}$
(as in sub-subsection~\ref{subsub:zerodim}), is still a
necessary step, which does not combine well with the vague
topology (test functions $f:\scrX\to\RR$ for the vague
topology on $M(\scrX)$ remain continuous when viewed
as $f:\hat{\scrX}\to\RR$ but compactness of their
supports in $\hat{\scrX}$ is lost in general).
% This problem would not present itself with the
% tight topology, but the latter is defined only on the subspace
% of \emph{bounded} Radon measures $M_b(\scrX)\subset M(\scrX)$. 

However, since $M(\scrX)$ with the vague topology is the inverse
limit of the spaces $M(K)$ for compact $K\subset\scrX$, we may
also limit attention to \emph{compact} subsets $K\subset\scrX$
initially, and then use theorem~\ref{thm:prokhorov} (with
a directed set of compact $K\subset\scrX$ labelling a coherent
inverse system of histogram limits $\Pi_K$) to define a
limiting Radon probability measure $\Pi$ on $M(\scrX)$ with the
vague topology.

We defer proof of existence for such a `vague inverse limit of
histogram limits on compacta' to future work (\cite{Kleijn25b})
and focus here on the case where $\scrX$ itself is a
\emph{compact} Polish space. Then $M(\scrX)=M_b(\scrX)$ and
the vague and tight topologies coincide. Although $\hat{\scrX}$
is not compact in general, $M_b(\hat{\scrX})$ with the tight
topology still stands in continuous bijective
correspondence with $M(\scrX)$ (as in (the first part of)
proposition~\ref{prop:M1hatXtoM1X}), and the histogram
projections $\hat{\varphi}_{*\,\al}:M_b(\hat{\scrX})\to
M(\scrX_\al)$ of the form (\ref{eq:msrproj}),
are continuous. This enables the use of theorem~\ref{thm:prokhorov}
to prove existence of histogram limits $\Pi$ that are Radon
probability measures on $M(\scrX)$ with the vague/tight
topology.

% NOTE: Polish spaces are Lindel\"of spaces, and any locally
% compact Lindel\"of space is `'exhausted by compacts', which means
% there exists a sequence of compacta $K_n\subset\scrX$, $(n\geq1)$, such
% that their (countable) union covers $\scrX$ and, for every
% $n\geq1$, $K_{n}$ is a subset of the \emph{interior} of
% $K_{n+1}$ (which makes the property slightly stronger than
% $\sigma$-compactness).

\begin{theorem}
\label{thm:vagueexistence} % was {thm:vagueGaussian}, {thm:tightGaussian}
Let $\scrX$ be a compact Polish space and let $\scrA$ be a
directed set of partitions that resolves
$\scrX$, generated by a basis hat gives rise to a zero-dimensional
$\hat{\scrX}$. Consider $M(\scrX)$ with the tight topology and a
coherent random histogram system $(\Pi_{\al},\varphi_{*\,\al\be})$. 
If,
\begin{itemize}
\item[(i.)] for every $\ep>0$ there is a constant $M>0$ such
  that for all $\al\in\scrA$,
  \[
    \tag{P1-signed}\label{property:P1signed} % {property:P1signed}
    \Pi_{\al}\bigl(\bigr\{ \Phi_\al\in M(\scrX_\al)\,:\,
      \|\Phi_\al\|_{1,\al}>M \bigr\}\bigr) < \ep;
  \]
\item[(ii.)] and, for every $\ep,\delta>0$ there is a compact
  $\hat{K}\subset\hat{\scrX}$ such that for all $\al\in\scrA$,
  \[
    \tag{P2-signed}\label{property:P2signed}
    \Pi_{\al}\bigl(\bigr\{ \Phi_\al\in M(\scrX_\al)\,:\,
      |\Phi_\al|(\varphi_\al(\hat{K}))<|\Phi_\al|(\scrX_\al)
        -\delta \bigr\}\bigr) < \ep,
  \]
\end{itemize}
then there exists a unique Radon probability distribution
$\Pi$ on $M(\scrX)$ projecting to $\Pi_{\al}$ for all $\al\in\scrA$.
\end{theorem}
\begin{proof}
According to proposition~\ref{prop:msrstaysmsr}, the Borel sets
on $\scrX$ and $\hat{\scrX}$, as well as (bounded) signed
Borel set functions and measures are the same. The
first part of proposition~\ref{prop:M1hatXtoM1X} remains
true, (the set-theoretic identity mapping
$i_*:M_b(\hat{\scrX})\to M(\scrX)$ is a tight-to-tight
continuous bijection), but the second part fails because
$M(\scrX)$ and $M_b(\hat{\scrX})$ are not necessarily Polish
spaces (see remark~\ref{rem:almostnecessary}). Like before,
the mappings $(\hat\varphi_{*\,\al},
\varphi_{*\,\al\be})$ form a coherent and separating family of
continuous mappings on $M_b(\hat{\scrX})$. Trivially, the linear
spaces $M(\scrX_\al)$ are finite-dimensional normed
spaces (the vague, tight and total-variational topology are all
equivalent) and the mappings $\varphi_{*\,\al\be}$ are
surjective and continuous, like in proposition~\ref{prop:invsysmsr}.
Accordingly, $(M(\scrX_\al),\varphi_{*\,\al\be})$ forms
an inverse system.

To show that condition~(\ref{property:P}) holds, let
$\ep>0$ be given and let $M>0$ be a constant such that,
property~(\ref{property:P1signed}) is satisfied for every
$\al\in\scrA$. Define,
\[
  H_1 = \bigcap\bigl\{ \Phi\in M(\hat{\scrX})\,:\,
    \bigl\|\Phi_\al\bigr\|_{1,\al}\leq M,\, \al\in\scrA \bigr\}
    = \bigl\{ \Phi\in M(\hat{\scrX})\,:\,
    \bigl\|\Phi\bigr\|_{1,\scrX}\leq M\},
\]
since $\|\Phi\|_{1,\scrX}= \sup_{\al\in\scrA} \|\Phi_\al\|_{1,\al}$,
\cf\ proposition~\ref{prop:alphatotalvariation}.

Let $\hat{K}$ be a compact subset of $\hat{\scrX}$. For any
$\Phi\in H_1$,
\[
  |\Phi|(\hat{K})\leq |\Phi|(\hat{\scrX}) =|\Phi|(\scrX)
    = \|\Phi\|_{1,\scrX}
    = \sup_{\al\in\scrA} \|\Phi_\al\|_{1,\al}\leq M,
\]
With $(\ep_n), (\delta_n)$ as in the proof of
theorem~\ref{thm:tightexistence}, and the corresponding
compacta $(\hat{K}_n)$ in $\hat{\scrX}$ \cf\
property~(\ref{property:P2signed}), we also
define,
\[
  H_2 = \bigcap\bigl\{ \Phi\in M(\hat{\scrX})\,:\,
    |\Phi_\al|(\varphi_\al(\hat{K}_n))\geq |\Phi_\al|(\scrX_\al)-\delta_n,\, 
    n\geq1,\, \al\in\scrA \bigr\}.
\]
Following steps analogous to those in the proof of
theorem~\ref{thm:tightexistence}, one then finds that
$H=H_1\cap H_2$ satisfies Prokhorov's conditions
(see~(\ref{eq:prokhorovscondition})), so that the
closure $\overline{H}$ forms a compact subset of
$M(\hat{\scrX})$ and for any $\al$, we have (by
monotony of $\al\mapsto|\Phi_\al|(B)$ for any
$B\in\sigma(\al)$, like in the proofs of
theorems~\ref{thm:WeakExistence}
and~\ref{thm:tightexistence}),
\[
  \begin{split}
  \Pi_\al\bigl(
    M(\scrX_\al&)\setminus\hat\varphi_{*\,\al}(\overline{H}) \bigr)
  \leq \Pi_\al\bigl( \bigl\{ \Phi\in M(\hat{\scrX})\,:\,
    \bigl\|\Phi_\al\bigr\|_{1,\al}> M \bigr\} \bigr)\\[3mm]
    &+ \sum_{n\geq1}\Pi_\al\bigl( \bigl\{ \Phi_\al\in M(\scrX_\al)\,:\,
        |\Phi_\al|(\varphi_\al(\hat{K}_n))<|\Phi_\al|(\scrX_\al)
        -\delta_n \bigr\}\bigr)
    < 2\ep,
  \end{split}
\]
which shows that condition~(\ref{property:P}) of
theorem~\ref{thm:prokhorov} is satisfied. Conclude that there exists a
unique Radon probability measure
$\hat\Pi$ on $M(\hat{\scrX})$ that projects to $\Pi_\al$ for all
$\al\in\scrA$. The continuous mapping $i_*:M(\hat{\scrX})\to M(\scrX)$
serves to define $\Pi=\hat\Pi\circ i_*^{-1}$, a Radon probability measure
on $M(\scrX)$, and $\Pi$ still projects to $\Pi_\al$ for all $\al\in\scrA$.
\end{proof}
As it stands, property~(\ref{property:P2signed}) is somewhat
unwieldy due to occurrence of compacta in $\hat{\scrX}$. Analogous
to corollary~\ref{cor:tightexistence}, we also provide a version
that refers only to compacta in $\scrX$. For brevity's sake, we
omit the proof (which follows the precise same steps of the proof of
corollary~\ref{cor:tightexistence}): if,
\begin{itemize}
\item[] for all $\al\in\scrA$, all $A\in\al$ and all $\ep,\delta>0$,
  there is a $K\subset A$, compact in $\scrX$, such that,
  \[
    \tag{P2-signed'}\label{property:P2signedprime}
    \Pi_\be\bigl(\bigr\{ \Phi_\be\in M(\scrX_\be)\,:\,
      |\Phi_\be|(\varphi_\be(K))< |\Phi_\be|(\varphi_\be(A))-\delta
        \bigr\}\bigr) < \ep,
  \]
  for all $\be\in\scrA$ such that $\al\leq\be$,
\end{itemize}
then property~(\ref{property:P2signed}) is satisfied.

\begin{remark}
\label{rem:almostnecessary}
Regarding the pair of properties~(\ref{property:P1signed})
and~(\ref{property:P2signed}), we remark that,
unlike earlier applications of theorem~\ref{thm:prokhorov}, our
conditions are sufficient but (perhaps) not necessary for the existence
of a histogram limit: note that $M(\scrX)$ is not necessarily complete
and not a Polish space generically (see \II{III}{1}{9}{Proposition~14}),
so that Borel measurability of the inverse mapping $i^{-1}_*$
can no longer be guaranteed. Accordingly, not every Radon probability
measure on $M(\scrX)$ can be extended to a Borel probability measure
on $M(\hat{\scrX})$ canonically, and there may exist coherent
histogram systems with an tight inverse limit $\Pi$ on $M(\scrX)$,
for which stated conditions do not hold.
\end{remark}

\subsection{Existence of tight Gaussian histogram limits}
\label{sub:gaussianlimits}

For the definition of Gaussian histogram systems, location
and covariance parameters are defined in a way comparable to
that of the base measure of the Dirichlet family.
\begin{definition}
\label{def:Gaussparameters}
Let $\scrX$ be a compact Polish space, let $\lambda$ be a
signed Radon measure on $\scrX$. For any Borel measurable
partition $\al$, let $\lambda_{\al}$ denote the
$\al$-histogram projection of $\lambda$,
$\lambda_{\al} = (\lambda(A_1),\ldots,\lambda(A_{{|\al|}}))$
in $M(\scrX_\al)$. Let $\Sigma$ be a signed,
symmetric Radon measure on $\scrX\times\scrX$ (symmetric meaning that
$\Sigma(A\times B) = \Sigma(B\times A)$, for all $A,B\in\scrB$).
Assume that for every $\al\in\scrA$, the
${|\al|}\times {|\al|}$-matrix $\Sigma_{\al}$, with entries,
\[
  \Sigma_{\al,ij} = \Sigma(A_i\times A_j),
\]
($1\leq i,j\leq {|\al|}$) is \emph{semi-positive definite}. We
refer to $\lambda$ as the \emph{centre measure}, to $\Sigma$ as
the \emph{covariance measure}, and to $(\lambda, \Sigma)$ as
\emph{Gaussian parameters}.
% For every $\al\in\scrA$, the 
% eigenvalues of $\Sigma_\al$ are denoted
% $\{\tau_{\phantom{|}\Sigma,\al,i}:1\leq i\leq {|\al|}\}$,
% and the largest eigenvalue  is denoted $\tau_{\phantom{|}\Sigma,\al}$. 
\end{definition}
The measure $\Sigma$ may be viewed equivalently as a linear
mapping that takes continuous functions on $\scrX$ into
signed Radon measures on $\scrX$,
% : for $f\in C(\scrX)$,
% and $B\in\scrB$,
% \[
%   \Sigma(f)(B)=\int_{B\times\scrX} f(y)\dd\Sigma(x,y),
% \]
or as a symmetric bi-linear form.
% : for $f,g\in C(\scrX)$,
% \[
%   \Sigma(f,g)=\int_{\scrX\times\scrX} f(x)g(y)\dd\Sigma(x,y)
% \]
To give examples, we may turn to the theory of
\emph{reproducing kernel Hilbert spaces}.
\begin{example}
\label{ex:kernelfunctions}
For $d\geq1$, let $\scrX$ be a compact subset of $\RR^d$, and
consider a so-called \emph{positive-definite symmetric
kernel function} $k:\scrX\times\scrX\to\RR$;
\cf\ the \emph{Moore–Aronszajn theorem}, every such
kernel function is the reproducing kernel for a unique Hilbert
space of functions on $\scrX$. We use $k$ to define,
\begin{equation}
  \label{eq:RKHSkernel}
  \Sigma_k(A\times B) = \int_{A\times B} k(x,y)\dd x\dd y,
\end{equation}
and note that such a $\Sigma_k$ is a covariance measure in the
sense of definition~\ref{def:Gaussparameters}. \emph{Mercer's
theorem} formulates the associated spectral theory, with
$\Sigma_k$ viewed as a (compact, self-adjoint, positive)
integral operator $L^2(\scrX)\to L^2(\scrX)$. Indeed, we
may define a kernel by choice of a countable orthonormal
subset of continuous functions $\{\phi_i:i\in I\}$, and
non-negative $\{\lambda_i:i\in I\}$, to define $k(x,y)=
\sum_{i\in I} \lambda_i\phi_i(x)\phi_i(y)$.
\end{example}
To extend the previous example for general covariance measures
$\Sigma$, note that if $\scrA$ consists of partitions generated
by a basis for $\scrX$, then for any continuous $f:\scrX\to\RR$,
\[
  \langle f,f\rangle_\Sigma:=
    \int_{\scrX\times\scrX}f(x)f(y)\dd\Sigma(x,y)\geq0.
\]
Polarization then defines a positive semi-definite bilinear
form on $C(\scrX)\times C(\scrX)$. Within $C(\scrX)$,
there is a linear space $Q_\Sigma$ of functions $f$ that are
$\Sigma$-almost-surely equal to zero ($f\sim 0$ if
$\langle f,f\rangle_\Sigma=0$; $f\sim g$ if $f-g\sim0$.),
and the quotient space $C(\scrX)/Q_\Sigma$ is a real pre-Hilbert
space (see \cite{Treves67}, def.~7.5) with Hilbert space
completion denoted $L^2(\Sigma)$, which generalizes
the Moore-Aronszajn Hilbert spaces associated to reproducing
kernel functions.

\begin{definition}
Given Gaussian parameters $(\lambda,\Sigma)$, we define a
Gaussian histogram system as follows: for all $\al\in\scrA$,
we choose normal probability distributions
$\Pi_{\lambda,\Sigma,\al}$ for random signed histograms
$\Phi_\al\in M(\scrX_\al)$, as follows:
\[
  \bigl(\Phi_\al(A_1),\ldots,\Phi_\al(A_{{|\al|}})\bigr)
    \sim \Pi_{\lambda,\Sigma,\al}
    = N\bigl( \lambda_\al, \Sigma_\al \bigr),
\]
where $N(\lambda_\al,\Sigma_\al)$ denotes the multivariate normal
distribution on $\RR^{{|\al|}}$ with expectation $\lambda_\al$
and covariance matrix $\Sigma_\al$. When $\lambda=0$, we speak of
a \emph{centred} Gaussian histogram distribution, denoted
$\Pi_{\Sigma,\al}=N(0,\Sigma_\al)$.
\end{definition}
For partitions $\al,\be\in\scrA$, where $\be$ refines $\al$,
let $\varphi_{*\,\al\be}$ be as in (\ref{eq:transitionmapping}),
the mapping that expresses finite additivity. Below, we show that
for any $\scrA$ and any Gaussian parameters $(\lambda,\Sigma)$,
the above histogram distributions define a coherent system,
referred to as the \emph{Gaussian histogram system}
$(\Pi_{\lambda,\Sigma,\al},\varphi_{*\,\al\be})$
associated with the parameters $(\lambda,\Sigma)$.

The inclusion of a centre measure $\lambda$ is not of influence
for the existence of Gaussian histogram limits: for all
$\al\in\scrA$ and all Borel sets $B$ in $M(\scrX_\al)$,
\[
  \Pi_{\lambda,\Sigma,\al}\Bigl( \bigl\{
    \Phi_\al\in M(\scrX_\al): \Phi_\al \in B\bigr\}\Bigr)
  = \Pi_{\Sigma,\al}\Bigl( \bigl\{
    \Phi_\al\in M(\scrX_\al): \Phi_\al \in B-\lambda_\al \bigr\}\Bigr)
\]
and hence, $\Pi_{\lambda,\Sigma}$ exists,
if and only if, $\Pi_{\Sigma}$ exists. Existence of
histogram limits therefore only concerns the $\Sigma$-parameter.
The existence conditions of theorem~\ref{thm:vagueexistence}
can be dominated by uniform bounds on
mixed second moments of the absolute histogram components
$|\Phi_\al(A)|$, which we denote by,
\begin{equation}
  \label{eq:TAB}
  \begin{split}
  T(A,B) &= \int_{M(\scrX_\al)}
      \bigl| \Phi_\al(A)\Phi_\al(B) \bigr|\dd\Pi_{\Sigma,\al}(\Phi_\al)\\
    &= C \, \Sigma(A\times B)
    + \frac{2}{\pi}
    \sqrt{\Sigma(A\times A)\Sigma(B\times B)-\Sigma(A\times B)^2},
  \end{split}
\end{equation}
where $0\leq C\leq 1$ is a constant
that depends only on the correlation coefficient between
$\Phi_\al(A)$ and $\Phi_\al(B)$ (see \cite{Kan17}, p.~933).
\begin{corollary}
\label{cor:ExistenceGaussTight}
Let $\scrX$ be a compact Polish space and let $\scrA$ be a
directed set of partitions generated by a basis as in
example~\ref{ex:partitionscountablebasis}. Let $\Sigma$
be a covariance measure on $\scrX\times\scrX$. Consider
$M(\scrX)$ with the tight topology and the centred Gaussian
histogram system $(\Pi_{\Sigma,\al},\varphi_{*\,\al\be})$.
If the covariance measure $\Sigma$ is such that,
\begin{itemize}
\item[(i.)] 
  \begin{equation}
    \label{eq:supEabssqrd} % was {eq:maxeigenvaluesigma}
    \sup_{\al\in\scrA} \,\sum_{A,B\in\al} T(A,B) < \infty;
  \end{equation}
\item[(ii.)] and, for any open $U_n\in\cup\{\sigma(\al):\al\in\scrA\}$
  that decrease to $\emptyset$,  
  \begin{equation}
    \label{eq:localsupEabssqrd}
    \lim_{n\to\infty} \, \sup_{\al\in\scrA}\,
      \sum\bigl\{T(A,B):
        \,A,B\in\al,\,A,B\subset U_n\bigl\} = 0,
  \end{equation}
\end{itemize}
then there exists a unique Radon probability distribution
$\Pi_{\Sigma}$ on $M(\scrX)$ projecting to $\Pi_{\Sigma,\al}$ for all
$\al\in\scrA$.
\end{corollary}
\begin{proof}
%Consider first the simplest of refinements of the type
% (\ref{eq:intersect}): denote $\al=(A_1,A_2)$ and
% $\be=(B_1,B_2,B_3)$, where $A_1=B_1\cup B_2$ and $A_2=B_3$.
% Because $\lambda$ is finitely additive,
% $\lambda(A_1)=\lambda(B_1)+\lambda(B_2)$ and
% $\lambda(A_2)=\lambda(B_3)$. Similarly,
% \[
%   \begin{split}
%   \bigl(\Sigma_\al\bigr)_{11} &= \bigl(\Sigma_\be\bigr)_{11}
%     + 2\bigl(\Sigma_\be\bigr)_{12} + \bigl(\Sigma_\be\bigr)_{22},\\
%   \bigl(\Sigma_\al\bigr)_{12} &= \bigl(\Sigma_\be\bigr)_{13}
%     + \bigl(\Sigma_\be\bigr)_{23},\\
%   \bigl(\Sigma_\al\bigr)_{22} &= \bigl(\Sigma_\be\bigr)_{33},
%   \end{split}
% \]
% This verifies that,
% \[
%   \lambda_\al = P \lambda_\be,\quad \Sigma_\al = P\,\Sigma_\be\,P^T,
% \]
% where $P:\RR^3\to\RR^2$ is the linear mapping,
% \[
%   P = \threebytwo%
%     {1}{1}{0}%
%     {0}{0}{1}.
% \]
% More generally, if
To use theorem~\ref{thm:vagueexistence}, we first verify coherence
of Gaussian histogram systems. (We do this first step of the proof
generically, that is, with $\lambda\neq0$.) If
$\al,\be\in\scrA$, $\al\leq\be$, and we write
$A_i=\cup_{k\in J_{\al\be}(i)}B_k$, then,
\[
  \bigl(\lambda_\al\bigr)_i = \sum_{k\in J_{\al\be}(i)}
    \bigl(\lambda_\be\bigr)_k,\quad
  \bigl(\Sigma_\al\bigr)_{ij} = \sum_{k\in J_{\al\be}(i)}
    \sum_{l\in J_{\al\be}(j)}\bigl(\Sigma_\be\bigr)_{kl},
\]
for $1\leq i,j\leq {|\al|}$. This can %again
be expressed in terms of a linear mapping
$P_{\al\be}:\RR^{{|\be|}}\to\RR^{{|\al|}}$ such that,
\[
  \lambda_\al = P_{\al\be} \lambda_\be,\,\,\,
  \Sigma_\al = P_{\al\be}
    \,\Sigma_\be\,(P_{\al\be})^t,
\]
(where $(P_{\al\be})^t$ denotes the matrix transpose of
$P_{\al\be}$). Recall that for any finite $d,d'\geq1$, any linear
$P:\RR^d\to\RR^{d'}$ and a random variable $\Phi$ distributed
$N(\lambda,\Sigma)$, the random variable $P\Phi$ is distributed
$N(P\lambda,P\Sigma P^T)$. So
$\Phi_\al\sim\Pi_{\lambda,k,\al}$ has the same distribution as
$\varphi_{*\,\al\be}(\Phi_\be)$ for $\Phi_\be\sim\Pi_{\lambda,k,\be}$,
\[
  \varphi_{*\,\al\be}(\Phi_\be)_i =\sum_{j\in J_{\al\be}(i)}(\Phi_{\be})_j
  = \bigl( P\Phi_\be \bigr)_i,
\]
for all $1\leq i\leq {|\al|}$. This verifies coherence of the histogram
system $(\Pi_{\lambda,\Sigma,\al},\varphi_{*\,\al\be})$.

In the rest of the proof, we assume that the histogram system is
centred: $\lambda=0$. To show that property~(\ref{property:P1signed})
holds, we use \emph{Chebyshev's inequality} to upper-bound its
left-hand side, for every $\al\in\scrA$ and $M>0$:
\[
  \begin{split}
  \Pi_{\al}\bigl(\bigr\{ \Phi_\al&\in M(\scrX_\al)\,:\,
      \|\Phi_\al\|_{1,\al}>M \bigr\}\bigr)
    = \Pi_{\Sigma,\al}\biggl(\Bigr\{ \Phi_\al\in M(\scrX_\al)\,:\,
      \sum_{A\in\al}|\Phi_\al(A)|>M \Bigr\}\biggr)\\
    &= \Pi_{\Sigma,\al}\biggl(\Bigr\{ \Phi_\al\in M(\scrX_\al)\,:\,
      \sum_{A,B\in\al}|\Phi_\al(A)\Phi_\al(B)|>M^2 \Bigr\}\biggr)\\
    &\leq \frac1{M^2}\int_{M(\scrX_\al)} \sum_{A,B\in\al}
     \bigl| \Phi_\al(A)\Phi_\al(B) \bigr|\dd\Pi_{\Sigma,\al}(\Phi_\al).
  \end{split}
\]
Assuming condition~(\ref{eq:supEabssqrd}) and choosing $M$ large
enough, we see that property~(\ref{property:P1signed}) is satisfied
for every $\ep>0$ and all $\al\in\scrA$.

Let $\al\in\scrA$ and $A\in\al$ be given. Because $\al$ is
generated by a basis, $A$ is the intersection $U\cap C$ of an (open)
finite intersection $U$ of basis sets and a (closed) finite
intersection $C$ of complements of basis sets. Because $\scrX$ is a
Polish space, $U$ is $F_\sigma$, \ie\ $U$ is equal to a
countable union of closed sets $U=\cup_{m\geq1}C_m$. Then
the closed sets $C'_n=C\cup(\cup_{1\leq m\leq n}C_m)$ increase to $A$
as $n\to\infty$. In the open complements of $C_n'$ in $A$, there
exists a decreasing sequence of basis elements $U_n$, and we may
define closed sets $K_n=A\setminus U_n$. The open sets
$U_n$ decrease to $\emptyset$ and the sets $K_n$ are closed
as subsets of $\scrX$ and therefore compact. By assumption (see
example~\ref{ex:partitionscountablebasis}), $\scrA$
is such that all sets in the basis occur as elements of
$\sigma(\al)$ for some $\al\in\scrA$. So for every $n\geq1$ and
every $\al\in\scrA$, there exists a $\beta\geq\al$ such that the
decomposition $A=U_n\cup K_n$ is such that $U_n,K_n\in\sigma(\be)$.

Let $\ep,\delta>0$ be given. Based on assumption (\ref{eq:localsupEabssqrd}),
choose $n\geq1$ large enough such that,
\[
  \sup_{\gamma\geq\al}
    \sum\bigl\{T(B,C):\,B,C\in\gamma,\,B,C\subset U_n\bigl\}
    < \delta^2\ep.
\]
Let $\gamma\geq\al$ be given. Since
$\varphi_\gamma(A)\setminus\varphi_\gamma(K_n)
\subseteq\varphi_\gamma(U_n)$, 
\[
  \begin{split}
  \Pi_{\Sigma,\gamma}\bigl(\bigr\{ \Phi_\gamma&\in M(\scrX_\gamma)\,:\,
      |\Phi_\gamma|(\varphi_\gamma(K_n))< |\Phi_\gamma|(\varphi_\gamma(A))-\delta
        \bigr\}\bigr)\phantom{\biggl(}\\
    &\leq \Pi_{\Sigma,\gamma}\bigl(\bigr\{ \Phi_\gamma\in M(\scrX_\gamma)\,:\,
      |\Phi_\gamma|(\varphi_\gamma(U_n)) > \delta\bigr\}\bigr)\phantom{\biggl(}\\
    &= \Pi_{\Sigma,\gamma}\biggl(\biggr\{ \Phi_\gamma\in M(\scrX_\gamma)\,:\,
      \sum\Bigl\{|\Phi_\gamma(B)\Phi_\gamma(C)|\,:
        \,B,C\in\gamma,\,B,C\subset U_n\Bigl\} >\delta^2 \biggr\}\biggr)\\
    &\leq \frac1{\delta^2}\int_{M(\scrX_\gamma)} \sum_{B,C\subset U_n}
     \bigl| \Phi_\gamma(B)\Phi_\gamma(C) \bigr|\dd\Pi_{\Sigma,\gamma}(\Phi_\gamma)
    = \frac1{\delta^2}\sum_{B,C\subset U_n} T(B,C) < \ep,\phantom{\biggl(}
  \end{split}
\]
showing that property~(\ref{property:P2signedprime}) holds.
\end{proof}
In the above proof, we satisfy property~(\ref{property:P2signedprime})
by roughly following the proof of theorem~\ref{thm:existencepolya}, but
with a different construction of the compacta $K_n$, which is more
generic and based on example~\ref{ex:partitionscountablebasis}. In
applications, control over the choice of $\al$ allows for convenient
constructions. Where the proof of theorem~\ref{thm:existencepolya}
depends on the Carath\'eodory-like condition that
$G_\al(\varphi_\al(U_n))\to0$ for (specific) sequences $(U_n)$ in
the generating ring that decrease to $\emptyset$, here the
second-absolute-moment set-functions of
condition~(\ref{eq:localsupEabssqrd})
are required to go to zero.

Based on condition~(\ref{eq:supEabssqrd}), we briefly come back to
the space $L^2(\Sigma)$ and indicate how it is related to the
covariance structure of a centred Gaussian histogram limit
$\Pi_\Sigma$. To that end, we define the real-valued stochastic
integrals $\Phi_f:=\int_{\scrX}f\dd\Phi$, for $f\in C(\scrX)$,
and consider the linear space of real-valued random variables
$L=\{\Phi_f:f\in C(\scrX)\}$ that they span. Assuming integrability,
on $L$ we define the bilinear form,
\[
  \langle \Phi_f,\Phi_g \rangle_L :=
    \int_{M(\scrX)}\Phi_f\,\Phi_g\,\dd\Pi_\Sigma(\Phi).
\]
With $Q_L=\{\Phi_f\in L:\langle \Phi_f,\Phi_f \rangle_L=0\}$,
the quotient space $L/Q_L$ with $\langle \cdot,\cdot \rangle_L$ as
an inner product is a real pre-Hilbert space, with Hilbert space
completion we denote by $\hat{L}$. The following proposition
involves histogram approximations of continuous functions: for
any $f\in C(\scrX)$ and any partition $\al$ generated by a basis,
let $f_\al(A)$, $(A\in\al)$ be real numbers such that
$\inf\{f(x):x\in A\}\leq f_\al(A)\leq \sup\{f(x):x\in A\}$, and
let $f_\al(x)=\sum_{A\in\al} f_\al(A)1_A(x)$, noting that, for
all $x\in\scrX$, $f_\al(x)\to f(x)$ as $\al$ refines within an
$\scrA$ that resolves $\scrX$.
\begin{proposition}
Let $\scrX$ be a compact Polish space and let $\scrA$ be a
directed set of partitions generated by a basis as in
example~\ref{ex:partitionscountablebasis}. Let $\Sigma$ be a
covariance measure for which conditions~(\ref{eq:supEabssqrd})
and~(\ref{eq:localsupEabssqrd}) hold. Then, for all $f,g\in C(\scrX)$,
\[
  \langle f,g\rangle_\Sigma
  = \lim_\al\sum_{A,B\in\al}f_\al(A)g_\al(B)\,\Sigma(A \times B)
  = \langle \Phi_f,\Phi_g \rangle_L,\phantom{\int}
\]
and the mapping $C(\scrX)\to L:f\mapsto\Phi_f$ extends
to an isometric isomorphism of Hilbert spaces $L^2(\Sigma)\to\hat{L}$.
\end{proposition}
\begin{proof}
We first show that the (semi-definite) inner-product spaces
$C(\scrX)$ and $L$ are isometrically isomorphic.
Let $f,g:\scrX\to\RR$ be continuous and denote their supremum norms
by $\|f\|_{\infty,\scrX},\|g\|_{\infty,\scrX}<\infty$. By
continuity of $\RR\to\RR:(x,y)\mapsto xy$,
\[
  \begin{split}
  \langle f,g\rangle_\Sigma
    &=\int_{\scrX^2}f(x)g(y)\dd\Sigma(x,y)
    =\int_{\scrX^2} \Bigl(\lim_\al\sum_{A\in\al}f_\al(A)1_A(x)\Bigr)\,
      \Bigl(\lim_\be\sum_{B\in\be}g_\be(B)1_B(y)\Bigr)\dd\Sigma(x,y)\\
    &=\int_{\scrX^2}\lim_\gamma\sum_{A,B\in\gamma}f_\gamma(A)g_\gamma(B)
      1_A(x)1_B(y)\dd\Sigma(x,y)
    =\lim_\gamma\sum_{A,B\in\gamma}f_\gamma(A)g_\gamma(B)\,\Sigma(A \times B),
  \end{split}
\]
and the last
step holds by Lebesgue's dominated convergence, based on the facts
that $\|f\|_{\infty,\scrX},\|g\|_{\infty,\scrX}<\infty$ and
$\Sigma(\scrX\times\scrX)<\infty$. Using the definition
of $\Phi_\al$, we may then write,
\[
  \begin{split}
  \langle f,g\rangle_\Sigma
    &= \lim_\al\sum_{A,B\in\al}\int_{M(\scrX_\al)}
      f_\al(A)\Phi_\al(A)\,g_\al(B)\Phi_\al(B)\dd\Pi_{\Sigma,\al}(\Phi_\al)\\
    &= \lim_\al\int_{M(\scrX)}\sum_{A,B\in\al}
      f_\al(A)\Phi(A)\,g_\al(B)\Phi(B)\dd\Pi_{\Sigma}(\Phi)
  \end{split}
\]
To complete the argument, note that
\[
  \biggl|\sum_{A,B\in\al}
    f_\al(A)\Phi(A)\,g_\al(B)\Phi(B)\biggr|
    \leq \|f\|_{\infty,\scrX}\,\|g\|_{\infty,\scrX}
      \sum_{A,B\in\al}\bigl|\Phi(A)\Phi(B)\bigr|
\]
and the right-hand side is monotone increasing in $\al$.
By monotone convergence,
\[
  \lim_{\al\in\scrA}\sum_{A,B\in\al}\int_{M(\scrX)}\bigl|\Phi(A)\Phi(B)\bigr|\dd\Pi_\Sigma(\Phi)
    = \int_{M(\scrX)}\lim_{\al\in\scrA}\sum_{A,B\in\al}\bigl|\Phi(A)\Phi(B)\bigr|\dd\Pi_\Sigma(\Phi),
\]
so that condition~(\ref{eq:supEabssqrd}) asserts integrability of
the function $\Phi\mapsto\lim_{\al}\sum_{A,B\in\al}|\Phi(A)\Phi(B)|$.
So, using again continuity of $\RR\to\RR:(x,y)\mapsto xy$ and dominated
convergence,
\[
  \begin{split}
    \langle f,g\rangle_\Sigma
      &= \int_{M(\scrX)}\lim_\al \sum_{A\in\al}
        f_\al(A)\Phi(A)\,\sum_{B\in\al}g_\al(B)\Phi(B)\dd\Pi_{\Sigma}(\Phi)\\
      &= \int_{M(\scrX)}\Bigl(\lim_\al\sum_{A\in\al}
        f_\al(A)\Phi(A)\Bigr)\,\Bigl(\lim_\be\sum_{B\in\be}
        g_\be(B)\Phi(B)\Bigr)\dd\Pi_{\Sigma}(\Phi)\\
      &= \int_{M(\scrX)}\Phi_f\,\Phi_g\,\dd\Pi_\Sigma(\Phi).
  \end{split}
\]
This implies that for any $f\in C(\scrX)$, $\langle f,f\rangle_\Sigma=0$,
if and only if, $\langle \Phi_f,\Phi_f \rangle_L=0$. Consequently the
pre-Hilbert spaces $C(\scrX)/Q_\Sigma$ and $L/Q_L$ are in
isometrically isomorphic correspondence, and so are their (unique)
Hilbert space completions.
\end{proof}

\subsection{Existence of weak Gaussian histogram limits}
\label{sub:weakgaussianlimits}

In subsection~\ref{sub:signedlimits}, we have seen
that for the existence of tight
histogram limits, condition~(\ref{property:P1signed}) 
(which is close to necessary, \cf\ 
remark~\ref{rem:almostnecessary}) says that the limiting
random $\Phi$ has a norm $\|\Phi\|_{1,\scrX}$ that is a
\emph{tight} real-valued random variable. We shall see in the
present subsection, that for the existence of weak Gaussian
histogram limits it is sufficient that $\|\Phi\|_{1,\scrX}$
is an \emph{integrable} random variable.

Let $\scrX$ be a compact Hausdorff space,
% We choose compact here, because 
% (1) the Dunford-Pettis-Grothendieck
% theorem is formulated for $L$- and $M$-spaces spaces, while we
% use spaces of Radon measures. They are equal in case of a compact
% $\scrX$ (see \cite{Kakutani41}) but I am not sure about the
% non-compact case.
% (2) Calculation and positive-definiteness of $\SIgma$ usually
% require that the space is bounded.
fix the topology on
$M(\scrX)$ to be the \emph{weak topology}, and choose a
directed set $\scrA$ of finite partitions in non-empty
Borel sets. Then $\scrX$, $\scrA$ and $M(\scrX)$ satisfy
the minimal conditions of definition~\ref{def:minimalconditions}.
Compactness of a subset of $M(\scrX)$ is still characterized by the
Dunford-Pettis-Grothendieck condition, but with $P$ replaced
by the positive measure $|\Phi|$: a subset $H$ of $M(\scrX)$ is
relatively compact in the weak topology, if and only if, for
some positive, bounded measure $Q\in M^+(\scrX)$,
\[
  \sup_{\Phi\in H}\,\bigl\| |\Phi|
    -|\Phi|\wedge L\,Q \bigr\|_{1,\scrX} \to 0,
\]
as $L\to\infty$. The proof of the existence theorem for
histogram limits that are Radon with respect to the weak
topology, does not differ substantially from that of
theorem~\ref{thm:WeakExistence}, so we omit explicit
statement.
\begin{theorem}
\label{thm:WeakExistenceSigned}
Let $\scrX$ be a compact Hausdorff space, consider $M(\scrX)$ with the
weak topology and choose a directed set of finite partitions
$\scrA$ in non-empty Borel sets that resolves $\scrX$. Let
$(\Pi_\al,\varphi_{*\,\al\be})$ be a coherent system
of Borel probability measures on the inverse system
$(M(\scrX_\al),\varphi_{*\,\al\be})$. There exists a
unique weakly-Radon probability measure $\Pi$ on $M(\scrX)$
projecting to $\Pi_\al$ for all $\al\in\scrA$, if and only if,
\begin{itemize}
\item[] there is a $Q\in M_b^+(\scrX)$ such that, for
  every $\ep,\delta>0$ there is a $L>0$ such that,
  \[
    \tag{P-weak-signed}\label{property:Pweaksigned}
    \Pi_\al\bigl(\{ \Phi_\al\in M(\scrX_\al):
      \bigl\||\Phi_\al|-|\Phi_\al|\wedge LQ_\al\bigr\|_{1,\scrX_\al}
      > \delta \}\bigr) < \ep,
  \]
  for all $\al\in\scrA$.
\end{itemize}
\end{theorem}
Given a Radon probability measure $\Pi$ on $M(\scrX)$,
the role of the mean measure $G$ of
definition~\ref{def:meanmeasure} as the dominating
measure, is taken over by the positive measure,
\[
  Q(B) = \int_{M(\scrX)} |\Phi|(B)\dd\Pi(\Phi),
\]
for $B\in\scrB$, and $Q(B)=0$
implies $\Pi(\{\Phi\in M(\scrX):|\Phi|(B)>0\})=0$
(as in lemma~\ref{lem:WeakSupport}). 
Proposition~\ref{prop:WeakSupport} can be adapted to the
signed case as well: $\{\Phi\in M(\scrX):|\Phi|\ll Q\}$ is closed
in $M(\scrX)$ and,
\[
  \supp_{\scrT_1}(\Pi) \subset \{\Phi\in M(\scrX):|\Phi|\ll Q\}.
\]
Moreover, if $\Phi\in M(\scrX)$ is such that for all
 $\al\in\scrA$, $|\Phi_\al|$ lies in
the support of $\Pi_\al$ in $M(\scrX_\al)$, then $\Phi$
lies in the weak support of $\Pi$.

Below, we apply theorem~\ref{thm:WeakExistenceSigned}
to Gaussian histogram systems. To prepare, it is noted
that for all $\al\in\scrA$ and all $A\in\al$,
\begin{equation}
  \label{eq:GaussianQalpha}
  Q_\al(A) =\int|\Phi_\al|(A)\dd\Pi_\al(\Phi_\al)
    = \sqrt{\ft{2}{\pi}\Sigma(A\times A)},
\end{equation}
% \begin{proposition}
% \label{prop:GaussianGalpha}
% Let $\Pi=\Pi_\Sigma$ be a (weak or tight) zero-mean
% Gaussian histogram limit, with covariance measure
% $\Sigma$. For all $\al\in\scrA$,
% \begin{equation}
% \label{eq:GaussianQalpha}
%   Q_\al(A) =\int|\Phi_\al|(A)\dd\Pi_\al(\Phi_\al)
%     = \sqrt{\ft{2}{\pi}\Sigma(A\times A)}.
% \end{equation}
% for all $A\in\al$.
% \end{proposition}
% \begin{proof}
% Let $\al\in\scrA$ and $A\in\al$ be given. Note that the
% component $\Phi_\al(A)$ of $\Phi_\al$ is distributed marginally
% according to the one-dimensional normal distribution
% with variance $\sigma^2=\Sigma(A\times A)$. Therefore,
% \[
%   Q_\al(A)
%     = \frac{1}{\sqrt{2\pi\sigma^2}}
%       \int_0^\infty 2x e^{-\ft12 x^2/\sigma^2} \dd x
%     = \sqrt{\ft{2\sigma^2}{\pi}}.
% \]
% \end{proof}
Clearly, the resulting set-functions $Q_\al:\sigma(\al)\to[0,\infty)$
are \emph{not the $\sigma(\al)$-restrictions} of the measure $Q$
(contrary to the cases of positive or probability histogram limits,
where all $G_\al$ are restrictions of the mean measure $G$). For
every $\al\in\scrA$ and all $A\in\al$, $|\Phi_\al|(A)\leq|\Phi|(A)$,
and hence $Q_\al(A)=\int|\Phi_\al|(A)\dd\Pi_\al(\Phi_\al)
\leq \int|\Phi|(A)\dd\Pi(\Phi)=Q(A)$.
If we assume that $\scrA$ resolves $\scrX$, the positive measures
$|\Phi_\al|$ (and the total-variational norms
$\|\Phi_\al\|_{1,\scrX_\al}=|\Phi_\al|(\scrX_\al)$)
increase to $|\Phi|$ (and the total-variational norm
$\| \Phi \|_{1,\scrX}=|\Phi|(\scrX)$).

\begin{corollary}
\label{cor:WeakExistenceGauss}
Let $\scrX$ be a compact Hausdorff space and let $\scrA$ be a
directed set of partitions generated by a basis, that resolves
$\scrX$. Let $\Sigma$ be a covariance measure on $\scrX$.
Consider $M(\scrX)$ with the weak topology and the centred Gaussian
histogram system $(\Pi_{\Sigma,\al},\varphi_{*\,\al\be})$.
If,
\begin{equation}
  \label{eq:Pweakgauss}
  Q(\scrX) = \sup_{\al\in\scrA}\sum_{A\in\al}
    \sqrt{\Sigma(A\times A)}<\infty,
\end{equation}
then then there exists a unique weakly Radon probability distribution
$\Pi_{\Sigma}$ on $M(\scrX)$ projecting to
$\Pi_{\Sigma,\al}$ for all $\al\in\scrA$.
\end{corollary}
\begin{proof}
Let $\ep,\delta>0$ be given and choose
$L=\sqrt{\pi}S/(\sqrt{2}\ep\delta)$, where
$S>0$ denotes any upper-bound for the left-hand side of
condition~(\ref{eq:Pweakgauss}). Due to the bound
(\ref{eq:boundforws}) and Markov's inequality, we have
for any $\al\in\scrA$,
\[
  \begin{split}
  \Pi_{\Sigma,\al}\bigl(\{ &\Phi_{\al}\in M(\scrX_{\al}):
    \bigl\||\Phi_{\al}|-|\Phi_{\al}|\wedge LQ_{\al}\bigr\|_{1,\scrX_{\al}}
    > \delta \}\bigr)\phantom{\biggm|}\\
  &\leq
  \Pi_{\Sigma,\al}\biggl(\Bigl\{ \Phi_{\al}\in M(\scrX_{\al}):
    \frac1{L} \sum_{A\in\al} \frac{\Phi_{\al}(A)^2}{Q_{\al}(A)}
    > \delta \Bigr\}\biggr)\\
  &\leq \frac1{L\delta}\int_{M(\scrX_{\al})} \sum_{A\in\al}
    \frac{\Phi_{\al}(A)^2}{Q_{\al}(A)} \dd\Pi_{\Sigma,\al}(\Phi_{\al})
  \leq \frac1{L\delta}
    \sum_{A\in\al} \sqrt{\ft{\pi}{2}\Sigma(A\times A)}<\ep,
  \end{split}
\]
where we use that, for every $\al\in\scrA$ and all $A\in\al$,
$\int\Phi_{\al}(A)^2\dd\Pi_{\Sigma,\al}(\Phi_{\al})=\Sigma(A\times A)$,
and (\ref{eq:GaussianQalpha}).
\end{proof}
Two remarks are in order: firstly, we relate
condition~(\ref{eq:Pweakgauss}) for the existence of weak
Gaussian limits, to condition~(\ref{eq:supEabssqrd}) for
existence of tight Gaussian limits, by the Cauchy-Schwartz
inequality:
\[
  \begin{split}
  \sum_{A,B\in\al} \int_{M(\scrX_\al)} &
     \bigl| \Phi_\al(A)\Phi_\al(B) \bigr|\dd\Pi_{\Sigma,\al}(\Phi_\al)\\
    &\leq \sum_{A,B\in\al} \sqrt{\Sigma(A\times A)} \sqrt{\Sigma(B\times B)}
    = \Bigl(\sum_{A\in\al} \sqrt{\Sigma(A\times A)}\Bigr)^2,
  \end{split}
\]
showing that~(\ref{eq:Pweakgauss}) implies~(\ref{eq:supEabssqrd}).
Second, we note that corollary~\ref{cor:TVExistence} stays valid
in the signed case, so if condition~(\ref{eq:Pweakgauss}) holds,
there exists a unique Radon probability measure $\Pi_\Sigma$ on
$M(\scrX)$ with the \emph{total-variational} topology, projecting
to $\Pi_{\Sigma,\al}$ for all $\al\in\scrA$.

\begin{example}
\label{ex:constgauss}
Let $\scrX$ be a compact subset of $\RR^d$. Consider
example~\ref{ex:kernelfunctions} with a kernel function
that is constant: $k(x,y)=c$ for some $c>0$. If we let
$\scrA$ consist of partitions $\al$ with the property
that for all $A,A'\in\al$, there exists an translation
vector $x$ in $\RR^d$ such that the Lebesgue measure of
$(A\setminus(A'-x))\cup((A'-x)\setminus A)$ is zero. Then,
for every $\al$, $\Sigma_\al$ is the
${|\al|}\times{|\al|}$-matrix with all entries equal to,
\[
  \Sigma_{\al,ij} = \int_{A_i\times A_j} c \dd x \dd y
    = c\mu(A_i)\mu(A_j) = c \mu_\al^2,
\]
(where $\mu_\al=\mu(A)$ for any $A\in\al$). Clearly, the
corresponding covariance measure $\Sigma$ gives rise to
positive semi-definite covariance matrices $\Sigma_\al$,
with random histogram components that are highly
dependent: in fact, the linear space of all
$\mu_\al\in M(\scrX_\al)$ with components that sum to
zero forms the kernel of $\Sigma_\al$,
\[
    \text{Ker}(\Sigma_\al) = \biggr\{\mu_\al\in M(\scrX_\al)\,:\,
    \sum_{i\in I(\al)}\mu_{\al,i}=0\biggr\},
\]
and a centred multivariate normal distribution is
supported on the range of its covariance matrix.
That means that $\Phi_\al$ lies on the diagonal
of $M(\scrX_\al)$ with probability one,
\begin{equation}
  \label{eq:randomlebesgue}
  \Pi_{\Sigma,\al}\bigl(\bigr\{ \Phi_\al\in M(\scrX_\al)\,:\,
    \Phi_{\al,1}=\ldots=\Phi_{\al,{|\al|}}\bigr\}\bigr)=1.
\end{equation}
Note that,
\begin{equation}
\label{eq:diagonalsqrtsummable}
  \sup_{\al\in\scrA}\sum_{A\in\al} \sqrt{\Sigma(A\times A)}
    = \sup_{\al\in\scrA}\sum_{A\in\al}\sqrt{c}\,\mu(A)
    = \sqrt{c}\,\mu(\scrX)<\infty,
\end{equation}
so according to corollary~\ref{cor:WeakExistenceGauss}, there
exists a weakly-Radon probability measure $\Pi_{\Sigma}$ on
$M(\scrX)$ projecting to $\Pi_{\Sigma,\al}$ for all
$\al\in\scrA$.
\end{example}
Additionally, a moment's thought shows that the above example
serves as a bound for a host of examples based on the reproducing
kernels of example~\ref{ex:kernelfunctions}.
\begin{corollary}
\label{cor:bndkernel}
Let $\scrX$ be a compact Hausdorff space and let $\scrA$ be a
directed set of partitions generated by a basis, that resolves
$\scrX$. Let $\Sigma$ be a covariance measure based on a
bounded kernel function $k:\scrX\times\scrX\to\RR$. Then
there exists a unique weakly-Radon probability distribution
$\Pi_{\Sigma}$ on $M(\scrX)$ projecting to
$\Pi_{\Sigma,\al}$ for all $\al\in\scrA$.
\end{corollary}
\begin{proof}
Assume first, that $\scrA$ is as in example~\ref{ex:constgauss}.
Since for some $c>0$ and all $x,y\in\scrX$, $|k(x,y)|\leq c$,
for all $\al$ and all $A\in\al$, $\Sigma(A\times A)\leq c\mu(A)^2$.
Summability then follows as in
inequality~(\ref{eq:diagonalsqrtsummable}), proving the
existence of a weakly-Radon histogram limit $\Pi_\Sigma$.
For any Borel-measurable partition $\al$ of $\scrX$, the
weakly continuous mapping $\varphi_{*\,\al}:M(\scrX)\to
M(\scrX_\al):\mu\mapsto(\mu(A_1),\ldots,\mu(A_{|\al|}))$
induces the Gaussian histogram distribution
$\Pi_{\Sigma,\al}= N(0,\Sigma_\al)$ on $M(\scrX_\al)$.
Uniqueness of the limit proves the assertion.
\end{proof}

Over recent years there has been considerable interest in
the so-called Gaussian Free Field (see, \eg,
\cite{Sheffield07,Werner20}); below we use Green's functions for
the harmonic operators in $d=1,2,\ldots$ as covariance
kernels to define Gaussian histogram systems, in the closure
$\scrX$ of a non-empty, bounded, open subset of $\RR^d$.
\begin{example}
We consider the existence question for the Gaussian Free
Field first in $d=1$: the Green's function
$G_1:\scrX\times\scrX\to\RR$ is of the form,
\[
  G_1(x,y) = -|x-y|+f(x,y),
\]
where $x\mapsto f(x,y)$ is harmonic in $x$
% (\ie\ lies in the kernel of $\laplace_x$)
for every $y$. Define,
\[
  \Sigma_1(A\times B) = \int_{A\times B} G_1(x,y) \dd x \dd y,
\]
(where we choose $f$ such that $\Sigma_1$ is symmetric and
positive-definite). Choose a directed set $\scrA$ of
partitions generated by a basis
that resolves $\scrX$. Based on corollary~\ref{cor:bndkernel},
we see immediately that the associated centred Gaussian
histogram system with histogram distributions
$\Pi_{\Sigma_1,\al}$ has a limit $\Pi_{\Sigma_1}$, that is
a weakly Radon probability measure on $M(\scrX)$. Then
$Q(B)=\int|\Phi|(B)\dd\Pi_\Sigma(\Phi)$, ($B\in\scrB$), is a
multiple of Lebesgue measure, due to translation
invariance, and,
\[
  \Pi_{\Sigma_1}\bigl(\{\Phi\in M(\scrX):\Phi\ll Q\}\bigr)=1,
\]
implying that the random element $\Phi\sim\Pi_\Sigma$ is of
the form,
\[
  \Phi(B) = \int_B \phi(x) \dd x,
\]
for all $B\in\scrB$, where $\phi$ is a random Radon-Nikodym
density function in $L^1(\scrX)$, the Banach space of
Lebesgue integrable functions on $\scrX$.

For $d=2,3,\ldots$, the situation changes drastically: Green's
functions $G_d:\scrX\times\scrX\to\RR$ are unbounded and
display singular behaviour in neighbourhoods of the diagonal,
namely $G_2(x,y)=-\log|x-y|$ and $G_d(x,y)=-|x-y|^{d-2}$ for
$d\geq3$.
\begin{figure}[ht]
\begin{center}
%\vspace*{0.50em}
%\hspace*{-6ex}
\parbox{3.75cm}{
  %\begin{lpic}{./graph/G2-8x8(0.3)}
  \begin{lpic}{./G2-8x8(0.3)}
    % \lbl[t]{24,22;{\footnotesize $0$}}
    % \lbl[t]{82,22;{\footnotesize $n=16$}}
    % \lbl[t]{138,22;{\footnotesize $1$}}  
  \end{lpic}
}
\parbox{3.75cm}{
  %\begin{lpic}{./graph/G2-16x16(0.3)}
  \begin{lpic}{./G2-16x16(0.3)}
    % \lbl[t]{24,22;{\footnotesize $0$}}
    % \lbl[t]{82,22;{\footnotesize $n=32$}}
    % \lbl[t]{138,22;{\footnotesize $1$}}    
  \end{lpic}
}
\parbox{3.75cm}{
  %\begin{lpic}{./graph/G2-32x32(0.3)}
  \begin{lpic}{./G2-32x32(0.3)}
    % \lbl[t]{24,22;{\footnotesize $0$}}
    % \lbl[t]{82,22;{\footnotesize $n=64$}}
    % \lbl[t]{138,22;{\footnotesize $1$}}    
  \end{lpic}
}
\parbox{3.75cm}{
  %\begin{lpic}{./graph/G2-64x64(0.3)}
  \begin{lpic}{./G2-64x64(0.3)}
    % \lbl[t]{24,22;{\footnotesize $0$}}
    % \lbl[t]{82,22;{\footnotesize $n=128$}}
    % \lbl[t]{138,22;{\footnotesize $1$}}    
  \end{lpic}
}\\%[-1.8em]
\parbox{13cm}{
  \caption{
    \label{fig:coarsenedhistograms}
    \footnotesize{A sample from a random histogram on a 64x64
    partitioned square patch of two-dimensional Euclidean space-time
    with the Green's function for the Laplacian to
    define the covariance measure, and its
    32x32, 16x16 and 8x8 coarsened histograms. Coherence of the histogram
    system says that the distributions of the random 8x8, 16x16
    and 32x32 histograms must equal the distributions implied
    by coarsening of the random 64x64 histogram. The histogram
    limit is the random object obtained by infinite refinement,
    to the infinite right of these four histograms.}
  }
}
\end{center}
%\vspace*{-1.5em}
\end{figure}
\begin{figure}[ht]
\begin{center}
%\vspace*{-1.0em}
%\hspace*{-6ex}
\parbox{3.75cm}{
  %\begin{lpic}{./graph/G2-64x64-2(0.3)}
  \begin{lpic}{./G2-64x64-2(0.3)}
    % \lbl[t]{24,22;{\footnotesize $0$}}
    % \lbl[t]{82,22;{\footnotesize $n=16$}}
    % \lbl[t]{138,22;{\footnotesize $1$}}  
  \end{lpic}
}
\parbox{3.75cm}{
  %\begin{lpic}{./graph/G3-64x64(0.3)}
  \begin{lpic}{./G3-64x64(0.3)}
    % \lbl[t]{24,22;{\footnotesize $0$}}
    % \lbl[t]{82,22;{\footnotesize $n=32$}}
    % \lbl[t]{138,22;{\footnotesize $1$}}    
  \end{lpic}
}
\parbox{3.75cm}{
  %\begin{lpic}{./graph/G4-64x64(0.3)}
  \begin{lpic}{./G4-64x64(0.3)}
    % \lbl[t]{24,22;{\footnotesize $0$}}
    % \lbl[t]{82,22;{\footnotesize $n=64$}}
    % \lbl[t]{138,22;{\footnotesize $1$}}    
  \end{lpic}
}
\parbox{3.75cm}{
  %\begin{lpic}{./graph/4D-yukawa-sample-64x64-delta=0.1-noframe(0.3)}
  \begin{lpic}{./4D-yukawa-sample-64x64-delta=0.1-noframe(0.3)}
    % \lbl[t]{24,22;{\footnotesize $0$}}
    % \lbl[t]{82,22;{\footnotesize $n=128$}}
    % \lbl[t]{138,22;{\footnotesize $1$}}    
  \end{lpic}
}\\%[-1.8em]
\parbox{13cm}{
  \caption{
    \label{fig:gaussianhistograms}
    \footnotesize{Samples from Gaussian random histograms on a 64x64
    partitioned
    square slice of Euclidean space-time, with the Green's function
    for the Laplacian to define the covariance measure, in two
    dimensions; in three
    dimensions; and in four dimensions; alongside a sample with the
    Yukawa potential of a massive scalar boson field in four
    dimensions.}
  }
}
\end{center}
\end{figure}
To apply corollary~\ref{cor:bndkernel} 
we modify $G_d$ for small length scales to regularize the
singular behaviour near the diagonal (\eg, for some small
$\ep>0$, replace $|x-y|^{d-2}$ by $|(|x-y|+\ep)^{d-2}|$,
which replaces the pole for $x=y$ with an upper bound
$1/\ep^{d-2}$). The modified $G_d$ are bounded kernel
functions, and corollary~\ref{cor:bndkernel} guarantees
that for every $\ep>0$, there exists a weak histogram limit
$\Pi_{\Sigma_d,\ep}$. One may hope that the limit
for $\Pi_{\Sigma_d,\ep}$ as $\ep\to0$ (which may exist only
as a tightly Radon probability measure) describes the
so-called \emph{Gaussian Free Field in
$d$ Euclidean dimensions}. In light of earlier explorations
(see, for example, \cite{Werner20} for a detailed
discussion of (mostly) the $d=2$ case), it is also possible
that the $\ep\to0$ limit exists only if we embed the space $M(\scrX)$
of Radon measures on $\scrX$, in spaces of distributions on
$\scrX$. The limiting probability distribution
$\Pi_{\Sigma_d}$ would then describe a \emph{random generalized
function} of the type discussed in (\cite{Gelfand64})
and \IIIthree{IX}{6}{10}.
\end{example}

\subsection{Completely random Gaussian histogram limits}
\label{sub:CRgaussianlimits}

The class of (centred) Gaussian histogram limits has a non-empty
intersection with the class of completely random measures,
characterized by covariance measures that place all
mass on the diagonal. Completely random Gaussian histogram
limits exist in the fixed or random atomic phase.
\begin{definition}
Let $\Sigma$ be a covariance measure on $\scrX\times\scrX$.
We say that $\Sigma$ is \emph{diagonal}, if for all $A,B\in\scrB$,
$\Sigma(A \times B) = \Sigma( A\cap B \times A \cap B)$.
\end{definition}
Note that with a diagonal $\Sigma$, the set-function
$\tau:\scrB\to[0,\infty)$ defined by $\tau(B)=\Sigma(B\times B)$ for all
$B\in\scrB$, is a positive Radon measure on $\scrX$, and the Hilbert
space $L^2(\Sigma)$ is isometrically isomorphic to
$L^2(\tau)$, the usual space of $\tau$-square-integrable
functions on $\scrX$.

A diagonal covariance measure leaves histogram components
independent, and leads to a completely random limit, which
is of a fixed- or random-atomic nature also in the case
of a \emph{signed} completely random measure.
If a Gaussian histogram system with a diagonal covariance measure
has a tight limit and $\Sigma$ distributes it mass in a uniformly
asymptotically negligible way (see \cite{Feller91}, section~XVII.7),
infinite divisibility of the distribution of the random variable
$\|\Phi\|_{1,\scrX}$ is implied.
\begin{corollary}
\label{cor:diaggauss}
Let $\scrX$ be a compact Polish space and let $\scrA$ be a
directed set of partitions that resolves $\scrX$ and is generated
by a basis. Assume that $\Sigma$ is a diagonal covariance measure.
Then the centred Gaussian histogram system
$(\Pi_{\Sigma,\al},\varphi_{*\,\al\be})$ has a unique tightly-Radon
probability distribution $\Pi_{\Sigma}$ on $M(\scrX)$ projecting
to $\Pi_{\Sigma,\al}$ for all $\al\in\scrA$. If, in addition,
\begin{equation}
  \label{eq:UnifAsympNegl}
  \lim_{\al\in\scrA}\,\max_{A\in\al}\,\tau(A) = 0,
\end{equation}
then the total-variational norm $\|\Phi\|_{1,\scrX}$ has a
probability distribution that is infinitely divisible.
\end{corollary}
\begin{proof}
For a diagonal covariance measure $\Sigma$, any $\al$
and any $t\in\RR$ the restrictions of cumulant measures
$\lambda_t$ of definition~\ref{def:cumulant}
to the $\sigma$-algebras $\sigma(\al)$ are given by,
\[
  \begin{split}
  \lambda_t(B) &= \log\int_{M(\scrX_\al)}
      e^{t\Phi_\al(B)} \dd\Pi_{\Sigma,\al}(\Phi_\al)
    = \log\int_{M(\scrX_\al)} \prod_{A\subset B}
      e^{t\Phi_\al(A)} \dd\Pi_{\Sigma,\al}(\Phi_\al)\\
    &= \log\prod_{A\subset B}\int_\RR e^{t\Phi_\al(A)}\dd\Pi_{\Sigma,\al}(\Phi_\al(A))
    = \log\prod_{A\subset B} e^{\ft12 t^2\tau(A)}
    = \frac{t^2}2 \sum_{A\subset B} \tau(A) = \frac{t^2}2\tau(B),
  \end{split}
\]
for any $B\in\sigma(\al)$. By the Carath\'eodory extension,
all cumulants are therefore finite positive measures, and hence
(\cite{Hellmund09}, theorem~4.4),
there exists a tight completely random limit described by
a marked Poisson process (\citep{daley07}, Ch.~9-10)
(with marks in $\RR$ rather than $(0,\infty)$), of the form
(\ref{eq:KingmanCRM}).

With a diagonal covariance measure $\Sigma$, the components
of $\Phi_\al$ are independent random variables for every
$\al\in\scrA$. Therefore the norms of the random histograms
$\Phi_\al$,
\[
  \bigl\|\Phi_\al\bigr\|_{1,\al}
    = \sum_{A\in\al}\bigl|\Phi_\al(A)\bigr|,
\]
are sums of \emph{independent terms} in a triangular array,
and \emph{uniform asymptotic negligibility}, as assumed in
(\ref{eq:UnifAsympNegl}), is sufficient for
tight convergence to an infinitely divisible limiting
probability distribution (\cite{Feller91}, section~XVII.7).
Since the total-variational norm $\|\Phi\|_{1,\scrX}$ is
the monotone limit of the norms $\|\Phi_\al\|_{1,\al}$
(proposition~\ref{prop:alphatotalvariation}), its
probability distribution is tight and infinitely divisible.
\end{proof}

\begin{example}
\label{ex:BrownianMotion}
Let $\scrX=[0,1]$, let $\scrA$ consist of partitions in
half-open intervals, generated by a basis and collectively
fine enough to resolve $\scrX$. Consider a centred
Gaussian histogram system with diagonal covariance measure
$\Sigma$ defined by choosing $\tau$ equal to Lebesgue measure.
With sets $A_{\al,i}$ of the form $(s_{\al,i},s_{\al,i+1}]$ (and
the singleton $\{0\}$, for which we set $\Phi_\al(\{0\})=0$
for all $\al\in\scrA$), the histogram system,
\[
  \bigl(\Phi_\al(A_{\al,1}),\ldots,\Phi_\al(A_{\al,{|\al|}})\bigr)
    \sim \prod_{A\in\al} N\bigl(0,(s_{\al,i+1}-s_{\al,i})\bigr),
\]
describes the independent, normally distributed increments
of Brownian motion started from $B(0)=0$, so the (random)
Stieltjes function $B:[0,1]\to\RR$ for the measure $\Phi$,
\[
  B(t)=\inf_{\al\in\scrA}\sum\{\Phi_\al(A_{\al,i}):1
    \leq i\leq{|\al|}, s_{\al,i}< t\},
\]
is a version of the sample path of
Brownian motion on $[0,1]$. Since $\scrA$ resolves $\scrX$,
the Lebesgue measures of all intervals
$(s_{\al,i+1},s_{\al,i}]$ goes to zero as $\al$ increases
in $\scrA$, and condition~(\ref{eq:UnifAsympNegl}) is
satisfied. (By extension, if we
replace normal distributions by stable distributions
in this construction appropriately, infinite divisibility
preserves coherence and existence of $\Phi$, \cf\
theorem~\ref{thm:vagueexistence}, implies random
Stieltjes functions corresponding to right-continuous
versions of L\'evy sample paths on $[0,1]$.)
\end{example}
Weak histogram limits with diagonal covariance measures
display a limitation similar to that of
theorem~\ref{thm:weakdirichlet}. To appreciate the
problem, note that for a diagonal covariance measure
with positive $\tau$ dominated by Lebesgue measure,
condition~(\ref{eq:Pweakgauss}) cannot be satisfied.
Purely atomic measures $\tau$, however, lead to Gaussian
histogram systems with weak limits.

Corollary~\ref{cor:diaggauss} and example~\ref{ex:constgauss}
form two extremes: in diagonal cases Gaussian histogram limits
manifest in the fixed- or random-atomic phase, while covariance
measures that spread their mass more homogeneously over
$\scrX\times\scrX$, dependence introduces a degree of
smoothness, a situation that we have seen in its most extreme
form in (the highly non-diagonal) example~\ref{ex:constgauss}.
Gaussian histogram limits for other covariance measures $\Sigma$
are somewhere in between: depending on the degree to which
$\Sigma$-mass is located away from the diagonal,
corresponding to the degree of dependence between Gaussian
histogram components, the histogram limit
may manifest in close-to-atomic (\ie\ highly concentrated) or
in smooth/closer-to-constant form.

To demonstrate the explanatory value of the phase structure
of Gaussian histogram limits described above, the last example,
which is analysed more comprehensively in forthcoming work
\citep{Kleijn25b}, suggests applicability of Gaussian histogram
limits in (Euclidean) quantum field theory \citep{Zuber12}.
\begin{example}
Let $d\geq1 $ be given and consider
$\scrX=\{p\in\RR^d:\|p\|\leq\Lambda\}$ for some constant
$\Lambda>0$, and $\Sigma$ diagonal with,
\[
  \tau(B)=\int_B\frac{1}{p^2}\dd p,
\]
for all $B\in\scrB$. The space $\scrX$ plays the role of
$d$-dimensional Euclidean `momentum space' (and the constant
$\Lambda$ that makes $\scrX$ compact, is known as the
\emph{UV-cutoff scale} in physics). The kernel defining
$\Sigma$ is interpreted as the (unregularized, Euclidean)
`propagator' of the massless scalar field (roughly, the
Green's function for the Laplace operator, which is
represented by the convolution kernel $-p^2\delta_d(p-q)$
in momentum space). The diagonal Gaussian histogram limit
exists, \cf\ corollary~\ref{cor:diaggauss}.

In anticipation of the analysis in \citep{Kleijn25b}, we point
out the following consequence of the phase of this Gaussian
histogram limit. In this case, \emph{`quantization'}, a
description of the field in terms of particles,
emerges as a consequence of complete randomness: the
Gaussian histogram limit is in the \emph{random-atomic
phase} and manifests as a random sum of discrete point
masses in Euclidean momentum space. Such configurations have
an immediate physical interpretation, as states describing
(off-shell) particles, point-like quanta of momentum.
It is noted that the \emph{emergence of quantization} is not
a feature of (second-quantized) quantum field theory: in the
physical theory of quantum fields, particles are axiomatic and
introduced by hand with the formal introduction of a
Fock space to describe quantum states of the field (see
the classical work \cite{Zuber12}, p.~106). 
% \begin{quotation}
%   \noindent\emph{At first sight a field theory does not seem to
%   have an interpretation of microobjects identified as particles.
%   The long story of the wave-particle duality in the theory of
%   light is a testimony of the confusion that may arise and
%   still remains in certain aspects up to this day. [...] The
%   structure of the Hilbert space of states constructed for
%   free fields---the Fock space---will reflect this aspect.}
% \end{quotation}
\end{example}

%%%%%%%%%%%%%%%%%%%%%%%%%%%%%%%%%%%%%%%%%%%%%%%%%%%%%%%%%%%%%%%%%%%%%%%%%%%%%%
% \appendix
% 
% \section{Appendix}
% \label{app:}
%%%%%%%%%%%%%%%%%%%%%%%%%%%%%%%%%%%%%%%%%%%%%%%%%%%%%%%%%%%%%%%%%%%%%%%%%%%%%%
\end{document}